\documentclass[a4paper,reqno]{amsart}

\usepackage{mathdots}
\usepackage{amssymb}
\usepackage{amstext}
\usepackage{amsmath}
\usepackage{amscd}
\usepackage{amsthm}
\usepackage{amsfonts}
\usepackage{enumerate}
\usepackage{graphicx}
\usepackage{latexsym}
\usepackage[all,color,dvips]{xy}
\xyoption{all}

\usepackage{pstricks}
\usepackage{lscape}
\usepackage{comment}

\newtheorem{theorem}{Theorem}[section]

\newtheorem{corollary}[theorem]{Corollary}
\newtheorem{lemma}[theorem]{Lemma}
\newtheorem{proposition}[theorem]{Proposition}
\newtheorem{question}[theorem]{Question}

\theoremstyle{definition}
\newtheorem{definition}[theorem]{Definition}
\newtheorem{remark}[theorem]{Remark}
\newtheorem{example}[theorem]{Example}

\newcommand{\arr}[1]{\overset{#1}{\rightarrow}}
\newcommand{\LM}[1]{\mathop{L_{#1}}}
\newcommand{\RM}[1]{\mathop{R_{#1}}}
\newcommand{\com}[1]{\operatorname{com}\nolimits_{#1}}

\newcommand{\tildetensor}{\tilde\otimes}
\newcommand{\C}{{\mathcal C}}
\renewcommand{\P}{{\mathcal P}}
\newcommand{\J}{{\mathcal J}}
\newcommand{\DD}{{\mathcal D}}

\newcommand{\N}{\mathbf{N}}
\newcommand{\Z}{\mathbf{Z}}
\newcommand{\R}{\mathbf{R}}
\newcommand{\per}{\operatorname{per}\nolimits}
\newcommand{\add}{\operatorname{add}\nolimits}
\newcommand{\proj}{\operatorname{proj}\nolimits}

\newcommand{\gl}{\operatorname{gl.dim}\nolimits}

\newcommand{\id}{\operatorname{id}\nolimits}
\newcommand{\injdim}{\operatorname{inj.dim}\nolimits}

\newcommand{\Hom}{\operatorname{Hom}\nolimits}
\newcommand{\End}{\operatorname{End}\nolimits}
\newcommand{\Ext}{\operatorname{Ext}\nolimits}

\newcommand{\RHom}{\mathbf{R}\strut\kern-.2em\operatorname{Hom}\nolimits}

\renewcommand{\mod}{\operatorname{mod}\nolimits}

\newcommand{\cut}{\ar@{.}}
\newcommand{\lat}{\ar@[|(4)]@{-}}
\newcommand{\point}[1]{\bullet}

\begin{document}
\title{Selfinjective quivers with potential and 2-representation-finite algebras}
\author{Martin Herschend and Osamu Iyama}
\address{O. Iyama: Graduate School of Mathematics, Nagoya University, Chikusa-ku, Nagoya,
464-8602 Japan}
\email{iyama@math.nagoya-u.ac.jp}
\address{M. Herschend: Graduate School of Mathematics, Nagoya University, Chikusa-ku, Nagoya,
464-8602 Japan}
\email{martin.herschend@gmail.com}

\begin{abstract}
We study quivers with potential (QPs) whose Jacobian algebras are finite dimensional selfinjective.
They are an analogue of the `good QPs' studied by Bocklandt whose Jacobian algebras are 3-Calabi-Yau.
We show that 2-representation-finite algebras are truncated Jacobian algebras of
selfinjective QPs, which are factor algebras of Jacobian algebras by certain sets of arrows called cuts.
We show that selfinjectivity of QPs is preserved under successive mutation with respect to orbits of the Nakayama permutation.
We give a sufficient condition for all truncated Jacobian algebras of a fixed QP to be derived equivalent.
We introduce planar QPs which provide us with a rich source of selfinjective QPs.
\end{abstract}
\maketitle
\tableofcontents

\section{Introduction}
Quivers with potential (QPs) play important roles in the theory of cluster algebras/categories \cite{BIRS,DWZ1,DWZ2,KY}
and also appear in Seiberg duality in physics \cite{BD,S}.
They give rise to Ginzburg DG algebras \cite{G}, which enjoy the 3-Calabi-Yau property and play a crucial role in the theory of generalized cluster categories \cite{Am,K2}.
The zero-th homology $\P(Q,W)$ of the Ginzburg DG algebra associated with a QP $(Q,W)$ is called the \emph{Jacobian algebra}.
Important classes of algebras appearing in representation theory are known to be Jacobian algebras of certain QPs. 
For example graded 3-Calabi-Yau algebras are Jacobian algebras of certain QPs \cite{Bo}.
Also every algebra $A$ of global dimension two gives a QPs $(Q_A,W_A)$ whose Jacobian algebra is the 3-preprojective algebras $\Pi_3(A)$ of $A$ \cite{IO2,K2}.
In particular cluster tilted algebras are Jacobian algebras of certain QPs \cite{BIRS}.

Jacobian algebras of QPs are Iwanaga-Gorenstein of dimension at most one if they are finite dimensional \cite{BIRS}.
A main subject in this paper is \emph{selfinjective} QPs, which are QPs whose Jacobian algebras are finite dimensional selfinjective algebras. 
It is known that a finite dimensional algebra $A$ of global dimension two has a cluster tilting modules if and only if
the 3-preprojective algebra $\Pi_3(A)$ is selfinjective \cite{IO2}.
In this case $A$ is called \emph{2-representation-finite}. 2-representation-finite algebras play a central role in 3-dimensional Auslander-Reiten theory \cite{HI,I1,I2,I3,I4,IO1,IO2}.
We shall give a structure theorem of 2-representation-finite algebras in terms of selfinjective QPs and their truncated Jacobian algebras,
which is another main subject in this paper.

For a QP $(Q,W)$, we call a set $C$ of arrows of $Q$ a \emph{cut}
if each cycle of $Q$ appearing in $W$ contains precisely one
arrow in $C$. In this case the factor algebra of the Jacobian
algebra by the ideal generated by arrows in $C$ is called
the \emph{truncated Jacobian algebra}.
For example algebras $A$ of global dimension two are truncated Jacobian algebras of QPs $(Q_A,W_A)$. 
This gives us a structure theory of tilted algebras in terms of cluster tilted algebras \cite{BMR,BFPPT}.
In \cite{IO1} a family of 2-representation-finite algebras are constructed as truncated Jacobian algebras of certain QPs.
Our first main result (Theorem \ref{QP}) asserts that 2-representation-finite algebras
are exactly truncated Jacobian algebras of selfinjective QPs.
This structure theorem provides us with a rich source of 2-representation-finite algebras.

Notice that selfinjective QPs are an analogue of the `good QPs' studied by Bocklandt \cite{Bo} whose Jacobian algebras are 3-Calabi-Yau.
In a forthcoming paper, the authors will study finite dimensional algebras called `2-representation-infinite' by using good QPs.

We study the relationship between selfinjectivity of QPs and
mutation of QPs introduced by Derksen-Weyman-Zelevinsky \cite{DWZ1}.
Our second main result (Theorem \ref{selfinjective QP}) asserts that selfinjectivity of QPs
is preserved under mutation along orbits of the Nakayama permutation.
This provides us a systematic method to construct selfinjective QPs in Sections \ref{Examples of selfinjective QPs} and \ref{planar}.

We also study derived equivalence between truncated Jacobian
algebras of a fixed QP. Although all truncated Jacobian algebras of
a fixed QP are cluster equivalent in the sense of Amiot \cite{Am2}, they are
not necessarily derived equivalent. So it is natural to ask when
all truncated Jacobian algebras of QPs are derived equivalent (e.g. \cite{AO}).
Our third main result (Theorem \ref{transitive}) provides a sufficient condition.
This is based on a combinatorial consideration of Galois coverings $\widetilde{Q}$ of quivers $Q$,
which generalizes Riedtmann's translation quivers $\Z Q$ associated with quivers $Q$. Key results are a bijection (\emph{cut-slice correspondence}) between certain sets of arrows in $Q$ and certain
full subquivers of $\widetilde{Q}$ called \emph{slices} (Theorem \ref{cut slice})
and transitivity of a certain combinatorial operation called \emph{cut-mutation} (Theorem \ref{transitive mutation cut}).
This is relevant for us as cut-mutation is a combinatorial interpretation of 2-APR tilting \cite{IO1} (Theorem \ref{2-APR tilting}).

In Section \ref{canvas}, we introduce a CW complex associated with every QP which we call the \emph{canvas}. 
We call a QP \emph{simply connected} if the canvas is simply connected. 
Simply connectedness is useful to show that all truncated Jacobian algebras are derived equivalent.
In Section \ref{planar} we introduce \emph{planar QPs}, which are analogue of QPs associated with dimer models (e.g. \cite{Br,D,IU}) and provide a rich source of simply connected QPs.
We shall give three mutation-equivalence classes of planar QPs which are selfinjective.
It is an open question whether they are all selfinjective planar QPs.

\medskip\noindent
{\bf Conventions }
Throughout this paper we denote by $K$ an algebracially closed field,
and by $D=\Hom_K(-,K)$ the $K$-dual.
For a quiver $Q$, we denote by $Q_0$ the set of vertices and by $Q_1$ the set of arrows.
We denote by $a:s(a)\to e(a)$ the start and end vertices.
The composition $ab$ of arrows $a$ and $b$ means first $a$ and then $b$.

\medskip
\noindent{\bf Acknowledgement }
The first author is grateful to JSPS for funding his stay at Nagoya University.
The second author was supported by JSPS Grant-in-Aid for Scientific Research 21740010 and 21340003.

\section{Preliminaries}

\subsection{Quivers with potential}
Let $Q$ be a quiver. We denote by $\widehat{KQ}$ the \emph{complete path algebra}
and by $J_{\widehat{KQ}}$ the Jacobson radical of $\widehat{KQ}$ \cite{BIRS}.
Then $\widehat{KQ}$ is a topological algebra with $J_{\widehat{KQ}}$-adic topology.
We denote
\[
\com Q := \overline{\left[\widehat{KQ}, \widehat{KQ}\right]}\subset \widehat{KQ},
\]
where $\overline{(\ )}$ denotes the closure.
Then $\widehat{KQ}/\com Q$ has a topological basis consisting of cyclic paths considered up to cyclic permutation. Thus there is a unique continuous linear map
\[
\sigma : \widehat{KQ} /\com Q \rightarrow \widehat{KQ}
\]
induced by $a_1\cdots a_n \mapsto \sum_{i} a_i\cdots a_n a_1 \cdots a_{i-1}$ for all cycles $a_1\cdots a_n$. For each $a\in Q_1$ define the continuous linear maps
\[
\LM{a^{-1}}, \RM{a^{-1}} : J_{\widehat{KQ}} \rightarrow \widehat{KQ},
\]
by $\LM{a^{-1}}(ap) = p$, $\RM{a^{-1}} (pa) = p$ and $\LM{a^{-1}}p' = 0 = \RM{a^{-1}}p''$ for all paths $p,p',p''$ such that $p'$ and $p''$ do not start, respectively end with $a$.

Composing these maps we obtain the cyclic derivative
\[
\partial_a = \LM{a^{-1}} \circ \sigma = \RM{a^{-1}} \circ \sigma:  J_{\widehat{KQ}} /\com Q \rightarrow \widehat{KQ}.
\]
We also define 
\[
\partial_{(a,b)} = \LM{a^{-1}} \circ \RM{b^{-1}} \circ \sigma = \RM{b^{-1}} \circ \LM{a^{-1}} \circ \sigma :
J_{\widehat{KQ}}^2 /(J_{\widehat{KQ}}^2\cap\com Q) \rightarrow \widehat{KQ}.
\]
A \emph{potential} is an element $W \in J_{\widehat{KQ}}^2 /(J_{\widehat{KQ}}^2\cap\com Q)$. We call $(Q,W)$ a \emph{quiver with potential} (\emph{QP}) and define the \emph{Jacobian ideal} and the \emph{Jacobian algebra} by
\begin{eqnarray*}
\J(Q,W) &:=& \overline{\langle \partial_a W \;|\; a \in Q_1 \rangle},\\
\P(Q,W) &:=& \widehat{KQ}/\J(Q,W)
\end{eqnarray*}
respectively.

\subsection{3-preprojective algebras and generalized cluster categories}\label{preliminary2}

The following notion will play an important role in this paper.

\begin{definition}\label{K2}
Let $A$ be a finite dimensional $K$-algebra of global dimension at most two.
\begin{itemize}
\item[(a)]\cite{K2} For a presentation $A = \widehat{KQ}\left/ \overline{\langle r_1, \ldots, r_l \rangle}\right.$ by a quiver $Q$ and a minimal set $\{r_1, \ldots, r_l\}$ of relations in $\widehat{KQ}$,
we define the QP $(Q_A,W_A)$ by
\begin{itemize}
\item[$\bullet$] $Q_{A,0} = Q_0$,
\item[$\bullet$] $Q_{A,1} = Q_1 \coprod C_A$ with $C_A:=\{\rho_i: e(r_i)\to s(r_i)\ |\ 1\le i\le l\}$,
\item[$\bullet$] $W_A = \sum_{i=1}^l \rho_ir_i$.
\end{itemize}
\item[(b)]\cite{IO2} We define the \emph{complete 3-preprojective algebra} as the complete tensor algebra
\[\Pi_3(A):=\prod_{i\ge0}\Ext^2_A(DA,A)^{\otimes_Ai}\]
of the $(A,A)$-module $\Ext^2_A(DA,A)$.
\end{itemize}
\end{definition}

The following relationship between $(Q_A,W_A)$ and $\Pi_3(A)$ can be shown easily.

\begin{proposition}\label{preproj}(e.g. \cite[Thm. 6.10]{K2}) 
The $K$-algebras $\Pi_3(A)$ and $\P(Q_A,W_A)$ are isomorphic.
\end{proposition}

%

Next let $(Q,W)$ be a QP.
We denote by $\Gamma(Q,W)$ the \emph{complete Ginzburg DG algebra} \cite{G}\cite[Sect. 2.6]{KY}.
We have an isomorphism $H^0(\Gamma(Q,W)) \simeq \P(Q,W)$ of algebras \cite[Lem. 2.8]{KY}.
The \emph{generalized cluster category} of $(Q,W)$ is defined by
\[\C_{(Q,W)}:=\per\Gamma(Q,W)/\DD^{\rm b}\Gamma(Q,W),\]
where $\per\Gamma(Q,W)$ is the perfect derived category and $\DD^{\rm b}\Gamma(Q,W)$ is the bounded derived category (see \cite[Sect. 7]{KY} for details).

We say that a $K$-linear triangulated category $\C$ is \emph{2-Calabi-Yau} (or \emph{2-CY}) if 
each morphism space is finite dimensional over $K$ and there exists a functorial isomorphism
\[\Hom_{\C}(X,Y)\simeq D\Hom_{\C}(Y,X[2])\]
for any $X,Y\in\C$.
In this case, we say that an object $T\in\C$ is \emph{cluster tilting} if 
\begin{eqnarray*}
\add T=\{X\in\C\ |\ \Hom_{\C}(T,X[1])=0\}.
\end{eqnarray*}

\begin{proposition}\label{amiot}\cite{Am}\cite[Thm. 7.21]{KY}
If $\P(Q,W)$ is a finite dimensional $K$-algebra, then $\C_{(Q,W)}$ is a 2-Calabi-Yau triangulated category with a cluster tilting object $\Gamma(Q,W)$ such that $\End_{\C_{(Q,W)}}(\Gamma(Q,W))\simeq\P(Q,W)$.
\end{proposition}

For a finite dimensional $K$-algebra $A$ of global dimension at most two, we define the \emph{generalized cluster category} of $A$ by
\[\C_A:=\C_{(Q_A,W_A)}.\]
We have the following description of $\Gamma(Q_A,W_A)$ in terms of the \emph{complete derived 3-preprojective DG algebra} ${\bf \Pi}_3(A)$ \cite[Sect. 4.1]{K2}.

\begin{proposition}\label{derived preproj}\cite[Thm. 6.10]{K2}
${\bf \Pi}_3(A)$ is quasi-isomorphic to $\Gamma(Q_A,W_A)$.
\end{proposition}

In particular $\C_A$ is independent of choice of the presentation of $A$.
Moreover we recover Proposition \ref{preproj} by applying $H^0$ to the above quasi-isomorphism since we have $H^0({\bf \Pi}_3(A))=\Pi_3(A)$.




\section{Truncated Jacobian algebras and selfinjective QPs}

Let $(Q,W)$ be a QP. To each subset $C \subset Q_1$ we associate a grading $g_C$ on $Q$ by
\[
g_C(a) = \begin{cases} 1 & a \in C \\
0 & a \not \in C.
\end{cases}
\]
Denote by $Q_C$, the subquiver of $Q$ with vertex set $Q_0$ and arrow set $Q_1\setminus C$. 

\begin{definition} 
A subset $C \subset Q_1$ is called a \emph{cut} if $W$ is homogeneous of degree $1$ with respect to $g_C$. 
\end{definition}

If $C$ is a cut, 
then $g_C$ induces a grading on $\P(Q,W)$ as well.
The degree $0$ part of $\P(Q,W)$ is denoted by $\P(Q,W)_C$ and called the \emph{truncated Jacobian algebra}. We have
\[\P(Q,W)_C = \P(Q,W)\left/\overline{\langle C\rangle}\right. = \widehat{KQ_C}\left/ \overline{\langle \partial_c W \;|\; c \in C \rangle}. \right.\] 

\begin{definition} 
A cut $C$ is called \emph{algebraic} if the following conditions are satisfied.
\begin{itemize}
\item $\P(Q,W)_C$ is a finite dimensional $K$-algebra with global dimension at most two.
\item $\{\partial_cW\}_{c \in C}$ is a minimal set of generators of the ideal $\overline{\langle \partial_c W \;|\; c \in C \rangle}$ of $\widehat{KQ_C}$.
\end{itemize}
\end{definition}

Truncating Jacobian algebras by algebraic cuts can be viewed as a kind of inverse to taking $3$-preprojective algebras. More precisely we have the following result.

\begin{proposition}\label{cut preproj}
\begin{itemize}
\item[(a)] Let $A$ be a finite dimensional $K$-algebra with $\gl A\le 2$ and $(Q_A,W_A)$ be the associated QP. 
Then $C_A$ is an algebraic cut of $(Q_A,W_A)$ and $A \simeq \P(Q_A,W_A)_{C_A}$.
\item[(b)] Let $C$ be an algebraic cut in a QP $(Q,W)$ and set $A = \P(Q,W)_{C}$. Then $(Q_A,W_A)$ and $(Q,W)$ are isomorphic QPs.
\end{itemize}
\end{proposition}

\begin{proof}
(a) Assume that $A$ is as in Definition \ref{K2}. 
The assertion follows from the observations $(Q_A)_{C_A}= Q$ and $\partial_{\rho_i}W_A=r_i$ for $1\le i\le l$.

(b) By construction there is an isomorphism $Q_A \simeq Q$ sending $W_A$ to $W$.
\end{proof}

\begin{example}
\begin{itemize}
\item[(a)] Consider the following QP.
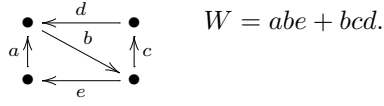
\[\xymatrix@R=.4cm@C=.4cm{
\bullet\ar^b[drr]&&\bullet\ar_d[ll]\\
\bullet\ar^a[u]&&\bullet\ar_c[u]\ar^e[ll]
}\ \ \ \ \ W=abe+bcd.\]
It has five cuts:
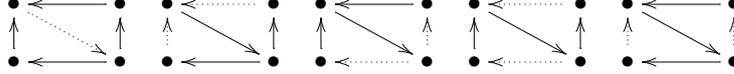
\[\xymatrix@R=.4cm@C=.4cm{
\bullet\ar@{.>}[drr]&&\bullet\ar[ll]\\
\bullet\ar[u]&&\bullet\ar[u]\ar[ll]
}\ \ \xymatrix@R=.4cm@C=.4cm{
\bullet\ar[drr]&&\bullet\ar@{.>}[ll]\\
\bullet\ar@{.>}[u]&&\bullet\ar[u]\ar[ll]
}\ \ \xymatrix@R=.4cm@C=.4cm{
\bullet\ar[drr]&&\bullet\ar[ll]\\
\bullet\ar[u]&&\bullet\ar@{.>}[u]\ar@{.>}[ll]
}\ \ \xymatrix@R=.4cm@C=.4cm{
\bullet\ar[drr]&&\bullet\ar@{.>}[ll]\\
\bullet\ar[u]&&\bullet\ar[u]\ar@{.>}[ll]
}\ \ \xymatrix@R=.4cm@C=.4cm{
\bullet\ar[drr]&&\bullet\ar[ll]\\
\bullet\ar@{.>}[u]&&\bullet\ar@{.>}[u]\ar[ll]
}\]
The last two cuts are not algebraic since $\P(Q,W)_C$
has global dimension three. The other cuts are algebraic.

\item[(b)] Next consider this QP:
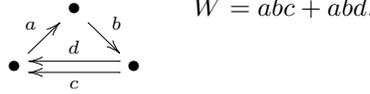
\[\xymatrix@R=.4cm@C=.4cm{
&\bullet\ar^b[dr]\\
\bullet\ar^a[ru]&&\bullet\ar@<.5ex>^c[ll]\ar@<-.5ex>_d[ll]
}\ \ \ \ \ W=abc+abd.
\]
It has three cuts:
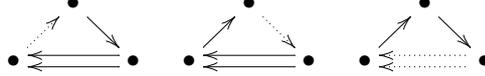
\[\xymatrix@R=.4cm@C=.4cm{
&\bullet\ar[dr]\\
\bullet\ar@{.>}[ru]&&\bullet\ar@<.5ex>[ll]\ar@<-.5ex>[ll]
}\ \ \ \xymatrix@R=.4cm@C=.4cm{
&\bullet\ar@{.>}[dr]\\
\bullet\ar[ru]&&\bullet\ar@<.5ex>[ll]\ar@<-.5ex>[ll]
}\ \ \ \xymatrix@R=.4cm@C=.4cm{
&\bullet\ar[dr]\\
\bullet\ar[ru]&&\bullet\ar@<.5ex>@{.>}[ll]\ar@<-.5ex>@{.>}[ll]
}\]
The last cut is not algebraic since $\partial_cW=\partial_dW$
contradicts minimality. The other cuts are algebraic.
\end{itemize}
\end{example}

\subsection{Selfinjective QPs}

Let $(Q,W)$ be a QP, and let $\Lambda = \P(Q,W)$ be the Jacobian algebra.
For each $i \in Q_0$, we denote by $P_i = \Lambda e_i$ the indecomposable projective left $\Lambda$-module,
by $I_i$ the maximal two-sided ideal of $\Lambda$,
and by $S_i:=\Lambda/I_i$ the simple $\Lambda$-bimodule. 
We denote by $(-)^*:=\Hom_\Lambda(-,\Lambda):\proj\Lambda\leftrightarrow\proj\Lambda^{\rm op}$ the duality between finitely generated projective modules.
Then $P_i^*=e_i\Lambda$ is the indecomposable projective right $\Lambda$-module.

It is known that Jacobian algebras satisfy the following Iwanaga-Gorenstein property.

\begin{proposition}\cite{BIRS}
If $\P(Q,W)$ is finite dimensional, then
\[\injdim_{\P(Q,W)}\P(Q,W)=\injdim_{\P(Q,W)^{\rm op}}\P(Q,W)\le 1.\]
\end{proposition}

This motivates us to introduce the following class of QPs, which are main subject of this paper.

\begin{definition}
\begin{itemize}
\item[(a)] We say that a QP $(Q,W)$ is \emph{finite dimensional} if
$\P(Q,W)$ is a finite dimensional $K$-algebra.
\item[(b)] We say that a QP $(Q,W)$ is \emph{selfinjective} if $\P(Q,W)$
is a finite dimensional selfinjective $K$-algebra. 
\end{itemize}
\end{definition}

Since 
\[
\sum_{a\in Q_1}a(\partial_{(a,b)} W) = \partial_{b}W\ \mbox{ and }\ 
\sum_{b\in Q_1}(\partial_{(a,b)} W)b = \partial_{a}W
\]
we have the following complexes of left respectively right $\Lambda$-modules:
\begin{eqnarray}\label{resolution1}
{\displaystyle
P_i \xrightarrow{[b]} \bigoplus_{\begin{smallmatrix}b\in Q_1\\ s(b)=i\end{smallmatrix}} P_{e(b)} 
\xrightarrow{[\partial_{(a,b)}W]} \bigoplus_{\begin{smallmatrix}a\in Q_1\\ e(a)=i\end{smallmatrix}} P_{s(a)} \xrightarrow{[a]} 
P_i\xrightarrow{} S_i \xrightarrow{} 0}\\ 
\label{resolution2}
{\displaystyle
P^*_i \xrightarrow{[a]} \bigoplus_{\begin{smallmatrix}a\in Q_1\\ e(a)=i\end{smallmatrix}} P^*_{s(a)}
\xrightarrow{[\partial_{(a,b)}W]} \bigoplus_{\begin{smallmatrix}b\in Q_1\\ s(b)=i\end{smallmatrix}} P^*_{e(b)} 
\xrightarrow{[b]} P^*_i\xrightarrow{} S_i \xrightarrow{} 0}
\end{eqnarray}
which are exact everywhere except possibly at $\bigoplus P_{e(b)} $ and $\bigoplus P^*_{s(a)}$ respectively. These sequences can be used to characterize selfinjective QPs.

\begin{theorem}\label{characterization}
The following conditions are equivalent for a QP $(Q,W)$.
\begin{itemize}
\item[(a)] $(Q,W)$ is selfinjective.
\item[(b)] $(Q,W)$ is finite dimensional and \eqref{resolution1} is exact.
\item[(c)] $(Q,W)$ is finite dimensional and \eqref{resolution2} is exact.
\end{itemize}
\end{theorem}

\begin{proof}
We only prove that (a) is equivalent to (b).

Condition (a) holds if and only if $\Ext^1_\Lambda(S_i,\Lambda)=0$ for all $i\in Q_0$.
This is equivalent to that we get an exact sequence by applying $\Hom_\Lambda(-,\Lambda)$ to the first part 
\[{\displaystyle
\bigoplus_{\begin{smallmatrix}a\in Q_1\\ e(a)=i\end{smallmatrix}} P^*_{s(a)}
\xrightarrow{[\partial_{(a,b)}W]} \bigoplus_{\begin{smallmatrix}b\in Q_1\\ s(b)=i\end{smallmatrix}} P^*_{e(b)} 
\xrightarrow{[b]} P^*_i\xrightarrow{} S_i \xrightarrow{} 0}
\]
of \eqref{resolution2}.
This is equivalent to that \eqref{resolution1} is exact.
\end{proof}

Immediately we have the following result.

\begin{proposition}\label{sc}
If a QP $(Q,W)$ is selfinjective, then each arrow in $Q$ appears in $W$.
\end{proposition}

\begin{proof}
Assume that there is an arrow $b:i\to e(b)$ that does not appear in $W$.
Then the middle map $[\partial_{(a,b)}W]$ in \eqref{resolution1} is zero on
the direct summand $P_{e(b)}$.
This contradicts to the exactness of \eqref{resolution1} since
$b:P_i\to P_{e(b)}$ is not surjective.
\end{proof}

\subsection{Structure theorem of 2-representation-finite algebras}

For a finite dimensonal algebra $A$, we say that $M\in\mod A$ is \emph{cluster tilting} \cite{I1,I3} if
\begin{eqnarray*}
\add M&=&\{X\in\mod A\ |\ \Ext^1_A(M,X)=0\}\\
&=&\{X\in\mod A\ |\ \Ext^1_A(X,M)=0\}.
\end{eqnarray*}
It is known that the category $\add M$ is a higher analogue
of the module category $\mod A$ in the sense that there exist
2-AR translation functors
\[\xymatrix{\add M\ar@<.5ex>^{\tau_2=D\Ext_A^2(-,A)}[rrr]&&&\add M\ar@<.5ex>^{\tau_2^-=\Ext_A^2(DA,-)}[lll]}\]
which give 2-almost split sequences in $\add M$.

We say that an algebra $A$ is \emph{2-representation-finite} if
there exists a cluster tilting $A$-module and $\gl A\le 2$.
In this case it is known that $\Pi_3(A)$ regarded as an $A$-module
gives a unique basic 2-cluster tilting $A$-module \cite{I4}.


We have the following criterion for 2-representation-finiteness.

\begin{proposition}\label{selfinj}
Let $A$ be a $K$-algebra such that $\gl A \leq 2$. Then the following conditions are equialent.
\begin{itemize}
\item[(a)] $A$ is $2$-representation-finite.
\item[(b)] $\Pi_3(A)$ is a finite dimensional selfinjective algebra.
\item[(c)] $(Q_A,W_A)$ is a selfinjective QP.
\end{itemize}
\end{proposition}

\begin{proof}
The equivalence of (a) and (b) is shown in \cite{IO2}.
The equivalence of (b) and (c) follows from Proposition \ref{preproj}.
\end{proof}

Propositions \ref{cut preproj} and \ref{selfinj} give a strong motivation for studying selfinjective QPs and their algebraic cuts. To do the latter, the following result is very useful.

\begin{proposition}\label{selfinjQP} 
Let $(Q,W)$ be a selfinjective QP. Then every cut $C \subset Q_1$ is algebraic.
\end{proposition}

\begin{proof}
We consider  $\Lambda:=\P(Q,W)$ to be graded by $g_C$. 
Let $A:=\P(Q,W)_C$. For each $i \in Q_0$ set $P_i = \Lambda e_i$ and $\overline{P}_i = A e_i$. 
Then $P_i$ is a graded $\Lambda$-module and the part in degree $0$ is $(P_i)_0 = \overline{P}_i$. Fix $i \in Q_0$. 
By Proposition \ref{characterization} we have an exact sequence
\[
{\displaystyle
P_i \xrightarrow{[b]} \bigoplus_{\begin{smallmatrix}b\in Q_1\\ s(b)=i\end{smallmatrix}} P_{e(b)} 
\xrightarrow{[\partial_{(a,b)}W]} \bigoplus_{\begin{smallmatrix}a\in Q_1\\ e(a)=i\end{smallmatrix}} P_{s(a)} \xrightarrow{[a]} 
P_i\xrightarrow{} S_i \xrightarrow{} 0.}
\]
Note that $\partial_{(a,b)}W$ is homogeneous of degree $1-g_C(a) - g_C
(b)$. Thus taking appropriate degree shifts we obtain the following exact sequence of graded $\Lambda$-modules and homogeneous morphisms of degree $0$.
\[
P_i(-1) \xrightarrow{\left[\begin{smallmatrix} [b] \\ [c] \end{smallmatrix}\right]} 
\begin{array}{l}
{\displaystyle\bigoplus_{\begin{smallmatrix}b \not\in C\\ s(b)=i\end{smallmatrix}} P_{e(b)}(-1)}
\\
\oplus {\displaystyle\bigoplus_{\begin{smallmatrix}c\in C\\ s(c)=i\end{smallmatrix}} P_{e(c)}}
\end{array}
\xrightarrow{\left[\begin{smallmatrix} [\partial_{(d,b)}W] & 0 \\ [\partial_{(a,b)}W]  & [\partial_{(a,c)}W] \end{smallmatrix}\right]}
\begin{array}{l}
\displaystyle{\bigoplus_{\begin{smallmatrix}d\in C \\ e(d)=i\end{smallmatrix}} P_{s(d)}(-1)}
\\
\oplus {\displaystyle\bigoplus_{\begin{smallmatrix}a \not\in C\\ e(a)=i\end{smallmatrix}} P_{s(a)}}
\end{array}
\xrightarrow{\left[\begin{smallmatrix} [d] & [a] \end{smallmatrix}\right]} 
P_i\xrightarrow{} S_i \xrightarrow{} 0.
\]
Taking the degree $0$ part of the above sequence we get an exact sequence
\[
0 \xrightarrow{}
{\displaystyle\bigoplus_{\begin{smallmatrix}c\in C\\ s(c)=i\end{smallmatrix}} \overline{P}_{e(c)}}
\xrightarrow{[\partial_{(a,c)}W]}
{\displaystyle\bigoplus_{\begin{smallmatrix}a \not\in C\\ e(a)=i\end{smallmatrix}} \overline{P}_{s(a)}}
\xrightarrow{[a]} 
\overline{P}_i\xrightarrow{} S_i \xrightarrow{} 0.
\]
Hence $\gl A \leq 2$. Moreover since the above sequence is exact it is also minimal and thus $C$ is algebraic.
\end{proof}

We now give a structure theorem of $2$-representation-finite algebras in terms of selfinjective QPs.

\begin{theorem}\label{QP}
\begin{itemize}
\item[(a)] For any selfinjective QP $(Q,W)$ and cut $C$, the truncated Jacobian algebra $\P(Q,W)_C$ is $2$-representation-finite.
\item[(b)] Every basic $2$-representation-finite algebra appears in this way.
\end{itemize}
\end{theorem}

In particular truncated Jacobian algebras of selfinjective QPs are exactly basic $2$-representation-finite algebras.

\begin{proof}
(a) By Proposition \ref{selfinjQP} we have that $C$ is algebraic and so $A := \P(Q,W)_C$ has global dimension at most two. By Proposition \ref{cut preproj}(b), we have 
\[
\P(Q_A,W_A) \simeq \P(Q,W),\]
which is selfinjective.
Therefore, by Proposition \ref{selfinj}, $A$ is $2$-representation-finite.

(b) Up to Morita equivalence every $2$-representation-finite algebra is of the form 
\[
A = KQ/\langle r_1, \ldots, r_n \rangle = \widehat{KQ}\left/ \overline{\langle r_1, \ldots, r_n \rangle}\right. 
\]
for some minimal set of admissible relations $\{r_1, \ldots, r_n\}$. 
By Proposition \ref{selfinj}, $(Q_A,W_A)$ is selfinjective. Moreover, $C_A$ is a cut and $A \simeq \P(Q_A,W_A)_{C_A}$ holds by Proposition \ref{cut preproj}(a).
\end{proof}

Theorem \ref{QP} reduces the description of $2$-representation-finite algebras to truncated Jacobian algebras $\P(Q,W)_C$ of selfinjective QPs $(Q,W)$ and their cuts $C$.

In order to construct $2$-representation-finite algebras we want to find selfinjective QPs. 
In the next section we show that new selfinjective QPs can be constructed from a given one by mutation.

\section{Mutation and selfinjectivity}

\subsection{Mutation of selfinjective QPs}

Let $(Q,W)$ be a QP and $\Lambda:=\P(Q,W)$. When $k\in Q_0$ is not contained
in any $2$-cycles, we have a new QP $\mu_k(Q,W)$ called the \emph{mutation} of $(Q,W)$ at $k$ (see \cite{DWZ1}).
Simple examples show that mutation does not preserve selfinjectivity of QPs.
In this section we define successive mutation with respect to orbits of the Nakayama permutation (see also \cite{BO}), and show that it preserves selfinjectivity.

\begin{definition}\label{succ mutation}
Let $(Q,W)$ be a selfinjective QP.

We denote by $\sigma:Q_0\to Q_0$ the Nakayama permutation,
so $D(e_i\Lambda)\simeq\Lambda e_{\sigma i}$ for any $i\in Q_0$.

For $k\in Q_0$, we denote by $(k) = \{k=\sigma^mk,\sigma k,\ldots,\sigma^{m-1}k\}$
the $\sigma$-orbit of $k$. We assume that
\begin{itemize}
\item[$\bullet$] $k$ is not contained in 2-cycles in $Q$,
\item[$\bullet$] there are no arrows between $k,\sigma k,\ldots,\sigma^{m-1}k$ in $Q$.
\end{itemize}
In this case we define the successive mutation
\[\mu_{(k)}(Q,W):=\mu_{\sigma^{m-1}k}\circ\cdots\circ\mu_k(Q,W)\]
with respect to the orbit $(k)$.
This is well-defined since $\sigma^ik$ is not contained in 2-cycles in the quiver of $\mu_{\sigma^{i-1}k}\circ\cdots\circ\mu_k(Q,W)$ for any $0\le i<m$ by our assumption.
Moreover $\mu_{(k)}(Q,W)$ is independent of the choice of order of mutations $\mu_k,\ldots,\mu_{\sigma^{m-1}k}$.
\end{definition}

The following key observation will be proved in the next section.

\begin{theorem}\label{selfinjective QP}
Let $(Q,W)$ be a selfinjective QP with the Nakayama permutation $\sigma:Q_0\to Q_0$.
Assume that $k\in Q_0$ satisfies the two conditions in Definition \ref{succ mutation}.
\begin{itemize}
\item[(a)] The above mutation $\mu_{(k)}(Q,W)$ is again a selfinjective QP.
\item[(b)] The Nakayama permutation of $\mu_{(k)}(Q,W)$ is again given by $\sigma$.
\end{itemize}
\end{theorem}

\subsection{Selfinjective cluster tilting objects}

Throughout this section, let $\C$ be a 2-CY triangulated category (see Section \ref{preliminary2}).
We consider a special class of cluster tilting objects in $\C$ and their mutation.

\begin{definition}
We say that a cluster tilting object $T\in\C$ is \emph{selfinjective} if $\End_{\C}(T)$ is a
finite dimensional selfinjective $K$-algebra.
\end{definition}

We have the following criterion for selfinjectivity.

\begin{proposition}\cite{IO2}\label{IO2}
Let $T=T_1\oplus\cdots\oplus T_n\in\C$ be a basic cluster tilting object with indecomposable summands $T_i$.
\begin{itemize}
\item[(a)] $T$ is selfinjective if and only if $T\simeq T[2]$.
\item[(b)] In this case we define a permutation $\sigma:\{1,\ldots,n\}\to\{1,\ldots,n\}$ by $T_i[2]\simeq T_{\sigma i}$.
Then $\sigma$ gives the Nakayama permutation of $\End_{\C}(T)$.
\end{itemize}
\end{proposition}

In the rest of this section we shall show that selfinjectivity of cluster tilting objects is preserved under certain successive cluster tilting mutation.

Let $T\in\C$ be a basic cluster tilting object, and let $T=U\oplus V$ be a decomposition.
We take triangles (called \emph{exchange triangles})
\begin{equation}\label{app}
\xymatrix{U^*\ar[r]&V'\ar^f[r]&U\ar[r]&U^*[1]}\ \mbox{ and }\ 
\xymatrix{U\ar^g[r]&V''\ar[r]&{}^*U\ar[r]&U[1]}
\end{equation}
with a minimal right $(\add V)$-approximation $f$ of $U$
and a minimal left $(\add V)$-approximation $g$ of $U$. Let
\[\mu^+_U(T):=U^*\oplus V\ \mbox{ and }\ \mu^-_U(T):={}^*U\oplus V.\]
Notice that $U$ is not assumed to be indecomposable.
A typical case of cluster tilting mutation is given by
$\mu^+_T(T)=T[-1]$ and $\mu^-_T(T)=T[1]$.

We have the following basic results.

\begin{proposition}\cite{IY}\label{basic results}
\begin{itemize}
\item[(a)] $\mu^+_U(T)$ and $\mu^-_U(T)$ are basic cluster tilting objects in $\C$.
\item[(b)] If $U$ is indecomposable, then $\mu^+_U(T)$ and $\mu^-_U(T)$
are isomorphic, which we denote by $\mu_U(T)$.
\end{itemize}
\end{proposition}


We have the following general result.

\begin{proposition}\label{selfinjective mutation}
If $U\simeq U[2]$, then $\mu^+_U(T)$ and $\mu^-_U(T)$ are selfinjective cluster tilting objects in $\C$.
\end{proposition}

\begin{proof}
Clearly we have $V\simeq V[2]$. Let $a:U\to U[2]$ be an isomorphism and
\[\xymatrix{U^*\ar[r]&V'\ar^{f}[r]&U\ar[r]&U^*[1]}\]
be an exchange triangle. Then $fa:V'\to U[2]$
is a minimal right $(\add V)$-approximation.
Since $V\simeq V[2]$, we have that $f[2]:V'[2]\to U[2]$ is also
a minimal right $(\add V)$-approximation.
Thus we have a commutative diagram
\[\xymatrix{
U^*\ar[r]&V'\ar^f[r]\ar^b[d]&U\ar[r]\ar^a[d]&U^*[1]\\
U^*[2]\ar[r]&V'[2]\ar^{f[2]}[r]&U[2]\ar[r]&U^*[3]
}\]
of triangles with isomorphisms $a$ and $b$.
By an axiom of triangulated categories, there exists
an isomorphism $c:U^*\to U^*[2]$ which keeps the above diagram
commutative. In particular we have $U^*\simeq U^*[2]$ and so $\mu^+_U(T)\simeq\mu^+_U(T)[2]$.
Thus $\mu^+_U(T)$ is a selfinjective cluster tilting object by Proposition \ref{IO2}.
\end{proof}

Now let us observe $\mu^+_U(T)$ and $\mu^-_U(T)$ more explicitly.
We take decompositions $U=T_1\oplus\cdots\oplus T_m$ and $V=T_{m+1}\oplus\cdots\oplus T_n$ into indecomposable summands.
For $i=1,\ldots,m$, we take triangles
\begin{equation*}
\xymatrix{T_i^\vee\ar[r]&V_i'\ar^{f_i}[r]&T_i\ar[r]&T_i^\vee[1]}\ \mbox{ and }\ 
\xymatrix{T_i\ar^{g_i}[r]&V_i''\ar[r]&{}^\vee T_i\ar[r]&T_i[1]}
\end{equation*}
with a minimal right $(\add V)$-approximation $f_i$ of $T_i$
and a minimal left $(\add V)$-approximation $g_i$ of $T_i$.
Then be definition we have decompositions
\begin{equation}\label{decomposition1}
\mu^+_U(T)\simeq T_1^\vee\oplus\cdots\oplus T_m^\vee\oplus T_{m+1}\cdots\oplus T_n\ \mbox{ and }\ 
\mu^-_U(T)\simeq{}^\vee T_1\oplus\cdots\oplus{}^\vee T_m\oplus T_{m+1}\cdots\oplus T_n
\end{equation}
into indecomposable summands.

\begin{proposition}\label{selfinjective mutation2}
In Proposition \ref{selfinjective mutation},
the Nakayama permutations of $\mu^+_U(T)$ and $\mu^-_U(T)$ with respect to the decompositions in \eqref{decomposition1} coincide with that of $T$.
\end{proposition}

\begin{proof}
Using an isomorphism $T_{\sigma i}\to T_i[2]$, we get isomorphisms
$T_{\sigma i}^\vee\to T_i^\vee[2]$ and ${}^\vee T_{\sigma i}\to {}^\vee T_i[2]$
by a similar argument as in the proof of Proposition \ref{selfinjective mutation}.
Thus the assertion follows.
\end{proof}

Although we do not assume that $U$ is indecomposable, sometimes we can decompose $\mu_U$
as the successive cluster tilting mutation with respect to indecomposable summands.

\begin{proposition}\label{successive mutation}
If there are no arrows between vertices $T_1,\ldots,T_m$ in the quiver of $\End_{\C}(T)$, then we have
$\mu^+_U(T)\simeq\mu^-_U(T)\simeq\mu_{T_m}\circ\cdots\circ\mu_{T_1}(T)$.
\end{proposition}

\begin{proof}
Inductively we shall show that $\mu_{T_k}\circ\cdots\circ\mu_{T_1}(T)$ is isomorphic to
\[T^{(k)}:=T_1^\vee\oplus\cdots\oplus T_k^\vee\oplus T_{k+1}\oplus\cdots\oplus T_n\]
for any $1\le k\le m$.
We assume that this is true for $k=i-1$, and we shall show that $\mu_{T_i}(T^{(i-1)})\simeq T^{(i)}$.
We only have to show that a right $(\add V)$-approximation of $T_i$ is also a right $(\add T^{(i-1)}/T_i)$-approximation.
This is true since we have a triangle $T_k^\vee\to V_k'\to T_k\to T_k^\vee[1]$, which implies that any morphism $T_k^\vee\to T_i$ factors through $V_k'\in\add V$.

Consequently we have $\mu_{T_m}\circ\cdots\circ\mu_{T_1}(T)\simeq T_1^\vee\oplus\cdots\oplus T_m^\vee\oplus T_{m+1}\oplus\cdots\oplus T_n\simeq\mu^+_U(T)$.
\end{proof}

The following result shows compatibility of cluster tilting mutation and QP mutation.

\begin{proposition}\label{birsm}\cite{BIRS}
Let $\C$ be a 2-CY triangulated category and $T=T_1\oplus\cdots\oplus T_n\in\C$ a basic cluster
tilting object with indecomposable summands $T_i$. If $\End_{\C}(T)\simeq\P(Q,W)$ for a QP $(Q,W)$ and
$k\in Q_0$ is not contained in 2-cycles in $Q$, then
$\End_{\C}(\mu_{T_k}(T))\simeq\P(\mu_k(Q,W))$.
\end{proposition}

Now we are ready to prove Theorem \ref{selfinjective QP}.

Let $\C:=\C_{(Q,W)}$ be a generalized cluster category in Section \ref{preliminary2}.
Then there exists a cluster tilting object $T\in\C$ such that $\End_{\C}(T)\simeq\P(Q,W)$ by Proposition \ref{amiot}.
Let $U$ be the summand of $T$ corresponding to the vertices $k,\sigma k,\ldots,\sigma^{m-1}k$.
By Proposition \ref{IO2}(b), we have $U\simeq U[2]$.
By Proposition \ref{selfinjective mutation}, we know that
$T':=\mu^+_U(T)$ is a selfinjective cluster tilting object.
Thus $\End_{\C}(T')$ is a selfinjective algebra.
On the other hand, $T'$ is isomorphic to
$\mu_{T_{\sigma^{m-1}k}}\circ\cdots\circ\mu_{T_k}(T)$
by Proposition \ref{successive mutation}.
Applying Proposition \ref{birsm} repeatedly, we have
\[\End_{\C}(T')\simeq\P(\mu_{\sigma^{m-1}k}\circ\cdots\circ\mu_k(Q,W)),\]
so we have the assertion.

The statement for Nakayama permutation follows from Proposition \ref{selfinjective mutation2}.
\qed

\section{Examples of selfinjective QPs}\label{Examples of selfinjective QPs}

\subsection{Examples from cluster categories}\label{ex from ctilt}
Cluster tilted algebras have been shown to be Jacobian algebras of QPs \cite{BIRS,K2}. The selfinjective ones have been classified by Ringel \cite{Rin}. 
The classification consists of the two families described below.

For each $n>3$ define the following QP
\[
\begin{xy} 0;<0.38pt,0pt>:<0pt,-0.38pt>:: 
(0,55) *+{n} ="0",
(50,0) *+{1} ="1",
(115,0) *+{2} ="2",
(160,55) *+{3} ="3",
(160,115) *+{4} ="4",
(115,165) *+{5} ="5",
(0,115) *+{n-1} ="6",
(60,140) *+{\ddots} ="7",
(-70,80) *+{Q^n:} ="8",
(340,80) *+{W^n=a_1\cdots a_n} ="8",
"1", {\ar"0"_{a_1}},
"0", {\ar"6"_{a_n}},
"2", {\ar"1"_{a_2}},
"3", {\ar"2"_{a_3}},
"4", {\ar"3"_{a_4}},
"5", {\ar"4"_{a_5}},
\end{xy}
\]
If $n$ is even then define also
\[
\begin{xy} 0;<0.38pt,0pt>:<0pt,-0.38pt>:: 
(0,55) *+{n} ="0",
(50,0) *+{1} ="1",
(115,0) *+{2} ="2",
(160,55) *+{3} ="3",
(160,115) *+{4} ="4",
(115,165) *+{5} ="5",
(0,115) *+{n-1} ="6",
(60,140) *+{\ddots} ="7",
(-100,80) *+{\tilde{Q}^{n}:} ="8",
(450,80) *+{\tilde{W}^n=a_1a_3\cdots a_{n-1} - \sum a_{i}b_{i-2}b_{i-1}} ="8",
"0", {\ar"1"^{b_n}},
"6", {\ar"0"^{b_{n-1}}},
"1", {\ar"2"^{b_1}},
"3", {\ar"1"^{a_3}},
"1", {\ar"6"^{a_1}},
"2", {\ar"3"^{b_2}},
"3", {\ar"4"^{b_3}},
"5", {\ar"3"^{a_5}},
"4", {\ar"5"^{b_4}},
\end{xy}
\]
These QPs are selfinjective with Nakayama automorphism induced by $i \mapsto i-2$.
Notice that mutating at the even vertices in $(\tilde{Q}^{n},\tilde{W}^{n})$ gives $(Q^{n},W^{n})$.
This is not surprising as the result must again be cluster tilted by \cite{BIRS} and 
selfinjective by Theorem \ref{selfinjective QP}. By \cite{Rin} there is only one alternative.

\medskip
Next we consider cluster categories associated with canonical algebras.
Selfinjective cluster tilting objects are classified in \cite{HILO} (see also \cite{BG}).
Let us give one example.
It is shown in \cite{IO2} that the tubular algebra of type $(2,2,2,2)$ given below is
2-representation-finite.
\[\xymatrix@R=.4cm@C=.4cm{
&&\bullet\ar_a[dll]\ar^b[dl]\ar_c[dr]\ar^d[drr]\\
\bullet\ar_{a'}[drr]&\bullet\ar^{b'}[dr]&&\bullet\ar_{c'}[dl]&\bullet\ar^{d'}[dll]&\ \ \ aa'+bb'+cc'=0=aa'+\lambda bb'+dd'\\
&&\bullet
}\]
In particular, we have the following selfinjective QP.
\[\xymatrix@R=.4cm@C=.4cm{
&&\bullet\ar_a[dll]\ar_b[dl]\ar^c[dr]\ar^d[drr]\\
\bullet\ar_{a'}[drr]&\bullet\ar_{b'}[dr]&&\bullet\ar^{c'}[dl]&\bullet\ar^{d'}[dll]&\ \ \ W=aa'e+bb'e+cc'e+aa'f+\lambda bb'f+dd'f\\
&&\bullet\ar@<.5ex>^e[uu]\ar@<-.5ex>_f[uu]
}\]
The Nakayama permutation is the identity.
Mutating at the leftmost vertex, we have the following selfinjective QP with $\lambda'=\frac{\lambda}{\lambda-1}$.
\[\xymatrix@R=.4cm@C=.4cm{
&&\bullet\ar_b[dl]\ar^c[dr]\ar^d[drr]\\
\bullet\ar^{a}[urr]&\bullet\ar_{b'}[dr]&&\bullet\ar^{c'}[dl]&\bullet\ar^{d'}[dll]&\ \ \ W=bb'e+cc'e+dd'e+\lambda'bb'a'a+dd'a'a\\
&&\bullet\ar^{a'}[ull]\ar^e[uu]
}\]

\subsection{Examples from pairs of Dynkin quivers}

In this section we construct a new class of selfinjective QPs.

Let $Q^1$ and $Q^2$ be finite quivers without oriented cycles. Define $Q = Q^1 \tildetensor Q^2$ by
\[
Q_0 = Q^1_0 \times Q^2_0
\]
and
\[
Q_1 = (Q^1_0 \times Q^2_1) \coprod (Q^1_1 \times Q^2_0) \coprod (Q^1_1\times Q^2_1),
\]
where 
\[\begin{matrix}
s(a, y) = (s(a), y), & s(x, b) = (x, s(b)), & s(a, b) = (e(a), e(b)), \\
e(a, y) = (e(a), y), & e(x, b) = (x, e(b)), & e(a, b) = (s(a), e(b))
\end{matrix}\]
for any $x\in Q^1_0$, $y\in Q^2_0$, $a\in Q^1_1$ and $b\in Q^2_1$.
Also define the potential
\[
W = W^{\tildetensor}_{Q^1,Q^2} = \sum_{\begin{smallmatrix}a \in Q^1_1 \\ b \in Q^2_1\end{smallmatrix}}(a,e(b))(s(a),b)(a,b) - (e(a),b)(a,s(b))(a,b).
\]
By the results in \cite{K3}, $C = Q^1_1\times Q^2_1$ is an algebraic cut and there is an algebra isomorphism
\[
\P(Q,W)_C \simeq KQ^1 \otimes KQ^2
\]
mapping $(x, y) \mapsto e_x \otimes e_y$, $(a, y) \mapsto a \otimes e_y$ and $(x, b) \mapsto e_x \otimes b$. Hence $\Pi_3(KQ^1 \otimes KQ^2) \simeq \P(Q,W)$.
As is common we denote $Q_C$ by $Q^1 \otimes Q^2$.

\medskip
Let us consider Dynkin diagrams:
\[\xymatrix@C0.4cm@R0.2cm{
A_n&1\ar@{-}[r]&2\ar@{-}[r]&3\ar@{-}[r]&\ar@{..}[rrr]&&&\ar@{-}[r]&n-2\ar@{-}[r]&n-1\ar@{-}[r]&n\\
&&&&&&&&n\\
D_n&1\ar@{-}[r]&2\ar@{-}[r]&3\ar@{-}[r]&\ar@{..}[rrr]&&&\ar@{-}[r]&n-2\ar@{-}[r]\ar@{-}[u]&n-1\\
&&&4\\
E_6&1\ar@{-}[r]&2\ar@{-}[r]&3\ar@{-}[r]\ar@{-}[u]&5\ar@{-}[r]&6\\
&&&&7\\
E_7&1\ar@{-}[r]&2\ar@{-}[r]&3\ar@{-}[r]&4\ar@{-}[r]\ar@{-}[u]&5\ar@{-}[r]&6\\
&&&&&8\\
E_8&1\ar@{-}[r]&2\ar@{-}[r]&3\ar@{-}[r]&4\ar@{-}[r]&5\ar@{-}[r]\ar@{-}[u]&6\ar@{-}[r]&7}\]
We define a canonical involution $\omega$ of each Dynkin diagram as follows:
\begin{itemize}
\item For $A_n$, we put $\omega(i)=n+1-i$.
\item For $D_n$ with odd $n$, we put $\omega(n-1)=n$, $\omega(n)=n-1$ and $\omega(i)=i$ for other $i$.
\item For $E_6$, we put $\omega(1)=6$, $\omega(2)=5$, $\omega(5)=2$, $\omega(6)=1$ and $\omega(i)=i$ for other $i$.
\item For other types, we put $\omega=\id$.
\end{itemize}
We call a Dynkin quiver \emph{stable} if it stable under $\omega$.

The \emph{Coxeter number} $h$ of each Dynkin diagram is given as follows:
\[\begin{array}{|c|c|c|c|c|}
\hline
A_n&D_n&E_6&E_7&E_8\\ \hline
n+1&2(n-1)&12&18&30\\ \hline
\end{array}\]

\begin{proposition}\label{tensor}
Let $Q^1$ and $Q^2$ be Dynkin quivers with canonical involutions $\omega_1$ respectively $\omega_2$. 
Assume that they are stable and have the same Coxeter numbers. 
Then $(Q^1 \tildetensor Q^2,W^{\tildetensor}_{Q^1,Q^2})$ is a selfinjective QP with Nakayama permutation induced by $\omega_1\times \omega_2$.
\end{proposition}

\begin{proof}
Under the assumptions $KQ^1 \otimes KQ^2$ is $2$-representation-finte by \cite{HI}. Since 
\[
\P(Q^1 \tildetensor Q^2,W^{\tildetensor}_{Q^1,Q^2}) \simeq \Pi_3(KQ^1 \otimes KQ^2),
\]
Proposition \ref{selfinj} implies that $(Q^1 \tildetensor Q^2,W^{\tildetensor}_{Q^1,Q^2})$ is selfinjective.
Moreover, by the results in \cite{HI} the Nakayama permutation is induced by $\omega_1\times \omega_2$.
\end{proof}

We use the class of selfinjective QPs in Proposition \ref{tensor} to construct another one by mutation. 

Let $Q^1$ and $Q^2$ be Dynkin quivers with alternating orientation (i.e., each vertex is either a sink or a source). 
Write $Q^i_1 = X_i \coprod Y_i$ where $X_i$ consists of sources and $Y_i$ of sinks, i.e, for each arrow $x \arr{a} y\in Q^i_1$ holds that $x \in X_i$ and $y \in Y_i$.

Recall from \cite{K3}, that the square product $Q^1\square Q^2$ is the quiver obtained from $Q^1 \otimes Q^2$ 
by replacing each arrow $a$ satisfying $s(a) \in Y_1\times X_2$ or $e(a) \in Y_1\times X_2$ with an arrow $a^*$ of opposite orientation. 
Let $W^{\square}_{Q^1,Q^2}$ be the potential defined by the sum of all cycles of the form
\[
(a,y')(x,b)(a,x')^*(y,b)^*
\]
where $x \arr a y \in Q^1_1$ and $x' \arr b y' \in Q^2_1$.

\begin{theorem}\label{square}
Let $Q^1$ and $Q^2$ be Dynkin quivers with alternating orientation such that their Coxeter numbers coincide. 
Let $\omega_1$ and $\omega_2$ be the corresponding canonical involutions. 
Then $(Q^1 \square Q^2,W^{\square}_{Q^1,Q^2})$ is a selfinjective QP with Nakayama permutation induced by $\omega_1 \times \omega_2$.
\end{theorem}

\begin{proof}
Let $Y_1\times X_2 = \{i_1\ldots i_s\}$, and set
\[
(Q,W) = \mu_{i_s}\cdots \mu_{i_1} (Q^1 \square Q^2,W^{\square}_{Q^1,Q^2}).
\]
We proceed to show that $(Q,W) \simeq (Q^1 \tildetensor Q^2, W^{\tildetensor}_{Q^1,Q^2})$.

The paths of length two through $(y,x') \in Y_1\times X_2$ are exactly those of the form $(a,x')^*(y,b)^*$, 
where $x \arr a y \in Q^1_1$ and $x' \arr b y' \in Q^2_1$ for some $x\in X_1$ and $y' \in Y_2$. Hence
\[
W = \sum_{a,b} (a,y')(x,b)[(a,x')^*(y,b)^*] + (y,b)^{**}(a,x')^{**}[(a,x')^*(y,b)^*].
\]
Now consider the isomorphism $\widehat{KQ} \simeq \widehat{K(Q^1 \tildetensor Q^2)}$ defined by $[(a,x')^*(y,b)^*] \mapsto (a,b)$, $(a,x')^{**}\mapsto (a,x')$, $(y,b)^{**} \mapsto -(y,b)$ and identity on all other arrows. 
It maps the potential $W$ to
\[
\sum (a,y')(x,b)(a,b) - (y,b)(a,x')(a,b).
\]
Thus the claim is proved.

Now observe that if the common Coxeter number is even, then every alternating orientation is in fact stable. 
Hence $(Q^1 \tildetensor Q^2, W^{\tildetensor}_{Q^1,Q^2})$ is selfinjective by Proposition \ref{tensor}. 
Moreover, $Y_1\times X_2$ is invariant under the Nakayama permutation which is induced by $\omega_1\times \omega_2$. Hence
\[
(Q^1 \square Q^2,W^{\square}_{Q^1,Q^2}) \simeq \mu_{i_1}\cdots \mu_{i_s} (Q,W) 
\]
is selfinjective with Nakayama permutation induced by $\omega_1\times \omega_2$ by Theorem \ref{selfinjective QP}.

It remains to consider the case when the Coxeter number is odd. It is covered by the following proposition.
\end{proof}

\begin{proposition}\label{coxeter_odd}
Let $Q$ be a Dynkin quiver of type $\mathbb{A}_{2s}$ with alternating orientation. 
Then $(Q \square Q,W^{\square}_{Q,Q})$ is a selfinjective QP with Nakayama permutation induced by $\omega \times \omega$.
\end{proposition}

\begin{proof}

We label the arrows in $Q$ as follows:
\[
\xymatrix{1 \ar[r]^{a_1} & 2 & 3 \ar[l]_{b_3} \ar[r]^{a_3} & \cdots \ar[r]^{a_{2s-1}} & 2s}.
\]
Thus $Q \square Q$ is
\[
\xymatrix@C1cm@R1cm{
(1,2s) \ar[r]^{(a_1,2s)} & (2,2s)\ar[d]^{(2,a_{2s-1})^*} & (3,2s) \ar[l]_{(b_3,2s)} \ar[r]^{(a_3,2s)} & \cdots \ar[r]^{(a_{2s-1},2s)\;\;\;\;} & (2s,2s)\ar[d]^{(2s,a_{2s-1})^*} \\
\vdots\ar[u]_{(1,a_{2s-1})} & \vdots\ar[d]^{(2,a_3)^*} & \vdots\ar[u]_{(3,a_{2s-1})} & \iddots & \vdots \ar[d]^{(2s,a_3)^*} \\
(1,3)\ar[u]_{(1,a_3)}\ar[d]^{(1,b_3)} & (2,3)\ar[l]_{(a_1,3)^*}\ar[r]^{(b_3,3)^*} & (3,3)\ar[d]^{(3,b_3)}\ar[u]_{(3,a_3)} & \ar[l]_{(a_3,3)^*} \cdots & \ar[l]_{(a_{2s-1},3)^*}(2s,3) \\
(1,2) \ar[r]^{(a_1,2)} & (2,2)\ar[d]^{(2,a_1)^*}\ar[u]_{(2,b_3)^*} & (3,2) \ar[l]_{(b_3,2)} \ar[r]^{(a_3,2)} & \cdots \ar[r]^{(a_{2s-1},2)} & (2s,2)\ar[d]^{(2s,a_1)^*}\ar[u]_{(2s,b_3)^*} \\
(1,1)\ar[u]_{(1,a_1)} & (2,1)\ar[l]_{(a_1,1)^*}\ar[r]^{(b_3,1)^*} & (3,1)\ar[u]_{(3,a_1)} & \ar[l]_{(a_3,1)^*} \cdots & \ar[l]_{(a_{2s-1},1)^*}(2s,1) \\
}
\]
The potential $W^{\square}_{Q,Q}$ is the sum of all cycles of length four that correspond to the small squares appearing above. 
Consider the automorphism $\widehat{K(Q \square Q)} \simeq \widehat{K(Q \square Q)}$ defined by multiplying all arrows of the form $(2k-1,a_{2i+1})$ and $(a_{2i+1},2k)$ by $-1$ and all other arrows by $1$. 
It sends $W^{\square}_{Q,Q}$ to the potential $W$ given as the sum of all clockwise oriented small squares minus the sum of all counterclockwise oriented small squares. 
Thus it suffices to show that $(Q \square Q,W)$ is selfinjective with Nakayama permutation induced by $\omega\times \omega$.

Set $\Lambda = \P(Q \square Q,W)$.
We say that a two paths $p$ and $p'$ in $Q \square Q$ are congruent (written $p \equiv p'$) if $p-p'\in \langle \partial_a W \;|\; a \in (Q \square Q)_1 \rangle$. 
Similarly we write $p \equiv 0$ if $p\in \langle \partial_a W \;|\; a \in (Q \square Q)_1 \rangle$. 
Since the relations $\partial_a W$ are all commutativity or zero relations, the set $B$ of nonzero congruence classes of paths forms a basis in $\Lambda$. 
For each $(i,j) \in (Q \square Q)_0$, let $B_{(i,j)}$ and ${}_{(i,j)}B$ be the subsets of $B$ corresponding to paths ending and starting at $(i,j)$ respectively. 
These sets form bases of $\Lambda e_{(i,j)}$ and $D(e_{(i,j)} \Lambda)$ respectively. 
We proceed to define a bijection ${}_{(i,j)}B \simeq B_{\sigma(i,j)}$ which induces an isomorphism $D(e_{(i,j)} \Lambda) \simeq \Lambda e_{\sigma(i,j)}$ thus completing the proof.

Consider the four cuts 
\[\begin{matrix}
R = \{(a_{2i+1},2k), (b_{2j+1},2k-1)^*\}_{i,j,k}, & 
U = \{(2k-1, a_{2i+1}), (2k,b_{2j+1})^*\}_{i,j,k}, \\ 
L = \{(a_{2i+1},2k-1)^*, (b_{2j+1},2k)\}_{i,j,k}, &
D = \{(2k, a_{2i+1})^*, (2k-1,b_{2j+1})\}_{i,j,k}.
\end{matrix}
\]
We define a $\Z^4$ grading on $(Q \square Q,W)$ by $\deg(a) = (g_{R}(a), g_{U}(a), g_{L}(a), g_{D}(a))$. 
Moreover, we associate to each path $p = a_1\cdots a_d$ in $Q$ an element $s(p) = s\in \{\pm 1\}^d$ by setting $s_i = 1$ if $a_i \in R\cup U$ and $s_i = -1$ if $a_i \in L\cup D$. 

At each vertex $x \in (Q \square Q)_0$ there are at most two arrows leaving that vertex. They both lie either in $R\cup L$ or $U \cup D$. 
The same is true for arrows entering $x$. Hence each path $p$ is determined by $s(p)$ together with either the starting point or end point of $p$. 
Moreover, if $s(p) = (u,v,-u)$ for some $u,v\in \{\pm 1\}$, then either $p\equiv 0$ or there is $p' \equiv p$ such that $s(p') = (-u,v,u)$. 

We conclude that for each path $p \not \equiv 0$ there is a path $p' \equiv p$ such that $s(p')_i \geq s(p')_{i+2}$ for all $i$. 
This condition determines $s(p')$ by $\deg(p)$ together with either the starting point or end point of $p$. 
Thus $p$ is determined up to $\equiv$ by $\deg(p)$ together with either the starting point or end point of $p$.

For each $(i,j) \in (Q \square Q)_0$ set $d_{ij} = (2s-i, 2s-j, i-1, j-1) \in \Z^4$. For $a, b \in \Z^4$ we write $a \leq b$ if $a_i \leq b_i$ for all $1\leq i \leq 4$. 
Let $p$ be a path starting at $(i,j) \in (Q \square Q)_0$. We claim that $p \equiv 0$ if and only if $\deg(p) \not \leq d_{ij}$. 
We only treat the case when $\deg(p)_1 > 2s-i$ and $i+j$ is even as the other cases are similar. 
By the arguments above there is a path $p' \equiv p$ such that $s(p')_{2k} = 1$ for $1\leq k \leq 2s-i$, $s(p')_{4s-2i+2} = -1$ and $s(p')_{4s-2i+4} = 1$. 
Then the arrows in position $4s-2i+2$, $4s-2i+3$, and $4s-2i+4$ in $p'$ form one of the zero relations in $\{\partial_a W \;|\; a \in (Q \square Q)_1\}$. 
Hence $p \equiv 0$. On the other hand, if none of the inequalities above hold, then it is impossible to reach any of the zero relations in $\{\partial_a W \;|\; a \in (Q \square Q)_1\}$.

Observe that for each $(i,j) \in (Q \square Q)_0$ there is a path $p_{ij}$ starting at $(i,j)$ of degree $d_{ij}$. 
In particular $p_{ij}$ ends at $(2s-i+1,2s-j+1)=(\omega\times \omega) (i,j)$. Now let $p$ be a path starting at $(i,j)$ satisfying $p \not \equiv 0$. 
Then $\deg(p) \leq \deg(p_{ij})$ and so there is a path $q$ ending at $(\omega\times \omega) (i,j)$ such that $pq \equiv p_{ij}$. 
Thus $\deg(q) = d_{ij} - \deg(p)$ and so $q$ is determined up to $\equiv$ by $p$. Similarly one can show that for each path $q \not \not \equiv 0$ 
ending at $(\omega\times \omega) (i,j)$ there is a path $p$ starting at $(i,j)$ unique up to $\equiv$ such that $pq \equiv p_{ij}$.

We define the bijection $\phi:{}_{(i,j)}B \rightarrow B_{\sigma(i,j)}$ by $\phi(p) = q$ where $pq \equiv p_{ij}$. 
It is easy to see that $\phi$ induces a $\Lambda$-module isomorphism $D(e_{(i,j)} \Lambda) \simeq \Lambda e_{\sigma(i,j)}$.
\end{proof}

\section{Covering theory of truncated quivers}\label{covering theory}

In the previous sections we have seen how mutation can be used to construct and organize selfinjective QPs. To construct $2$-representation-finite algebras we also need to consider their cuts. 
In this section we will be working in the setting when $Q$ is a quiver and $C$ is an arbitrary set of arrows in $Q$. 
In that case we call the pair $(Q,C)$ a \emph{truncated quiver}. Later we will specialize to the setting of QPs and their cuts.
Our results in this section generalize those in \cite{IO1} given for `quivers of type A'.

First recall some classical notions. The \emph{double} of $Q$ is the quiver $\overline{Q}$ defined by $\overline{Q}_0 = Q_0$ and $\overline{Q}_1 = Q_1 \coprod \{a^{-1} \: | \: a \in Q_1 \}$, with $s(a^{-1}) = e(a)$ and $e(a^{-1}) = s(a)$.
A \emph{walk} in $Q$ is a path in $\overline{Q}$. For each walk $p = a_1^{s_1}\cdots a_n^{s_n}$ in $Q$, its inverse walk is defined as $p^{-1} = a_n^{-s_n}\cdots a_1^{-s_1}$. 
Let $\sim$ be the equivalence relation on walks generated by $ap p^{-1}b \sim ab$ for all walks $a$, $p$ and $b$ such that
the end point of $a$ coincides with the starting point of $p$ and $b$. 
We call a walk \emph{reduced} if it is shortest in its equivalence class.

The grading $g_C$ on $Q$ extends to $\overline{Q}$ by $g_C(a^{-1}) = -g_C(a)$. This grading is invariant under $\sim$.

\subsection{Cut-slice correspondence}
We now give a general construction of a 
Galois covering of arbitrary quiver $Q$
with the Galois group $\Z$.

\begin{definition}
Let $(Q,C)$ be a truncated quiver.
\begin{itemize}
\item[(a)] We define a new quiver $\Z(Q,C)$ as follows:
\begin{itemize}
\item[$\bullet$] $\Z(Q,C)_0:=Q_0\times\Z$.
\item[$\bullet$] $\Z(Q,C)_1$ consists of arrows
$(a,\ell):(x,\ell)\to(y,\ell)$ for any arrow $a:x\to y$ in $Q_1\backslash C$
and $(a,\ell):(x,\ell)\to(y,\ell-1)$ for any arrow $a:x\to y$ in $C$.
\end{itemize}
\item[(b)] We define a morphism $\pi:\Z(Q,C)\to Q$ of quivers by
$\pi(x,\ell):=x$ for any $x\in Q_0$ and $\pi(a,\ell):=a$ for any $a\in Q_1$.
\item[(c)] We define an automorphism $\tau:\Z(Q,C)\to\Z(Q,C)$ by
$\tau(x,\ell):=(x,\ell+1)$ for any $x\in Q_0$ and
$\tau(a,\ell):=(a,\ell+1)$ for any $a\in Q_1$.
\end{itemize}
\end{definition}

\begin{example}\label{example of truncated quiver}
\begin{itemize}
\item[(a)] Let $Q$ be a quiver and $\overline{Q}$ be the double of $Q$.
Then the quiver $\Z(\overline{Q},Q_1)$ coincides with the translation quiver $\Z Q$
constructed by Riedtmann \cite{Rie}.
\item[(b)] Consider the following truncated quiver $(Q,C)$.
\[\xymatrix@R=.4cm@C=.4cm{
1\ar@<.5ex>^a[rr]\ar@{.>}@<-.5ex>_b[rr]&&2\ar[rr]^c&&3}\]
Then $\Z(Q,C)$ is given by the following quiver.
\[\xymatrix@R=.4cm@C=.4cm{
&&&\cdots&&&\\
(1,-1)\ar^{(a,-1)}[rrr]\ar^{(b,-1)}[urrr]&&&(2,-1)\ar^{(c,-1)}[rrr]&&&(3,-1)\\
(1,0)\ar^{(a,0)}[rrr]\ar^{(b,0)}[urrr]&&&(2,0)\ar^{(c,0)}[rrr]&&&(3,0)\\
(1,1)\ar^{(a,1)}[rrr]\ar^{(b,1)}[urrr]&&&(2,1)\ar^{(c,1)}[rrr]&&&(3,1)\\
&&&\cdots&&&
}\]
\end{itemize}
\end{example}

The following observation is immediate.

\begin{lemma}\label{covering}
For any $x\in Q_0$ and $\widetilde{x}\in\pi^{-1}(x)$, we have that $\pi$ induces
a bijection between arrows starting (respectively, ending)
at $\widetilde{x}$ in $\Z(Q,C)$ and those of $x$ in $Q$.
\end{lemma}

The following easy observation is useful.

\begin{lemma}\label{lifting path}
Let $p$ be a walk in $Q$ with $s(p)=x_0$ and $e(p)=y_0$. Let $\widetilde{x}_0\in\pi^{-1}(x_0)$. 
\begin{itemize}
\item[(a)] There exists a unique walk $\widetilde{p}$ in $\widetilde{Q}$ such that $s(\widetilde{p})=\widetilde{x}_0$ and $\pi(\widetilde{p})=p$.
\item[(b)] If $p$ is a cyclic walk, then we have $e(\widetilde{p})=\tau^{-g_C(p)}(s(\widetilde{p}))$.
\end{itemize}
\end{lemma}

We say that $\widetilde{p}$ is a \emph{lift} of $p$.

\begin{proof}
(a) This is immediate from Lemma \ref{covering}.

(b) Let $a:x\to y$ be an arrow in $Q$ and $\widetilde{a}:(x,\ell)\to (y,\ell')$ be a lift of $a$.
If $a\notin C$, then we have $\ell'=\ell$. If $a\in C$, then we have $\ell'=\ell-1$.
Thus the assertion follows.
\end{proof}

The quiver $\Z(Q,C)$ depends on $C$. So it is natural to ask when $\Z(Q,C)$ and $\Z(Q,C')$ are isomorphic for subsets $C$ and $C'$ of $Q_1$.
The key notion is the following:

\begin{definition}\label{compatible}
We say that two subsets $C,C' \subset Q_1$ are \emph{compatible}
if for any cyclic walk $p$ in $Q$ the equality $g_C(p) = g_{C'}(p)$ holds. In that case we write $C \sim C'$.
\end{definition}

\begin{proposition}\label{cover iso}
For subsets $C$ and $C'$ of $Q_1$, the following conditions are equivalent.
\begin{itemize}
\item[(a)] $C$ and $C'$ are compatible.
\item[(b)] There is an isomorphism $f:\Z(Q,C) \rightarrow \Z(Q,C')$ of quivers such that the following diagrams commute.
\[
\xymatrix@C0.3cm@R0.6cm{
\Z(Q,C) \ar[rr]^f \ar[dr]_{\pi} & & \Z(Q,C')\ar[dl]^{\pi} \\
& Q  &
}\ \ \ \ \ \xymatrix@C0.3cm@R0.6cm{
\Z(Q,C) \ar[rr]^f \ar[d]_{\tau} & & \Z(Q,C')\ar[d]^{\tau} \\
\Z(Q,C) \ar[rr]^f & & \Z(Q,C')
}\]
\end{itemize}
\end{proposition}

\begin{proof}
(b)$\Rightarrow$(a) Let $p$ be a cyclic walk in $Q$ starting at $x$.
Fix $\widetilde{x}\in\pi^{-1}(x)$ and $\widetilde{x}':=f(\widetilde{x})$.
Let $\widetilde{p}$ (respectively, $\widetilde{p}'$) be a lift of $p$ to $\Z(Q,C)$ (respectively, $\Z(Q,C')$)
starting at $x$ (respectively, $\widetilde{x}'$).
Then we have $\widetilde{p}'=f(\widetilde{p})$. Thus we have
\[\tau^{-g_{C'}(p)}(\widetilde{x}')\stackrel{{\rm Lem.}\ref{lifting path}}{=}e(\widetilde{p}')=f(e(\widetilde{p}))
\stackrel{{\rm Lem.}\ref{lifting path}}{=}f(\tau^{-g_C(p)}(\widetilde{x}))=\tau^{-g_C(p)}(f(\widetilde{x})) = \tau^{-g_C(p)}(\widetilde{x}').\]
Consequently we have $g_C(p)=g_{C'}(p)$.

(a)$\Rightarrow$(b) Without loss of generality, we can assume that $Q$ is connected.
Fix vertices $x\in\Z(Q,C)_0$ and $x'\in\Z(Q,C')_0$ such that $\pi(x)=\pi(x')$.

First we define a map $f:\Z(Q,C)_0\to\Z(Q,C')_0$.
For any $y\in\Z(Q,C)_0$, we take a walk $p:x\to y$ in $\Z(Q,C)$.
We take a lift $p':x'\to y'$ of $\pi(p)$ to $\Z(Q,C')$ starting at $x'$.
We shall show that $f(y):=y'$ is independent of choice of $p$.
Let $q:x\to y$ be another walk, and let $q':x'\to y''$ be a lift of $\pi(q)$
to $\Z(Q,C')$ starting at $x'$.
Then both $p^{-1}q:y\to y$ and $p'{}^{-1}q':y'\to y''$ are lifts of a cyclic walk $\pi(p^{-1}q):\pi(y)\to \pi(y)$ to $\Z(Q,C)$ and $\Z(Q,C')$ respectively.
By Lemma \ref{lifting path}, we have $y=\tau^{-g_{C}(\pi(p^{-1}q))}(y)$ and so $g_C(\pi(p^{-1}q))=0$. Thus we have
\[y''\stackrel{{\rm Lem.}\ref{lifting path}}{=}\tau^{-g_{C'}(\pi(p^{-1}q))}(y')
\stackrel{C\sim C'}{=}\tau^{-g_C(\pi(p^{-1}q))}(y')=y'.\]
Thus we have a well-defined map $f:\Z(Q,C)_0\to\Z(Q,C')_0$, which is clearly a bijection.

Next we define a map $f:\Z(Q,C)_1\to\Z(Q,C')_1$.
For any arrow $a:y\to z$ in $\Z(Q,C)$, we define $f(a)$ as a lift of $\pi(a)$
starting at $f(y)$ to $\Z(Q,C')$.
Clearly $f$ gives an isomorphism $\Z(Q,C)\to\Z(Q,C')$ of quivers.
\end{proof}

The next aim of this section is to give a one-to-one correspondence between
subsets of $Q_1$ which are compatible with $C$ and certain full subquivers of $\Z(Q,C)$ defined as follows:

\begin{definition}
A \emph{slice} is a full subquiver $S$ of $\Z(Q,C)$ satisfying the following conditions.
\begin{itemize}
\item Any $\tau$-orbit in $\Z(Q,C)_0$ contains precisely one vertex which belongs to $S$.
\item For any arrow $a:x\to y$ in $\Z(Q,C)$ with $x\in S_0$, we have $y\in S_0\cup\tau^{-1}S_0$.
\end{itemize}
\end{definition}

For example, the full subquiver of $\Z(Q,C)$ with the set $Q_0\times\{0\}$ of vertices is a slice.
If $S$ is a slice, then clearly so is $\tau^\ell S$ for any $\ell\in\Z$.

\begin{theorem}[{\bf cut-slice correspondence}]\label{cut slice}
Let $(Q,C)$ be a truncated quiver such that $Q$ is connected. Then we have a bijection
\[
\{S \subset \Z(Q,C) \:|\: S \mbox{ is a slice}\}/\tau^{\Z} \rightarrow \{C' \subset Q_1 \:|\: C' \sim C \}
\]
induced by $S \mapsto C_S:=Q_1 \setminus \pi(S_1)$.
Moreover $\pi$ induces an isomorphism $S\to Q_{C_S}$ of quivers.
\end{theorem}

\begin{proof}
(i) We shall show that $C_S$ is compatible with $C$.

We only have to construct an isomorphism $\Z(Q,C_S)\to\Z(Q,C)$ of quivers satisfying the conditions in Proposition \ref{cover iso}(b).
For each $x\in Q_0$, we denote by $\widetilde{x}$ the unique lift of $x$ in $S_0$.
Define a bijection $f:\Z(Q,C_S)_0\to\Z(Q,C)_0$ by $f(x,\ell):=\tau^\ell\widetilde{x}$.

We shall show that $f$ gives a quiver isomorphism.
Arrows of $\Z(Q,C_S)$ has the form $(a,\ell):(x,\ell)\to(y,\ell)$ with $a\in Q_{C_S}$
or $(a,\ell):(x,\ell)\to(y,\ell-1)$ with $a\in C_S$.
If $a\in Q_{C_S}$, then we have a lift $\widetilde{a}:\widetilde{x}\to\widetilde{y}$ in $S_1$.
Thus we can define $f(a,\ell):=\tau^\ell\widetilde{a}:(x,\ell)\to(y,\ell)$.
If $a\in C_S$, then we have a lift $\widetilde{a}:\widetilde{x}\to\tau^{-1}\widetilde{y}$ in $S_1$. 
Thus we can define $f(a,\ell):=\tau^\ell\widetilde{a}:(x,\ell)\to(y,\ell-1)$.

(ii) We shall show that the map $C\mapsto C_S$ is surjective.

Let $C'$ be a subset of $Q_1$ which is compatible with $C$.
Then we have an isomorphism $f:\Z(Q,C)\to\Z(Q,C')$ given in Proposition \ref{cover iso}.
Let $S'$ be the full subquiver of $\Z(Q,C')$ with the set of vertices $Q_0\times\{0\}$.
Then $S'$ is a slice of $\Z(Q,C')$.
Since $f$ is an isomorphism of quivers, we have that $S:=f^{-1}(S')$ is a slice of $\Z(Q,C)$.
Clearly we have $C_S=C'$.

(iii) We shall show that the map $S\mapsto C_S$ is injective.

Let $S$ and $S'$ be slices in $\Z(Q,C)$ such that $C_S=C_{S'}$.
Without loss of generality, we can assume $S_0\cap S'_0\neq\emptyset$
by replacing $S$ be $\tau^\ell S$ for some $\ell\in\Z$.

Let $X:=\pi(S_0\cap S'_0)$. Then $X$ is a non-empty subset of $Q_0$.
For any arrow $a:x\to y$ in $Q$, we will show that $x$ belongs to $X$ if and only if $y$ belongs to $X$.
Then we have $X=Q_0$ and $S=S'$ since $Q$ is connected.

Assume $x\in X$. Let $\widetilde{x}\in\pi^{-1}(x)$ be a vertex in $S_0\cap S'_0$ and 
$\widetilde{a}$ be a lift of $a$ starting at $\widetilde{x}$.
If $a\notin C_S=C_{S'}$, then $\widetilde{a}$ belongs to $S_1\cap S'_1$.
Thus $e(\widetilde{a})$ belongs to both $S_0$ and $S'_0$, and so $y=e(a)=\pi(e(\widetilde{a}))$
belongs to $X$.
If $a\in C_S=C_{S'}$, then $\widetilde{a}$ belong to neither $S_1$ nor $S'_1$.
Since $S$ and $S'$ are slices, we have that $e(\widetilde{a})$ belongs to both $\tau^{-1}S_0$ and $\tau^{-1}S'_0$,
and so $y=e(a)=\pi(e(\widetilde{a}))$ belongs to $X$.
Thus we have $y\in X$ in both cases.

Similarly one can check that $y\in X$ implies $x\in X$, so we are done.
\end{proof}

\begin{example}
Consider the truncated quiver $(Q,C)$ and its covering $\Z(Q,C)$
in Example \ref{example of truncated quiver}.
Then $C$ has the two compatible sets and $\Z(Q,C)$ has two slices
up to $\tau^\Z$.
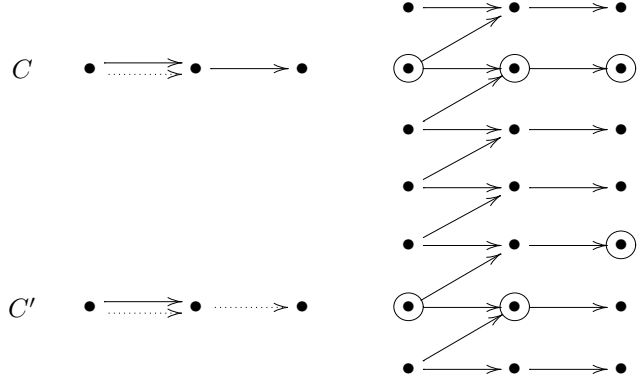
\[\xymatrix@R=.4cm@C=.4cm{
&&&&&&&
  \bullet\ar[rr]&&\bullet\ar[rr]&&\bullet\\
C&\bullet\ar@<.5ex>[rr]\ar@{.>}@<-.5ex>[rr]&&\bullet\ar[rr]&&\bullet&&
  *+[o][F]{\bullet}\ar[rr]\ar[urr]&&*+[o][F]{\bullet}\ar[rr]&&*+[o][F]{\bullet}\\
&&&&&&&
  \bullet\ar[rr]\ar[urr]&&\bullet\ar[rr]&&\bullet\\
&&&&&&&
  \bullet\ar[rr]\ar[urr]&&\bullet\ar[rr]&&\bullet\\
&&&&&&&
  \bullet\ar[rr]\ar[urr]&&\bullet\ar[rr]&&*+[o][F]{\bullet}\\
C'&\bullet\ar@<.5ex>[rr]\ar@{.>}@<-.5ex>[rr]&&\bullet\ar@{.>}[rr]&&\bullet&&
  *+[o][F]{\bullet}\ar[rr]\ar[urr]&&*+[o][F]{\bullet}\ar[rr]&&\bullet\\
&&&&&&&
  \bullet\ar[rr]\ar[urr]&&\bullet\ar[rr]&&\bullet
}\]
\end{example}

\subsection{Transitivity of cut-mutation and slice-mutation}

For a given subset $C$ of $Q_1$, we have operation to construct subsets of $Q_1$ which are
compatible with $C$.

\begin{definition}
Let $C$ be a subset of $Q_1$.
\begin{itemize}
\item[(a)] We say that a vertex $x$ of $Q$ is a \emph{strict source} of $(Q,C)$
if all arrows ending at $x$ belong to $C$ and all arrows starting at $x$ do not belong to $C$. 
\item[(b)] For a strict source $x$ of $(Q,C)$, we define the subset $\mu^+_x(C)$ of $Q_1$ by removing all
arrows in $Q$ ending at $x$ from $C$ and adding all arrows in $Q$ starting at $x$ to $C$.
\item[(c)] Dually we define a \emph{strict sink} and $\mu^-_x(C)$.
\end{itemize}
We call these operations \emph{cut-mutation}.
\end{definition}

We have the following easy observation.

\begin{lemma}\label{compatible mutation}
We have $\mu^+_x(C)\sim C$ (respectively, $\mu^-_x(C)\sim C$).
Moreover $x$ is a strict sink (respectively, strict source) of $(Q,\mu^+_x(C))$ (respectively, $(Q,\mu^-_x(C))$) and
we have $\mu^-_x(\mu^+_x(C))=C$ (respectively, $\mu^+_x(\mu^-_x(C))=C$).
\end{lemma}

\begin{proof}
Clearly we have $g_C(p)=g_{\mu^{\pm}_x(C)}(p)$ for any cyclic walk $p$ in $Q$.
Thus we have $\mu^{\pm}_x(C)\sim C$.
The other assertions are clear from the definition.
\end{proof}

It is natural to ask whether successive cut-mutation acts transitively
on the set of compatible subsets.

We need the following operation for slices.

\begin{definition}
Let $S$ be a slice in $\Z(Q,C)$.
\begin{itemize}
\item[(a)] We say that a vertex $x$ of $S$ is a \emph{strict source} if all arrows in $\Z(Q,C)$ ending at $x$ do not belong to $S$
and  all arrows in $\Z(Q,C)$ starting at $x$ belong to $S$. 

\item[(b)] For a strict source $x$ of $S$, we define the full subquiver $\mu^+_x(S)$ of $\Z(Q,C)$ by $\mu^+_x(S)_0=(S_0\backslash\{x\})\cup\{\tau^{-1}x\}$.
\item[(c)] Dually we define a \emph{strict sink} and $\mu^-_x(S)$.
\end{itemize}
We call these operations \emph{slice-mutation}.
\end{definition}

The following observation is clear.

\begin{lemma}\label{slice-mutation is a slice}
We have that $\mu^+_x(S)$ (respectively, $\mu^-_x(S)$) is again a slice of $\Z(Q,C)$.
Moreover $\tau^{-1}x$ is a strict sink (respectively, $\tau x$ is a strict source) of $\mu^+_x(S)$ (respectively, $\mu^-_x(S)$) and
we have $\mu^-_{\tau^{-1}x}(\mu^+_x(S))=S$ (respectively, $\mu^+_{\tau x}(\mu^-_x(S))=S$).
\end{lemma}

\begin{proof}
Let $a:y\to z$ be an arrow in $\Z(Q,C)$ such that either $y\in\mu^+_x(S)_0$.

Assume $y=\tau^{-1}x$. Then we have an arrow $\tau a:x\to\tau z$ in $\Z(Q,C)$.
Since $x$ is a strict source of $S$, we have $\tau z\in S_0\backslash\{x\}\subset\mu^+_x(S)_0$. Thus we have $z\in\tau^{-1}\mu^+_x(S)_0$.

Assume $y\neq\tau^{-1}x$. Then $y\in S_0$, and we have $z\in S_0\cup\tau^{-1}S_0$ since $S$ is a slice.
Since $x$ is a strict source of $S$, we have $z\neq x$. Thus we have
$z\in (S_0\backslash\{x\})\cup\tau^{-1}S_0\subset\mu^+_x(S)_0\cup\tau^{-1}\mu^+_x(S)_0$.

Consequently we have that $\mu^+_x(S)$ is a slice.
The other assertions can be checked easily.
\end{proof}

Under the bijection in Theorem \ref{cut slice}, slice-mutation corresponds to cut-mutation by the following observation.

\begin{proposition}\label{mutation cut slice}
Let $S$ and $Q_{C_S}$ be in Theorem \ref{cut slice}. Let $x\in S_0$.
\begin{itemize}
\item[(a)] $x$ is a strict source (respectively, strict sink) of $S$ if and only if so is $\pi(x)$ in $(Q,C_S)$.
\item[(b)] In this case we have $\mu^+_{\pi(x)}(C_S)=C_{\mu^+_x(S)}$ (respectively, $\mu^-_{\pi(x)}(C_S)=C_{\mu^-_x(S)}$).
\end{itemize}
\end{proposition}

\begin{proof}
(a) Since $\pi:S\to Q_{C_S}$ is an isomorphism of quives, this is clear from Lemma \ref{covering}.

(b) The difference of $C_{\mu^+_x(S)}$ and $C_S$ are arrows starting or ending at $\pi(x)$.
Since $\tau^{-1}x$ is a strict sink of $\mu^+_x(S)$ by Lemma \ref{slice-mutation is a slice}, so is $\pi(x)$ in $C_{\mu^+_x(S)}$ by (a).
Thus we have $C_{\mu^+_x(S)}=\mu^+_{\pi(x)}(C_S)$.
\end{proof}

In the rest of this subsection, we study transitivity of successive cut-mutations (respectively, slice-mutations).

\begin{definition}
Let $(Q,C)$ be a truncated quiver. We say that $(Q,C)$
\begin{itemize}
\item[(a)] has \emph{enough compatibles} if $Q_1 = \bigcup_{C'\sim C} C'$
\item[(b)] is \emph{sufficiently cyclic} if for each $a \in Q_1$ there is a cycle $p$ containing $a$ satisfying $g_C(p) \leq 1$.
\end{itemize}
If $(Q,C)$ has any of these properties, then so does $(Q,C')$ for any $C' \sim C$.
\end{definition}


We have a useful criterion for enough compatibility.

We call a numbering $x_1,\ldots,x_N$ of vertices of $Q$
a \emph{$C$-source sequence} (respectively, \emph{strict $C$-source sequence}) 
if $x_{i+1}$ is a source (respectively, strict source) in $Q_{\mu^+_{x_i}\circ\cdots\circ\mu^+_{x_1}(C)}$.
Dually we define a \emph{$C$-sink sequence} (respectively, \emph{strict $C$-sink sequence}).

\begin{proposition}\label{acyclic}
For a truncated quiver $(Q,C)$, the following conditions are equivalent.
\begin{itemize}
\item[(a)] $C$ has enough compatibles.
\item[(a')] For any cycle $p$ in $Q$, at least one arrow in $p$ is contained in some $C'\sim C$.
\item[(b)] $Q_{C'}$ is an acyclic quiver for any $C'\sim C$.
\item[(b')] $Q_C$ is an acyclic quiver.
\item[(c)] There exists a $C'$-source sequence for any $C'\sim C$.
\item[(c')] There exists a $C$-source sequence.
\item[(d)] There exists a $C'$-sink sequence for any $C'\sim C$.
\item[(d')] There exists a $C$-sink sequence.
\end{itemize}
If $(Q,C)$ is sufficiently cyclic, then we can replace the above `source' (respectively, `sink')
by `strict source' (respectively, `strict sink').
\end{proposition}

\begin{proof}
(a)$\Rightarrow$(a'), (b)$\Rightarrow$(b') and (c)$\Rightarrow$(c') Clear.

(a')$\Rightarrow$(b) Assume that $Q_{C'}$ has a cycle $p$.
By (a'), there exists a cut $C''$ of $Q$ such that
some arrow in $p$ is contained in $C''$.
Since $C'\sim C\sim C''$ we have $0=g_{C'}(p)=g_{C''}(p)>0$, a contradiction.

(b)$\Rightarrow$(c) and (b')$\Rightarrow$(c')
Assume that $x_1,\ldots,x_{i-1}$ are defined.
Since $Q_C$ is acyclic, so is the quiver
$Q_C\backslash\{x_1,\ldots,x_{i-1}\}$.
We define $x_i$ as a source of this new quiver.
It is easily checked that $x_1,\ldots,x_N$ is a $C$-source sequence.

(c')$\Rightarrow$(a) Assume that there is an arrow $a:x\to y$ that does not belong to any $C'\sim C$.
Then $x$ is not a source of $Q_{C'}$ for any $C'\sim C$ since $a:x\to y$ belongs to $Q_{C'}$.
Take a $C$-source sequence $x_1,\ldots,x_N$. 
Since $x$ is not a source of $Q_{\mu^+_{x_i}\circ\cdots\circ\mu^+_{x_1}(C)}$ for any $i$,
we have that $x$ is none of $x_1,\ldots,x_N$, a contradiction.

The assertions for the sufficiently cyclic case follows from the following Lemma.
\end{proof}

\begin{lemma}\label{mutable}
Assume that $(Q,C)$ is sufficiently cyclic. Then each source (respectively, sink) in $Q_C$ is strict.
\end{lemma}

\begin{proof}
Let $x$ be a source (respectively, sink) in $Q_C$.
We only have to show that any arrow $a \in Q_1$ starting (respectively, ending) at $x$ does not belong to $C$.
Since $(Q,C)$ is sufficiently cyclic there is a cycle $p$ in $Q$ containing $a$ such that  $g_{C}(p) \leq 1$. Since $p$ is a cycle it must contain an arrow $b$ ending (respectively, starting) at $x$. 
Thus $b \in C$ and $g_{C}(p) \geq g_{C}(a) + g_{C}(b) = g_{C}(a) + 1$. Hence $g_{C}(a) =0$ and $a \not \in C$.
\end{proof}

Now we are ready to prove the following main result in this subsection.

\begin{theorem}[{\bf transitivity}]\label{transitive mutation cut}
Let $(Q,C)$ be a truncated quiver which has enough compatibles and is sufficiently cyclic.
\begin{itemize}
\item[(a)] The set $\{C' \subset Q_1 \:|\: C' \sim C \}$ is transitive under successive cut-mutations.
\item[(b)] 
The set of all slices in $\Z(Q,C)$ is transitive under successive slice-mutations.
\end{itemize}
\end{theorem}

\begin{proof}
By Proposition \ref{mutation cut slice} it suffices to show (b).

Let $S$ be a slice in $\Z(Q,C)$. Our aim is to show that $S$ is slice-mutation equivalent to $S^0 := Q_C \times \{0\}$, 
which completes the proof. Set $C' := Q_1\setminus \pi(S_1)$. Since $(Q,C)$ has enough compatibles $Q_{C'}$ is acyclic and has both a strict $C'$-source sequence and a strict $C'$-sink sequence by Proposition \ref{acyclic}. 
Thus by Proposition \ref{mutation cut slice}, $S$ is acyclic and slice-mutation equivalent to $\tau^lS$ for any $l \in \Z$. In particular, we may assume that $S_0 \subset Q_0\times\Z_{\geq 0}$. Observe also that each sink and source in $S$ is strict by Proposition \ref{mutation cut slice} and Lemma \ref{mutable}.

Let $x \in Q_0$ and define the height at $x$ by $h_x(S) = l$ where $(x,l) \in S_0$. 
Moreover, define the volume under $S$ to be $V(S) = \sum_{x\in Q_0} h_x(S)$. 
Let $a : (x,l) \rightarrow (y,m)$ be an arrow in $S$. Since $m \in \{l, l-1\}$, 
we have $h_x(S) = l \geq m = h_y(S)$. Because $S$ is acyclic it follows that $h_x(S)$ 
takes its maximal value at some $x$ such that $(x,h_x(S))$ is a source in $S$. 
In particular if $V(S) >0$, there is a source $(x,l)$ in $S$ with $l >0$. 
Let $S' = \mu^+_x(S)$. Then $h_y(S') = h_y(S)$ for all $y \not = x$ and $h_x(S') = h_x(S) -1$. 
Thus $h_y(S') \geq 0$ for all $y \in Q_0$ and $V(S') = V(S) -1 < V(S)$. Hence we may assume that $V(S) =0$. 
But this implies $S_0 = S^0_0$ and thus $S = S^0$.
\end{proof}

\section{Derived equivalence of truncated Jacobian algebras}\label{derived equivalence}

In this section we study the relationship between truncated Jacobian algebras of a fixed QP by algebraic cuts.
Let us start with general observations on cut-mutation.

The following observations are easy.

\begin{lemma}
Let $(Q,W)$ be a QP with a cut $C$.
\begin{itemize}
\item[(a)] If $C'$ is a subset of $Q_1$ such that $C\sim C'$, then $C'$ is also a cut of $(Q,W)$.
\item[(b)] If $x$ is a strict source (respectively, strict sink) of $(Q,C)$, then $\mu^+_x(C)$ (respectively, $\mu^-_x(C)$) is again a cut of $(Q,W)$.
\end{itemize}
\end{lemma}

\begin{proof}
(a) This is clear. 

(b) By Lemma \ref{compatible mutation}, we have $\mu^\pm_x(C)\sim C$. Thus they are cuts of $(Q,W)$ by (a).
\end{proof}

The following result is an immediate consequence of Theorem \ref{transitive mutation cut}.

\begin{theorem}\label{independence0}
Let $(Q,W)$ be a sufficiently cyclic QP with a cut $C$ having enough compatibles. Then
the set of all cuts of $(Q,W)$ compatible with $C$ is transitive under successive cut-mutations.
\end{theorem}

Thus it is natural to ask when a cut of a QP is sufficiently cyclic and has enough compatibles.

We have the following sufficient conditions for sufficient cyclicity.

\begin{proposition} \label{selfinj sc}
Let $(Q,W)$ be a QP with a cut $C$. Then:
\begin{itemize}
\item[(a)] If any arrow in $Q$ appears in cycles in $W$, then $(Q,C)$ is sufficiently cyclic.
\item[(b)] If $(Q,W)$ is selfinjective, then $(Q,C)$ is sufficiently cyclic.
\end{itemize}
\end{proposition}

\begin{proof}
(a) For any $a\in Q_1$, take a cycle $p$ in $W$. Then $g_C(p)=1$.

(b) This is immediate from (a) and Proposition \ref{sc}.
\end{proof}

We separate enough compatibility into the following two conditions.

\begin{definition}
Let $(Q,W)$ be a QP. We say that $(Q,W)$
\begin{itemize}
\item[(a)] is \emph{fully compatible} if all cuts are compatible with each other.
\item[(b)] has \emph{enough cuts} if each arrow in $Q$ is contained in a cut.
\end{itemize}
\end{definition}

If $(Q,W)$ is fully compatible and has enough cuts, then for each cut $C$ in $(Q,W)$ the truncated quiver $(Q,C)$ has enough compatibles. 
Thus we obtain the following result.

\begin{corollary}\label{independence}
Let $(Q,W)$ be a selfinjective, fully compatible QP with enough cuts. Then
the set of all cuts of $(Q,W)$ is transitive under successive cut-mutations.
\end{corollary}

\begin{proof}
Since $(Q,W)$ is selfinjective Proposition \ref{selfinj sc} implies that $(Q,C)$ is sufficiently cyclic for any cut $C$. Thus the assertion follow from Theorem \ref{independence0}.
\end{proof}

In Section \ref{canvas} we give a sufficient condition for QPs to be fully compatible.

\subsection{Cut-mutation and 2-APR tilting}

We say that two algebras $A$ and $A'$ of global dimension at most two are \emph{cluster equivalent} if generalized cluster categories $\C_A$ and $\C_{A'}$ are triangle equivalent.
This is the case if $A$ and $A'$ are derived equivalent.

\begin{proposition}\label{cluster equivalence}
Let $(Q,W)$ be a QP.
Then all truncated Jacobian algebras $\P(Q,W)_C$ given by algebraic cuts $C$ are cluster equivalent.
\end{proposition}

\begin{proof} 
By Proposition \ref{derived preproj}, the derived 3-preprojective DG algebra of $\P(Q,W)_C$ is quasi-isomorphic to $\Gamma(Q,W)$.
Thus the assertion follows.
\end{proof}

In \cite{AO} Amiot-Oppermann studied when two cluster equivalent algebras are derived equivalent (see also \cite{Am2}). 
Thus it is natural to ask when the algebras $\P(Q,W)_C$ appearing above are derived equivalent to each other. 
In this section we give a sufficient condition on $(Q,W)$.


The main point of cut-mutation is its relation to 2-APR tilting modules \cite{IO1} defined as follows:
Let $x$ be a sink and $S$ the corresponding simple projective $A$-module.
If $\Ext_A^1(DA,S)=0$, then we have a cotilting $A$-module
\[T:=\tau^-\Omega^-S\oplus(A/S)\]
which we call a \emph{2-APR tilting} $A$-module.
Dually a \emph{2-APR cotilting module} $T$ is defined.
In both cases we call $\End_A(T)$ a \emph{2-APR tilt} of $A$.

\begin{theorem}\label{2-APR tilting}
Let $(Q,W)$ be a QP and $C$ an algebraic cut.
Let $T$ be a 2-APR tilting (respectively, cotilting) $\P(Q,W)_C$-module associated with a source (respectively, sink) $x$ of the quiver $Q_C$.
Then $\End_{\P(Q,W)_C}(T)$ is isomorphic to $\P(Q,W)_{\mu^+_x(C)}$ (respectively, $\P(Q,W)_{\mu^-_x(C)}$).
\end{theorem}

\begin{proof}
Let $A:=\P(Q,W)_C$ and $B:=\End_A(T)$.
The change of quivers with relations via 2-APR tilt was explicitly given in \cite[Theorem 3.11]{IO1}.
In particular we have that the QPs $(Q_A,W_A)$ and $(Q_B,W_B)$ are isomorphic,
and the isomorphism $Q_A\simeq Q_B$ induces a bijection between $\mu^+_x(C_A)$ and $C_B$.
Consequently we have
$B\simeq\P(Q_B,W_B)_{C_B}\simeq\P(Q_A,W_A)_{\mu^+_x(C_A)}$.
\end{proof}

Combining Corollary \ref{independence} with Theorem \ref{2-APR tilting} obtain the main result of this section.

\begin{theorem}\label{transitive}
Let $(Q,W)$ be a selfinjective, fully compatible QP with enough cuts.
Then all truncated Jacobian algebras of $(Q,W)$ are iterated 2-APR tilts of each other.
In particular they are derived equivalent.
\end{theorem}

\begin{proof}
The former assertion is an immediate consequence of Theorem \ref{2-APR tilting} since any source of 2-representation-finite algebras admits a 2-APR tilting module.
The latter assertion is a direct consequence of the former one.
\end{proof}

To apply Theorem \ref{transitive} it is important to have sufficient conditions for a QP to be fully compatible and have enough cuts. 
In Section \ref{canvas} we provide a sufficient condition for full compatibility.

\begin{example}
Let $Q = \bullet \leftarrow \bullet \rightarrow \bullet$. By Proposition \ref{tensor} and Theorem \ref{square}
the QPs $(Q \tildetensor Q,W^{\tildetensor}_{Q,Q})$ and $(Q \square Q,W^{\square}_{Q,Q})$ are selfinjective.
Their cut-mutations are displayed in Figure \ref{Cut mutation lattice}. 
For each of the two QPs, all cuts are connected by successive cut-mutation. 
By Theorem \ref{2-APR tilting} the corresponding $2$-representation finite algebras are derived equivalent.
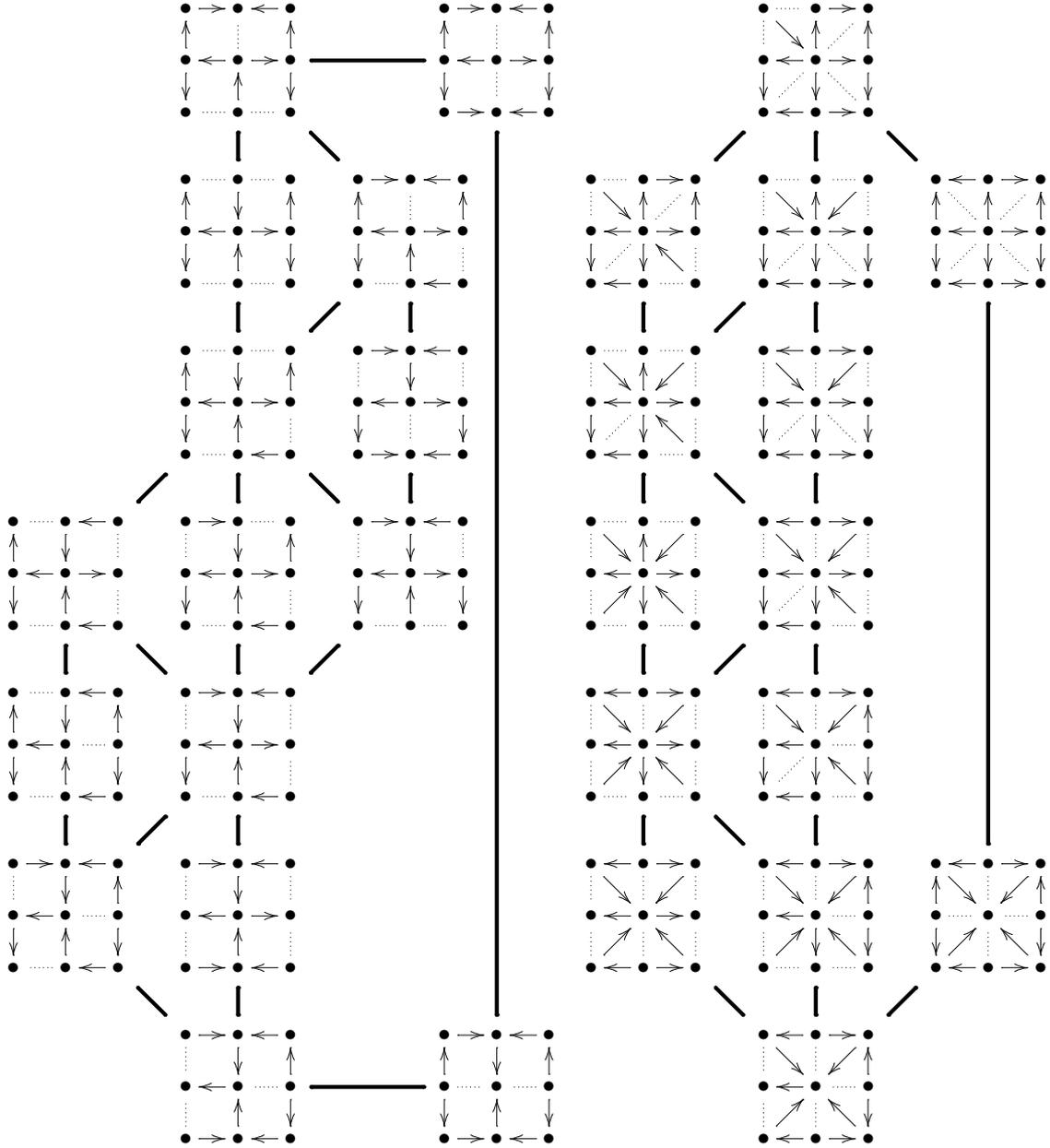
\begin{figure}[h]\label{Cut mutation lattice}
\[
\begin{xy} 0;<1.4pt,0pt>:<0pt,-1.4pt>:: 
(50,150) *+{
\begin{xy} 0;<0.17pt,0pt>:<0pt,-0.17pt>:: 
(0,0) *+{\bullet} ="0",
(125,0) *+{\bullet} ="1",
(250,0) *+{\bullet} ="2",
(0,125) *+{\bullet} ="3",
(125,125) *+{\bullet} ="4",
(250,125) *+{\bullet} ="5",
(0,250) *+{\bullet} ="6",
(125,250) *+{\bullet} ="7",
(250,250) *+{\bullet} ="8",
"0", {\ar"1"},
"3", {\cut"0"},
"2", {\cut"1"},
"1", {\ar"4"},
"5", {\ar"2"},
"4", {\ar"3"},
"3", {\ar"6"},
"4", {\ar"5"},
"7", {\ar"4"},
"5", {\cut"8"},
"6", {\cut"7"},
"8", {\ar"7"},
\end{xy}
} ="0",
(50,100) *+{
\begin{xy} 0;<0.17pt,0pt>:<0pt,-0.17pt>:: 
(0,0) *+{\bullet} ="0",
(125,0) *+{\bullet} ="1",
(250,0) *+{\bullet} ="2",
(0,125) *+{\bullet} ="3",
(125,125) *+{\bullet} ="4",
(250,125) *+{\bullet} ="5",
(0,250) *+{\bullet} ="6",
(125,250) *+{\bullet} ="7",
(250,250) *+{\bullet} ="8",
"0", {\cut"1"},
"3", {\ar"0"},
"2", {\cut"1"},
"1", {\ar"4"},
"5", {\ar"2"},
"4", {\ar"3"},
"3", {\ar"6"},
"4", {\ar"5"},
"7", {\ar"4"},
"5", {\cut"8"},
"6", {\cut"7"},
"8", {\ar"7"},
\end{xy}
} ="1",
(100,50) *+{
\begin{xy} 0;<0.17pt,0pt>:<0pt,-0.17pt>:: 
(0,0) *+{\bullet} ="0",
(125,0) *+{\bullet} ="1",
(250,0) *+{\bullet} ="2",
(0,125) *+{\bullet} ="3",
(125,125) *+{\bullet} ="4",
(250,125) *+{\bullet} ="5",
(0,250) *+{\bullet} ="6",
(125,250) *+{\bullet} ="7",
(250,250) *+{\bullet} ="8",
"0", {\ar"1"},
"3", {\ar"0"},
"2", {\ar"1"},
"1", {\cut"4"},
"5", {\ar"2"},
"4", {\ar"3"},
"3", {\ar"6"},
"4", {\ar"5"},
"7", {\ar"4"},
"5", {\cut"8"},
"6", {\cut"7"},
"8", {\ar"7"},
\end{xy}
} ="2",
(50,50) *+{
\begin{xy} 0;<0.17pt,0pt>:<0pt,-0.17pt>:: 
(0,0) *+{\bullet} ="0",
(125,0) *+{\bullet} ="1",
(250,0) *+{\bullet} ="2",
(0,125) *+{\bullet} ="3",
(125,125) *+{\bullet} ="4",
(250,125) *+{\bullet} ="5",
(0,250) *+{\bullet} ="6",
(125,250) *+{\bullet} ="7",
(250,250) *+{\bullet} ="8",
"0", {\cut"1"},
"3", {\ar"0"},
"2", {\cut"1"},
"1", {\ar"4"},
"5", {\ar"2"},
"4", {\ar"3"},
"3", {\ar"6"},
"4", {\ar"5"},
"7", {\ar"4"},
"5", {\ar"8"},
"6", {\cut"7"},
"8", {\cut"7"},
\end{xy}
} ="3",
(50,0) *+{
\begin{xy} 0;<0.17pt,0pt>:<0pt,-0.17pt>:: 
(0,0) *+{\bullet} ="0",
(125,0) *+{\bullet} ="1",
(250,0) *+{\bullet} ="2",
(0,125) *+{\bullet} ="3",
(125,125) *+{\bullet} ="4",
(250,125) *+{\bullet} ="5",
(0,250) *+{\bullet} ="6",
(125,250) *+{\bullet} ="7",
(250,250) *+{\bullet} ="8",
"0", {\ar"1"},
"3", {\ar"0"},
"2", {\ar"1"},
"1", {\cut"4"},
"5", {\ar"2"},
"4", {\ar"3"},
"3", {\ar"6"},
"4", {\ar"5"},
"7", {\ar"4"},
"5", {\ar"8"},
"6", {\cut"7"},
"8", {\cut"7"},
\end{xy}
} ="4",
(125,0) *+{
\begin{xy} 0;<0.17pt,0pt>:<0pt,-0.17pt>:: 
(0,0) *+{\bullet} ="0",
(125,0) *+{\bullet} ="1",
(250,0) *+{\bullet} ="2",
(0,125) *+{\bullet} ="3",
(125,125) *+{\bullet} ="4",
(250,125) *+{\bullet} ="5",
(0,250) *+{\bullet} ="6",
(125,250) *+{\bullet} ="7",
(250,250) *+{\bullet} ="8",
"0", {\ar"1"},
"3", {\ar"0"},
"2", {\ar"1"},
"1", {\cut"4"},
"5", {\ar"2"},
"4", {\ar"3"},
"3", {\ar"6"},
"4", {\ar"5"},
"7", {\cut"4"},
"5", {\ar"8"},
"6", {\ar"7"},
"8", {\ar"7"},
\end{xy}
} ="5",
(125,300) *+{
\begin{xy} 0;<0.17pt,0pt>:<0pt,-0.17pt>:: 
(0,0) *+{\bullet} ="0",
(125,0) *+{\bullet} ="1",
(250,0) *+{\bullet} ="2",
(0,125) *+{\bullet} ="3",
(125,125) *+{\bullet} ="4",
(250,125) *+{\bullet} ="5",
(0,250) *+{\bullet} ="6",
(125,250) *+{\bullet} ="7",
(250,250) *+{\bullet} ="8",
"0", {\ar"1"},
"3", {\ar"0"},
"2", {\ar"1"},
"1", {\ar"4"},
"5", {\ar"2"},
"4", {\cut"3"},
"3", {\ar"6"},
"4", {\cut"5"},
"7", {\ar"4"},
"5", {\ar"8"},
"6", {\ar"7"},
"8", {\ar"7"},
\end{xy}
} ="6",
(50,300) *+{
\begin{xy} 0;<0.17pt,0pt>:<0pt,-0.17pt>:: 
(0,0) *+{\bullet} ="0",
(125,0) *+{\bullet} ="1",
(250,0) *+{\bullet} ="2",
(0,125) *+{\bullet} ="3",
(125,125) *+{\bullet} ="4",
(250,125) *+{\bullet} ="5",
(0,250) *+{\bullet} ="6",
(125,250) *+{\bullet} ="7",
(250,250) *+{\bullet} ="8",
"0", {\ar"1"},
"3", {\cut"0"},
"2", {\ar"1"},
"1", {\ar"4"},
"5", {\ar"2"},
"4", {\ar"3"},
"3", {\cut"6"},
"4", {\cut"5"},
"7", {\ar"4"},
"5", {\ar"8"},
"6", {\ar"7"},
"8", {\ar"7"},
\end{xy}
} ="7",
(0,250) *+{
\begin{xy} 0;<0.17pt,0pt>:<0pt,-0.17pt>:: 
(0,0) *+{\bullet} ="0",
(125,0) *+{\bullet} ="1",
(250,0) *+{\bullet} ="2",
(0,125) *+{\bullet} ="3",
(125,125) *+{\bullet} ="4",
(250,125) *+{\bullet} ="5",
(0,250) *+{\bullet} ="6",
(125,250) *+{\bullet} ="7",
(250,250) *+{\bullet} ="8",
"0", {\ar"1"},
"3", {\cut"0"},
"2", {\ar"1"},
"1", {\ar"4"},
"5", {\ar"2"},
"4", {\ar"3"},
"3", {\ar"6"},
"4", {\cut"5"},
"7", {\ar"4"},
"5", {\ar"8"},
"6", {\cut"7"},
"8", {\ar"7"},
\end{xy}
} ="8",
(50,250) *+{
\begin{xy} 0;<0.17pt,0pt>:<0pt,-0.17pt>:: 
(0,0) *+{\bullet} ="0",
(125,0) *+{\bullet} ="1",
(250,0) *+{\bullet} ="2",
(0,125) *+{\bullet} ="3",
(125,125) *+{\bullet} ="4",
(250,125) *+{\bullet} ="5",
(0,250) *+{\bullet} ="6",
(125,250) *+{\bullet} ="7",
(250,250) *+{\bullet} ="8",
"0", {\ar"1"},
"3", {\cut"0"},
"2", {\ar"1"},
"1", {\ar"4"},
"5", {\cut"2"},
"4", {\ar"3"},
"3", {\cut"6"},
"4", {\ar"5"},
"7", {\ar"4"},
"5", {\cut"8"},
"6", {\ar"7"},
"8", {\ar"7"},
\end{xy}
} ="9",
(0,200) *+{
\begin{xy} 0;<0.17pt,0pt>:<0pt,-0.17pt>:: 
(0,0) *+{\bullet} ="0",
(125,0) *+{\bullet} ="1",
(250,0) *+{\bullet} ="2",
(0,125) *+{\bullet} ="3",
(125,125) *+{\bullet} ="4",
(250,125) *+{\bullet} ="5",
(0,250) *+{\bullet} ="6",
(125,250) *+{\bullet} ="7",
(250,250) *+{\bullet} ="8",
"0", {\cut"1"},
"3", {\ar"0"},
"2", {\ar"1"},
"1", {\ar"4"},
"5", {\ar"2"},
"4", {\ar"3"},
"3", {\ar"6"},
"4", {\cut"5"},
"7", {\ar"4"},
"5", {\ar"8"},
"6", {\cut"7"},
"8", {\ar"7"},
\end{xy}
} ="10",
(50,200) *+{
\begin{xy} 0;<0.17pt,0pt>:<0pt,-0.17pt>:: 
(0,0) *+{\bullet} ="0",
(125,0) *+{\bullet} ="1",
(250,0) *+{\bullet} ="2",
(0,125) *+{\bullet} ="3",
(125,125) *+{\bullet} ="4",
(250,125) *+{\bullet} ="5",
(0,250) *+{\bullet} ="6",
(125,250) *+{\bullet} ="7",
(250,250) *+{\bullet} ="8",
"0", {\ar"1"},
"3", {\cut"0"},
"2", {\ar"1"},
"1", {\ar"4"},
"5", {\cut"2"},
"4", {\ar"3"},
"3", {\ar"6"},
"4", {\ar"5"},
"7", {\ar"4"},
"5", {\cut"8"},
"6", {\cut"7"},
"8", {\ar"7"},
\end{xy}
} ="11",
(0,150) *+{
\begin{xy} 0;<0.17pt,0pt>:<0pt,-0.17pt>:: 
(0,0) *+{\bullet} ="0",
(125,0) *+{\bullet} ="1",
(250,0) *+{\bullet} ="2",
(0,125) *+{\bullet} ="3",
(125,125) *+{\bullet} ="4",
(250,125) *+{\bullet} ="5",
(0,250) *+{\bullet} ="6",
(125,250) *+{\bullet} ="7",
(250,250) *+{\bullet} ="8",
"0", {\cut"1"},
"3", {\ar"0"},
"2", {\ar"1"},
"1", {\ar"4"},
"5", {\cut"2"},
"4", {\ar"3"},
"3", {\ar"6"},
"4", {\ar"5"},
"7", {\ar"4"},
"5", {\cut"8"},
"6", {\cut"7"},
"8", {\ar"7"},
\end{xy}
} ="12",
(100,150) *+{
\begin{xy} 0;<0.17pt,0pt>:<0pt,-0.17pt>:: 
(0,0) *+{\bullet} ="0",
(125,0) *+{\bullet} ="1",
(250,0) *+{\bullet} ="2",
(0,125) *+{\bullet} ="3",
(125,125) *+{\bullet} ="4",
(250,125) *+{\bullet} ="5",
(0,250) *+{\bullet} ="6",
(125,250) *+{\bullet} ="7",
(250,250) *+{\bullet} ="8",
"0", {\ar"1"},
"3", {\cut"0"},
"2", {\ar"1"},
"1", {\ar"4"},
"5", {\cut"2"},
"4", {\ar"3"},
"3", {\ar"6"},
"4", {\ar"5"},
"7", {\ar"4"},
"5", {\ar"8"},
"6", {\cut"7"},
"8", {\cut"7"},
\end{xy}
} ="13",
(100,100) *+{
\begin{xy} 0;<0.17pt,0pt>:<0pt,-0.17pt>:: 
(0,0) *+{\bullet} ="0",
(125,0) *+{\bullet} ="1",
(250,0) *+{\bullet} ="2",
(0,125) *+{\bullet} ="3",
(125,125) *+{\bullet} ="4",
(250,125) *+{\bullet} ="5",
(0,250) *+{\bullet} ="6",
(125,250) *+{\bullet} ="7",
(250,250) *+{\bullet} ="8",
"0", {\ar"1"},
"3", {\cut"0"},
"2", {\ar"1"},
"1", {\ar"4"},
"5", {\cut"2"},
"4", {\ar"3"},
"3", {\ar"6"},
"4", {\ar"5"},
"7", {\cut"4"},
"5", {\ar"8"},
"6", {\ar"7"},
"8", {\ar"7"},
\end{xy}
} ="14",
"0", {\lat"1"},
"11", {\lat"0"},
"1", {\lat"2"},
"1", {\lat"3"},
"12", {\lat"1"},
"13", {\lat"1"},
"2", {\lat"4"},
"14", {\lat"2"},
"3", {\lat"4"},
"4", {\lat"5"},
"5", {\lat"6"},
"6", {\lat"7"},
"7", {\lat"8"},
"7", {\lat"9"},
"8", {\lat"10"},
"8", {\lat"11"},
"9", {\lat"11"},
"10", {\lat"12"},
"11", {\lat"12"},
"11", {\lat"13"},
"13", {\lat"14"},
\end{xy}
\begin{xy} 0;<1.4pt,0pt>:<0pt,-1.4pt>:: 
(100,50) *+{
\begin{xy} 0;<0.17pt,0pt>:<0pt,-0.17pt>:: 
(0,0) *+{\bullet} ="0",
(125,0) *+{\bullet} ="1",
(250,0) *+{\bullet} ="2",
(0,125) *+{\bullet} ="3",
(125,125) *+{\bullet} ="4",
(250,125) *+{\bullet} ="5",
(0,250) *+{\bullet} ="6",
(125,250) *+{\bullet} ="7",
(250,250) *+{\bullet} ="8",
"1", {\ar"0"},
"3", {\ar"0"},
"0", {\cut"4"},
"1", {\ar"2"},
"4", {\ar"1"},
"2", {\cut"4"},
"5", {\ar"2"},
"4", {\ar"3"},
"3", {\ar"6"},
"4", {\ar"5"},
"6", {\cut"4"},
"4", {\ar"7"},
"8", {\cut"4"},
"5", {\ar"8"},
"7", {\ar"6"},
"7", {\ar"8"},
\end{xy}
} ="0",
(50,0) *+{
\begin{xy} 0;<0.17pt,0pt>:<0pt,-0.17pt>:: 
(0,0) *+{\bullet} ="0",
(125,0) *+{\bullet} ="1",
(250,0) *+{\bullet} ="2",
(0,125) *+{\bullet} ="3",
(125,125) *+{\bullet} ="4",
(250,125) *+{\bullet} ="5",
(0,250) *+{\bullet} ="6",
(125,250) *+{\bullet} ="7",
(250,250) *+{\bullet} ="8",
"1", {\cut"0"},
"3", {\cut"0"},
"0", {\ar"4"},
"1", {\ar"2"},
"4", {\ar"1"},
"2", {\cut"4"},
"5", {\ar"2"},
"4", {\ar"3"},
"3", {\ar"6"},
"4", {\ar"5"},
"6", {\cut"4"},
"4", {\ar"7"},
"8", {\cut"4"},
"5", {\ar"8"},
"7", {\ar"6"},
"7", {\ar"8"},
\end{xy}
} ="1",
(50,50) *+{
\begin{xy} 0;<0.17pt,0pt>:<0pt,-0.17pt>:: 
(0,0) *+{\bullet} ="0",
(125,0) *+{\bullet} ="1",
(250,0) *+{\bullet} ="2",
(0,125) *+{\bullet} ="3",
(125,125) *+{\bullet} ="4",
(250,125) *+{\bullet} ="5",
(0,250) *+{\bullet} ="6",
(125,250) *+{\bullet} ="7",
(250,250) *+{\bullet} ="8",
"1", {\cut"0"},
"3", {\cut"0"},
"0", {\ar"4"},
"1", {\cut"2"},
"4", {\ar"1"},
"2", {\ar"4"},
"5", {\cut"2"},
"4", {\ar"3"},
"3", {\ar"6"},
"4", {\ar"5"},
"6", {\cut"4"},
"4", {\ar"7"},
"8", {\cut"4"},
"5", {\ar"8"},
"7", {\ar"6"},
"7", {\ar"8"},
\end{xy}
} ="2",
(0,50) *+{
\begin{xy} 0;<0.17pt,0pt>:<0pt,-0.17pt>:: 
(0,0) *+{\bullet} ="0",
(125,0) *+{\bullet} ="1",
(250,0) *+{\bullet} ="2",
(0,125) *+{\bullet} ="3",
(125,125) *+{\bullet} ="4",
(250,125) *+{\bullet} ="5",
(0,250) *+{\bullet} ="6",
(125,250) *+{\bullet} ="7",
(250,250) *+{\bullet} ="8",
"1", {\cut"0"},
"3", {\cut"0"},
"0", {\ar"4"},
"1", {\ar"2"},
"4", {\ar"1"},
"2", {\cut"4"},
"5", {\ar"2"},
"4", {\ar"3"},
"3", {\ar"6"},
"4", {\ar"5"},
"6", {\cut"4"},
"4", {\ar"7"},
"8", {\ar"4"},
"5", {\cut"8"},
"7", {\ar"6"},
"7", {\cut"8"},
\end{xy}
} ="3",
(0,100) *+{
\begin{xy} 0;<0.17pt,0pt>:<0pt,-0.17pt>:: 
(0,0) *+{\bullet} ="0",
(125,0) *+{\bullet} ="1",
(250,0) *+{\bullet} ="2",
(0,125) *+{\bullet} ="3",
(125,125) *+{\bullet} ="4",
(250,125) *+{\bullet} ="5",
(0,250) *+{\bullet} ="6",
(125,250) *+{\bullet} ="7",
(250,250) *+{\bullet} ="8",
"1", {\cut"0"},
"3", {\cut"0"},
"0", {\ar"4"},
"1", {\cut"2"},
"4", {\ar"1"},
"2", {\ar"4"},
"5", {\cut"2"},
"4", {\ar"3"},
"3", {\ar"6"},
"4", {\ar"5"},
"6", {\cut"4"},
"4", {\ar"7"},
"8", {\ar"4"},
"5", {\cut"8"},
"7", {\ar"6"},
"7", {\cut"8"},
\end{xy}
} ="4",
(50,100) *+{
\begin{xy} 0;<0.17pt,0pt>:<0pt,-0.17pt>:: 
(0,0) *+{\bullet} ="0",
(125,0) *+{\bullet} ="1",
(250,0) *+{\bullet} ="2",
(0,125) *+{\bullet} ="3",
(125,125) *+{\bullet} ="4",
(250,125) *+{\bullet} ="5",
(0,250) *+{\bullet} ="6",
(125,250) *+{\bullet} ="7",
(250,250) *+{\bullet} ="8",
"1", {\ar"0"},
"3", {\cut"0"},
"0", {\ar"4"},
"1", {\ar"2"},
"4", {\cut"1"},
"2", {\ar"4"},
"5", {\cut"2"},
"4", {\ar"3"},
"3", {\ar"6"},
"4", {\ar"5"},
"6", {\cut"4"},
"4", {\ar"7"},
"8", {\cut"4"},
"5", {\ar"8"},
"7", {\ar"6"},
"7", {\ar"8"},
\end{xy}
} ="5",
(0,150) *+{
\begin{xy} 0;<0.17pt,0pt>:<0pt,-0.17pt>:: 
(0,0) *+{\bullet} ="0",
(125,0) *+{\bullet} ="1",
(250,0) *+{\bullet} ="2",
(0,125) *+{\bullet} ="3",
(125,125) *+{\bullet} ="4",
(250,125) *+{\bullet} ="5",
(0,250) *+{\bullet} ="6",
(125,250) *+{\bullet} ="7",
(250,250) *+{\bullet} ="8",
"1", {\cut"0"},
"3", {\cut"0"},
"0", {\ar"4"},
"1", {\cut"2"},
"4", {\ar"1"},
"2", {\ar"4"},
"5", {\cut"2"},
"4", {\ar"3"},
"3", {\cut"6"},
"4", {\ar"5"},
"6", {\ar"4"},
"4", {\ar"7"},
"8", {\ar"4"},
"5", {\cut"8"},
"7", {\cut"6"},
"7", {\cut"8"},
\end{xy}
} ="6",
(50,150) *+{
\begin{xy} 0;<0.17pt,0pt>:<0pt,-0.17pt>:: 
(0,0) *+{\bullet} ="0",
(125,0) *+{\bullet} ="1",
(250,0) *+{\bullet} ="2",
(0,125) *+{\bullet} ="3",
(125,125) *+{\bullet} ="4",
(250,125) *+{\bullet} ="5",
(0,250) *+{\bullet} ="6",
(125,250) *+{\bullet} ="7",
(250,250) *+{\bullet} ="8",
"1", {\ar"0"},
"3", {\cut"0"},
"0", {\ar"4"},
"1", {\ar"2"},
"4", {\cut"1"},
"2", {\ar"4"},
"5", {\cut"2"},
"4", {\ar"3"},
"3", {\ar"6"},
"4", {\ar"5"},
"6", {\cut"4"},
"4", {\ar"7"},
"8", {\ar"4"},
"5", {\cut"8"},
"7", {\ar"6"},
"7", {\cut"8"},
\end{xy}
} ="7",
(0,200) *+{
\begin{xy} 0;<0.17pt,0pt>:<0pt,-0.17pt>:: 
(0,0) *+{\bullet} ="0",
(125,0) *+{\bullet} ="1",
(250,0) *+{\bullet} ="2",
(0,125) *+{\bullet} ="3",
(125,125) *+{\bullet} ="4",
(250,125) *+{\bullet} ="5",
(0,250) *+{\bullet} ="6",
(125,250) *+{\bullet} ="7",
(250,250) *+{\bullet} ="8",
"1", {\ar"0"},
"3", {\cut"0"},
"0", {\ar"4"},
"1", {\ar"2"},
"4", {\cut"1"},
"2", {\ar"4"},
"5", {\cut"2"},
"4", {\ar"3"},
"3", {\cut"6"},
"4", {\ar"5"},
"6", {\ar"4"},
"4", {\ar"7"},
"8", {\ar"4"},
"5", {\cut"8"},
"7", {\cut"6"},
"7", {\cut"8"},
\end{xy}
} ="8",
(50,200) *+{
\begin{xy} 0;<0.17pt,0pt>:<0pt,-0.17pt>:: 
(0,0) *+{\bullet} ="0",
(125,0) *+{\bullet} ="1",
(250,0) *+{\bullet} ="2",
(0,125) *+{\bullet} ="3",
(125,125) *+{\bullet} ="4",
(250,125) *+{\bullet} ="5",
(0,250) *+{\bullet} ="6",
(125,250) *+{\bullet} ="7",
(250,250) *+{\bullet} ="8",
"1", {\ar"0"},
"3", {\cut"0"},
"0", {\ar"4"},
"1", {\ar"2"},
"4", {\cut"1"},
"2", {\ar"4"},
"5", {\ar"2"},
"4", {\ar"3"},
"3", {\ar"6"},
"4", {\cut"5"},
"6", {\cut"4"},
"4", {\ar"7"},
"8", {\ar"4"},
"5", {\ar"8"},
"7", {\ar"6"},
"7", {\cut"8"},
\end{xy}
} ="9",
(50,250) *+{
\begin{xy} 0;<0.17pt,0pt>:<0pt,-0.17pt>:: 
(0,0) *+{\bullet} ="0",
(125,0) *+{\bullet} ="1",
(250,0) *+{\bullet} ="2",
(0,125) *+{\bullet} ="3",
(125,125) *+{\bullet} ="4",
(250,125) *+{\bullet} ="5",
(0,250) *+{\bullet} ="6",
(125,250) *+{\bullet} ="7",
(250,250) *+{\bullet} ="8",
"1", {\ar"0"},
"3", {\cut"0"},
"0", {\ar"4"},
"1", {\ar"2"},
"4", {\cut"1"},
"2", {\ar"4"},
"5", {\ar"2"},
"4", {\ar"3"},
"3", {\cut"6"},
"4", {\cut"5"},
"6", {\ar"4"},
"4", {\ar"7"},
"8", {\ar"4"},
"5", {\ar"8"},
"7", {\cut"6"},
"7", {\cut"8"},
\end{xy}
} ="10",
(0,250) *+{
\begin{xy} 0;<0.17pt,0pt>:<0pt,-0.17pt>:: 
(0,0) *+{\bullet} ="0",
(125,0) *+{\bullet} ="1",
(250,0) *+{\bullet} ="2",
(0,125) *+{\bullet} ="3",
(125,125) *+{\bullet} ="4",
(250,125) *+{\bullet} ="5",
(0,250) *+{\bullet} ="6",
(125,250) *+{\bullet} ="7",
(250,250) *+{\bullet} ="8",
"1", {\ar"0"},
"3", {\cut"0"},
"0", {\ar"4"},
"1", {\ar"2"},
"4", {\cut"1"},
"2", {\ar"4"},
"5", {\cut"2"},
"4", {\ar"3"},
"3", {\cut"6"},
"4", {\ar"5"},
"6", {\ar"4"},
"4", {\cut"7"},
"8", {\ar"4"},
"5", {\cut"8"},
"7", {\ar"6"},
"7", {\ar"8"},
\end{xy}
} ="11",
(50,300) *+{
\begin{xy} 0;<0.17pt,0pt>:<0pt,-0.17pt>:: 
(0,0) *+{\bullet} ="0",
(125,0) *+{\bullet} ="1",
(250,0) *+{\bullet} ="2",
(0,125) *+{\bullet} ="3",
(125,125) *+{\bullet} ="4",
(250,125) *+{\bullet} ="5",
(0,250) *+{\bullet} ="6",
(125,250) *+{\bullet} ="7",
(250,250) *+{\bullet} ="8",
"1", {\ar"0"},
"3", {\cut"0"},
"0", {\ar"4"},
"1", {\ar"2"},
"4", {\cut"1"},
"2", {\ar"4"},
"5", {\ar"2"},
"4", {\ar"3"},
"3", {\cut"6"},
"4", {\cut"5"},
"6", {\ar"4"},
"4", {\cut"7"},
"8", {\ar"4"},
"5", {\ar"8"},
"7", {\ar"6"},
"7", {\ar"8"},
\end{xy}
} ="12",
(100,250) *+{
\begin{xy} 0;<0.17pt,0pt>:<0pt,-0.17pt>:: 
(0,0) *+{\bullet} ="0",
(125,0) *+{\bullet} ="1",
(250,0) *+{\bullet} ="2",
(0,125) *+{\bullet} ="3",
(125,125) *+{\bullet} ="4",
(250,125) *+{\bullet} ="5",
(0,250) *+{\bullet} ="6",
(125,250) *+{\bullet} ="7",
(250,250) *+{\bullet} ="8",
"1", {\ar"0"},
"3", {\ar"0"},
"0", {\ar"4"},
"1", {\ar"2"},
"4", {\cut"1"},
"2", {\ar"4"},
"5", {\ar"2"},
"4", {\cut"3"},
"3", {\ar"6"},
"4", {\cut"5"},
"6", {\ar"4"},
"4", {\cut"7"},
"8", {\ar"4"},
"5", {\ar"8"},
"7", {\ar"6"},
"7", {\ar"8"},
\end{xy}
} ="13",
"1", {\lat"0"},
"0", {\lat"13"},
"2", {\lat"1"},
"3", {\lat"1"},
"4", {\lat"2"},
"5", {\lat"2"},
"4", {\lat"3"},
"6", {\lat"4"},
"7", {\lat"4"},
"7", {\lat"5"},
"8", {\lat"6"},
"8", {\lat"7"},
"9", {\lat"7"},
"10", {\lat"8"},
"11", {\lat"8"},
"10", {\lat"9"},
"12", {\lat"10"},
"12", {\lat"11"},
"13", {\lat"12"},
\end{xy}
\]
\caption{Cut-mutation lattices.}
\end{figure}
\end{example}

\section{The canvas of a QP}\label{canvas}
In this section we provide a sufficient condition for a QP to be fully compatible.

In the theory of dimer models QPs arise from finite bipartite tilings on compact Riemann surfaces (see \cite{Br,D,IU} and references there).
Of main interest is the case when the surface is a torus. This is because the topological properties of the torus yield desirable properties of the corresponding QP.

The idea behind this section is to start with the QP instead of the surface. To each QP $(Q,W)$ we assign a CW complex $X_{(Q,W)}$ called the canvas of $(Q,W)$. One may think of $(Q,W)$ as painted on $X_{(Q,W)}$ hence the name canvas. In the dimer case we recover the torus in this way. Therefore it it reasonable to hope that the topological properties of $X_{(Q,W)}$ will reflect interesting properties of $(Q,W)$.

It turns out that for many examples of selfinjective QPs the canvas can be realised as a bounded simply connected region of the plane. Since a bipartite tiling on the torus can be thought of as an infinite periodic tiling on the plane, this is in stark contrast to dimer models.

\subsection{Fundamental groups of QPs and simply connected QPs}
In order to define CW complexes we recall some basic notation. Denote by $D^n$ the closed $n$-dimensional disc for each $n>0$.
Thus $\partial D^n = S^{n-1}$, the $n-1$-dimensional sphere. In particular $D^1 = [0,1]$ and $\partial D^1 = S^0 = \{0,1\}$. 
We use the inductive definition of CW complexes. Thus a CW complex $X = \bigcup_{i \in \N} X^i$ with the weak topology, 
where $X^0$ is a discrete space and $X^{n}$ is obtained from $X^{n-1}$ by attaching $n$-dimensional discs $D^n_a$ along continuous maps $\phi_{a} : S^{n-1} \rightarrow X^{n-1}$. 
Hence there are continuous maps $\epsilon^n_{a} : D^n_a \rightarrow X$, embedding the interior of $D^n_a$ in $X$, whose restriction to $\partial D^n_a$ equals $\phi_{a}$. 
The images of the interior under these embedding will be referred to as $n$-cells. 
The elements in $X^0$ will be referred to as $0$-cells.

\begin{definition}
\begin{itemize}
\item[(a)] Given a quiver $Q$, define the CW complex $X_Q$ as follows. The $0$-cells are the vertices of $Q$, i.e., $X_Q^0 := Q_0$. 
The $1$-cells are indexed by the arrows $a \in Q_1$ with attaching maps $\phi_a : S^0 \rightarrow Q_0$ defined by $\phi_a(0) = s(a)$ and $\phi_a(1) = e(a)$. 
These are all cells, i.e, $X_Q = X_Q^1$.
\item[(b)] Let $(Q,W)$ be a QP. Write $W = \sum_{c}\lambda_c c$ where $c$ runs through all cycles in $Q$ and $\lambda_c \in K$. 
Let $Q_2$ be the set of cycles $c$ such that $\lambda_c \not = 0$. 
Define the CW complex $X_{(Q,W)}$ which we call the \emph{canvas} of $(Q,W)$ by attaching to $X_Q$ one $2$-cell for each $c = a_0\cdots a_{s-1}\in Q_2$ via the attaching map $\phi_c:S^1 \rightarrow X_Q$ defined by 
\[
\phi_c\left(\cos\left(\frac{2\pi}{s}(i+ t)\right),\sin\left(\frac{2\pi}{s}(i+ t)\right)\right) = \epsilon^1_{a_i}(t).
\]
for each integer $0\leq i < s$ and real number $0 \leq t < 1$.
\end{itemize}
\end{definition}

\begin{remark}
The space $X_{(Q,W)}$ is not invariant under right equivalence of QP.
\end{remark}

We proceed to compute the fundamental groups of $X_Q$ and $X_{(Q,W)}$. Fix a base point $x_0 \in Q_0 \subset X_Q$. Let $G_{x_0}(Q)$ be the group of equivalence classes of cyclic walks in $Q$ starting and ending at $x_0$. For each arrow $a \in Q_1$ we have two paths $f_{a}, f_{a^{-1}} : [0,1] \rightarrow X_Q$ defined by $f_{a}(t)= \epsilon^1_{a}(t)$ and $f_{a^{-1}}(t)= \epsilon^1_{a}(1-t)$.

For each cyclic walk $p = a_1^{s_1}\cdots a_\ell^{s_\ell}$ starting at $x_0$ ($s_i \in \{\pm1\}$) we get a closed curve
\[
f_p: [0,1]\rightarrow X_Q
\]
by composing the paths $f_{a_1^{s_1}},\ldots ,f_{a_\ell^{s_\ell}}$. This induces a group morphism $f : G_{x_0}(Q)\rightarrow \Pi^1_{x_0}(X_Q)$.

\begin{proposition}\label{CW1}
The group morphism 
\[
f : G_{x_0}(Q)\rightarrow \Pi^1_{x_0}(X_Q)
\]
is an isomorphism.
\end{proposition}
\begin{proof}
Choose a maximal subquiver $T \subset Q$ such that the underlying graph is a tree. Let $A = Q_1 \setminus T_1$. Since $X_T$ is a contractible closed subcomplex of $X_Q$ the quotient map $X_Q \rightarrow X_Q/X_T$ is a homotopy equivalence. Moreover, $X_Q/X_T$ is a wedge sum of loops, one for each element in $A$. Let $F_A$ be the free group on $A$ and $g : \Pi^1_{x_0}(X_Q) \rightarrow F_A$, the induced isomorphism. 

Every element in $G_{x_0}(Q)$ is written uniquely on the form $p = t_1a_1^{s_1}\cdots t_\ell a_\ell ^{s_\ell }t_{\ell+1}$, where $t_i$ is a reduced walk in $T$, $a_i \in A$ and $s_i \in \{\pm1\}$ such that $a_i^{s_i} \not = a_{i+1}^{-s_{i+1}}$. Moreover, $g(f(p)) = a_1^{s_1}\cdots a_\ell ^{s_\ell }$. Since every reduced word in $F_A$ appears as $g(f(p))$ for some $p \in G_{x_0}(Q)$, the morphism $g\circ f$ is an isomorphism, and therefore, so is $f$.
\end{proof}

Now choose for each cycle $c \in Q_2$ a walk $p_c$ from $x_0$ to the starting point of $c$. Let $N$ be the normal subgroup of $G_{x_0}(Q)$ generated by $\{ p_c c p_c^{-1} \; | \; c \in Q_2\}$.

\begin{proposition}\label{CW2}
The group morphism $f$ induces an isomorphism
\[
G_{x_0}(Q)/N\rightarrow \Pi^1_{x_0}(X_{(Q,W)}).
\]
\end{proposition}
\begin{proof}
According to \cite[Proposition~1.26]{Ha}, the inclusion of $X_Q$ in $X_{(Q,W)}$ induces an epimorphism $\Pi^1_{x_0}(X_Q) \rightarrow \Pi^1_{x_0}(X_{(Q,W)})$ with kernel generated by $\{ f_{p_c c p_c^{-1}} \; | \; c \in Q_2\}$. Hence the kernel is $f(N)$.
\end{proof}

\begin{definition}
We call a QP $(Q,W)$ \emph{simply connected} if $X_{(Q,W)}$ is simply connected. 
\end{definition}

We now obtain a sufficient condition for full compatibility.

\begin{proposition}\label{admissible degree}
Every simply connected QP is fully compatible.
\end{proposition}
\begin{proof}
Let $(Q,W)$ be a simply connected QP and $p$ a cyclic walk in $Q$ at $x \in Q_0$. Since $X_{(Q,W)}$ is simply connected $N = G_{x}(Q)$ by Proposition \ref{CW2} and thus $p \in N$. Hence there are cycles $c_i \in Q_2$, walks $p_i$ and $s_i\in\{\pm1\}$ such that
\[
p\sim (p_1c_1^{s_1}p_1^{-1})(p_2c_2^{s_2}p_2^{-1})
\cdots(p_\ell c_\ell^{s_\ell}p_\ell^{-1}).
\]
For any cut $C$ in $(Q,W)$ we have
\[
g_C(p)=\sum_{i=1}^\ell g_C(c_i^{s_i})=\sum_{i=1}^\ell s_i,
\]
which does not depend on $C$.
\end{proof}


The following is an immediate consequence of Theorem \ref{transitive}.

\begin{theorem}\label{transitive2}
Let $(Q,W)$ be a selfinjective, simply connected QP with enough cuts.
Then all truncated Jacobian algebras of $(Q,W)$ are iterated 2-APR tilts of each other.
In particular they are derived equivalent.
\end{theorem}

\section{Planar QPs}\label{planar}
In this section we introduce a class of QPs called planar. Many examples of selfinjective planar QPs are provided at the end of this section.

\begin{definition} 
A QP $(Q,W)$ is called \emph{planar} if $X_{(Q,W)}$ is simply connected and equipped with an embedding $\epsilon : X_{(Q,W)}\rightarrow \R^2$ 
\end{definition}

Observe that since planar QPs are simply connected they are fully compatible by Proposition \ref{admissible degree}.
In particular we have the following result.

\begin{theorem}
Let $(Q,W)$ be a selfinjective planar QP with enough cuts.
Then all truncated Jacobian algebras of $(Q,W)$ are iterated 2-APR tilts of each other.
In particular they are derived equivalent.
\end{theorem}

\begin{proposition}\label{generic}
Let $(Q,W)$ be a planar QP. Write $W = \sum_{c} \lambda_c c$ and define $W' = \sum_{\lambda_c \not =0} c$. Then $X_{(Q,W)} = X_{(Q,W')}$ and $(Q,W) \simeq (Q,W')$.
\end{proposition}

\begin{proof}
If $W =0$ there is nothing to show. Assume that this is not the case. Let $c$ be a cycle in $W$ containing an arrow $a$ on the boundary of $X_{(Q,W)}$. Consider the subquiver $\widehat{Q} = Q \setminus\{a\}$ and the potential $\widehat{W}$ on $\widehat{Q}$ satisfying $\widehat{W} + \lambda_c c = W$. The QP $(\widehat{Q},\widehat{W})$ is planar but $\widehat{W}$ contains fewer cycles than $W$. By induction we may assume that all coefficients in $\widehat{W}$ are $0$ or $1$. Now consider the automorphism $KQ \rightarrow KQ$ defined by $a \mapsto \lambda_c^{-1}a$ and $b \mapsto b$ for all $b \not = a$. It sends $W$ to $\widehat{W} + c$.
\end{proof}

\begin{remark}\label{planar convention}
Motivated by Proposition \ref{generic} we assume in the sequel that all planar QPs have potentials where every coefficient is $0$ or $1$. In particular, the potential of a planar QP $(Q,W)$ is determined by the restricted embedding $X_{Q}\rightarrow \R^2$. Therefore we will omit the potential and simply write $(Q,\epsilon)$ for a planar QP $(Q,W)$ with restricted embedding $\epsilon$.
\end{remark}

Let us consider 2-reduction of planar QPs.

\begin{definition}\label{planar 2-reduction}
Let $(Q,\epsilon)$ be a planar QP with at 2-cycle $ab$ not lying entirely on the boundary and set $Q' := Q\setminus\{a,b\}$ 
Then we define a \emph{planar 2-reduction} $(Q',\epsilon')$ of $(Q,\epsilon)$ as follows:
\begin{itemize}
\item[(a)] Both $a$ and $b$ are in the interior.
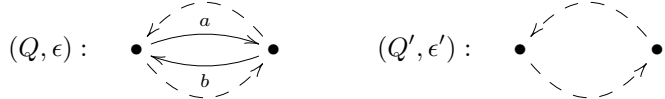
\[(Q,\epsilon):\ \ \ 
\xymatrix@R=.6cm@C=.6cm{
\bullet \ar@/^0.5pc/[rr]^a \ar@{-->}@/_1.5pc/[rr] & &\bullet  \ar@/^0.5pc/[ll]^b \ar@{-->}@/_1.5pc/[ll]
}\ \ \ \ \ \ \ \ \ \ 
(Q',\epsilon'):\ \ \ 
\xymatrix@R=.6cm@C=.6cm{
\bullet \ar@{-->}@/_1.5pc/[rr] & &\bullet \ar@{-->}@/_1.5pc/[ll]
}
\]
\item[(b)] The arrow $a$ is in the interior and $b$ is on the boundary (the case $b$ is in the interior and $a$ is on the boundary is similar).
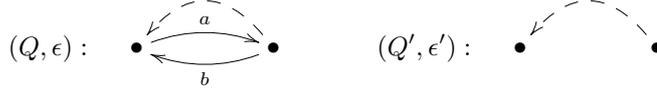
\[(Q,\epsilon):\ \ \ 
\xymatrix@R=.6cm@C=.6cm{
\bullet \ar@/^0.5pc/[rr]^a& &\bullet  \ar@/^0.5pc/[ll]^b \ar@{-->}@/_1.5pc/[ll]
}\ \ \ \ \ \ \ \ \ \ 
(Q',\epsilon'):\ \ \ 
\xymatrix@R=.6cm@C=.6cm{
\bullet & &\bullet \ar@{-->}@/_1.5pc/[ll]
}
\]
\end{itemize}
\end{definition}

\begin{proposition}
$(Q,\epsilon)$ is right equivalent to a direct sum of $(Q',\epsilon')$ and a trivial QP.
In particular we have $\P(Q,\epsilon)\simeq\P(Q',\epsilon')$.
\end{proposition}

\begin{proof}
We only show the case (a). We name the paths in $Q$ as follows:
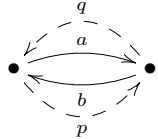
\[\xymatrix@R=.6cm@C=.6cm{
\bullet \ar@/^0.5pc/[rr]^a \ar@{-->}@/_1.5pc/[rr]_p & &\bullet  \ar@/^0.5pc/[ll]^b \ar@{-->}@/_1.5pc/[ll]_q
}\]
Then we can write $W=W'+ap-ab+qb$ where any cycle in $W'$ does not contain $a$ and $b$.
Since $W=W'+qp+(a-q)(p-b)$, the automorphism $f:\widehat{KQ}\to\widehat{KQ}$ defined by
$f(a)=a-b$, $f(b)=b-p$ and $f(c)=c$ for any $c\in Q_1\backslash\{a,b\}$
satisfies $f(W'+qp+ab)=W$.
Since the QP $(Q,W'+qp+ab)$ is a direct sum of $(Q',\epsilon')$ and a trivial QP, we have the assertion.
\end{proof}

We proceed to interpret mutation in the planar setting.

\subsection{Planar mutation}
Let $(Q,\epsilon)$ be a planar QP.
We fix a vertex $k\in Q_0$ satisfying one of the following conditions:
\begin{itemize}
\item[(i)] $k$ is in the interior and exactly 4 arrows start or end at $k$.
\item[(ii)] $k$ is at the boundary and at most 4 arrows start or end at $k$.
\end{itemize}

\begin{definition}\label{planar mutation}
We regard $\mu_k(Q,\epsilon)$ as a planar QP as follows (dashed arrows denote to paths):
\begin{itemize}
\item[(a)] $k$ is in the interior and exactly 4 arrows start or end at $k$.
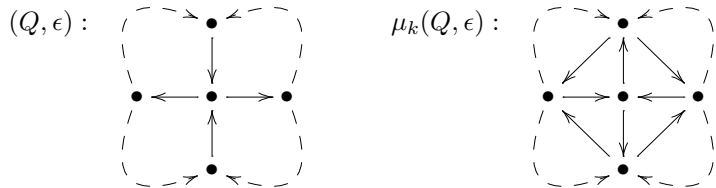
\[(Q,\epsilon):\ \ \ \xymatrix@R=.6cm@C=.6cm{
&\bullet\ar[d]&\\
\bullet\ar@{-->}@/^2pc/[ur]\ar@{-->}@/_2pc/[dr]&\bullet\ar[r]\ar[l]&\bullet\ar@{-->}@/_2pc/[ul]\ar@{-->}@/^2pc/[dl]\\
&\bullet\ar[u]&
}\ \ \ \ \ \ \ \ \ \ 
\mu_k(Q,\epsilon):\ \ \ \xymatrix@=.6cm@C=.6cm{
&\bullet\ar[ld]\ar[rd]&\\
\bullet\ar@{-->}@/^2pc/[ur]\ar@{-->}@/_2pc/[dr]\ar[r]&\bullet\ar[u]\ar[d]&\bullet\ar@{-->}@/_2pc/[ul]\ar[l]\ar@{-->}@/^2pc/[dl]\\
&\bullet\ar[lu]\ar[ru]&
}
\]
\item[(b)] $k$ is at the boundary and exactly 4 arrows start or end at $k$.
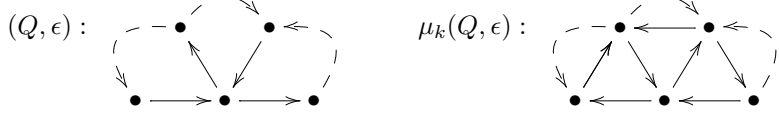
\[(Q,\epsilon):\ \ \ \xymatrix@R=.6cm@C=.2cm{
&\bullet\ar@{-->}@/_1.5pc/[dl]\ar@{-->}@/^1pc/[rr]&&\bullet\ar[dl]\\
\bullet\ar[rr]&&\bullet\ar[rr]\ar[ul]&&\bullet\ar@{-->}@/_1.5pc/[ul]
}\ \ \ \ \ \ \ \ \ \ 
\mu_k(Q,\epsilon):\ \ \ \xymatrix@R=.6cm@C=.2cm{
&\bullet\ar[dr]\ar@{-->}@/_1.5pc/[dl]\ar@{-->}@/^1pc/[rr]&&\bullet\ar[dr]\ar[ll]\\
\bullet\ar[ru]\ar[ur]&&\bullet\ar[ll]\ar[ur]&&\bullet\ar[ll]\ar@{-->}@/_1.5pc/[ul]
}\]
\item[(c)] $k$ is at the boundary and exactly 3 arrows start or end at $k$.
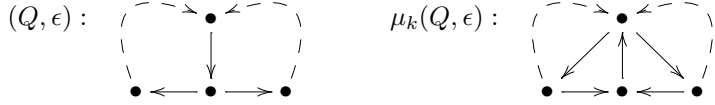
\[(Q,\epsilon):\ \ \ \xymatrix@R=.6cm@C=.6cm{
&\bullet\ar[d]&\\
\bullet\ar@{-->}@/^2pc/[ur]&\bullet\ar[r]\ar[l]&\bullet\ar@{-->}@/_2pc/[ul]
}\ \ \ \ \ \ \ \ \ \ 
\mu_k(Q,\epsilon):\ \ \ \xymatrix@=.6cm@C=.6cm{
&\bullet\ar[ld]\ar[rd]&\\
\bullet\ar@{-->}@/^2pc/[ur]\ar[r]&\bullet\ar[u]&\bullet\ar@{-->}@/_2pc/[ul]\ar[l]
}
\]
\item[(d)] $k$ is at the boundary and exactly 2 arrows start or end at $k$.
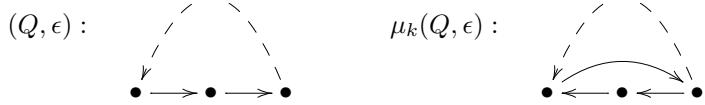
\[(Q,\epsilon):\ \ \ \xymatrix@R=.6cm@C=.6cm{
\\
\bullet\ar[r]&\bullet\ar[r]&\bullet\ar@{-->}@/_3pc/[ll]
}\ \ \ \ \ \ \ \ \ \ 
\mu_k(Q,\epsilon):\ \ \ \xymatrix@=.6cm@C=.6cm{
\\
\bullet\ar@/^1pc/[rr]&\bullet\ar[l]&\bullet\ar[l]\ar@{-->}@/_3pc/[ll]
}
\]
\end{itemize}
One can easily check that the above planar mutation of QPs is compatible with usual Derksen-Weyman-Zelevinsky mutation of QPs.
\end{definition}



\begin{theorem}\label{selfinj planar mutation}
Assume that $(Q,\epsilon)$ is a selfinjective planar QP and $k \in Q_0$ satisfies the conditions in Definition \ref{succ mutation}. 
Then $\mu_{(k)}(Q,\epsilon)$ is again a selfinjective planar QP.
\end{theorem}

\begin{proof}
Applying Definition \ref{planar mutation} repeatedly, we have that $\mu_{(k)}(Q,\epsilon)$ is planar. Also $\mu_{(k)}(Q,\epsilon)$ is selfinjective by Theorem \ref{selfinjective QP}.
\end{proof}

We say that two selfinjective planar QPs are related by planar mutation if one is obtained from the other by iterating the process described in Theorem \ref{selfinj planar mutation}.

In the rest of this section we exhibit explicit selfinjective planar QPs by considering already known examples and applying planar mutation.

\subsection{Examples from triangles}\label{planar 3}
The first family of selfinjective planar QPs that we consider are the $3$-preprojective algebras of Auslander algebras of quivers of type $A$ with linear orientation. 
We recall their description from \cite{IO1}. Define $Q = Q^{(s)}$ as follows:
\[
Q_0 := \{ (x_1,x_2,x_3) \in \Z^3_{\geq 0} \; | \; x_1+x_2+x_3 = s-1\}
\]
\[
Q_1 := \{x\arr{i} x+f_i \; | \; 1\leq i\leq 3, \;  x,x+f_i \in Q_0\}
\]
where $f_1 := (-1,1,0)$, $f_2 := (0,-1,1)$ and $f_3 := (1,0,-1)$. The vertices $Q_0$ lie in the plane in $\R^3$ defined by the equation $x_1+x_2+x_3 = s-1$. Identifying this with $\R^2$ we get an induced embedding $\epsilon^{(s)} : X_{Q^{(s)}} \rightarrow \R^2$. For instance, $(Q^{(5)}, \epsilon^{(5)})$ is the following planar QP.
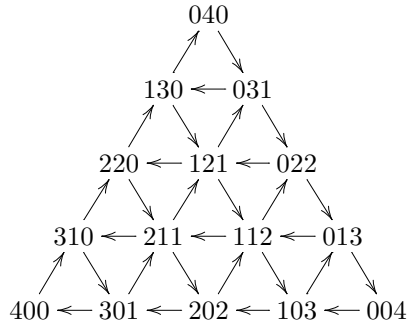
\[
\begin{xy} 
<0pt,0pt>;<0.5pt,0pt>:<0pt,-0.5pt>::
(134,0) *+{040} ="0",
(100,56) *+{130} ="1",
(167,56) *+{031} ="2",
(67,112) *+{220} ="3",
(134,112) *+{121} ="4",
(200,112) *+{022} ="5",
(33,167) *+{310} ="6",
(100,167) *+{211} ="7",
(167,167) *+{112} ="8",
(234,167) *+{013} ="9",
(0,223) *+{400} ="10",
(67,223) *+{301} ="11",
(134,223) *+{202} ="12",
(200,223) *+{103} ="13",
(267,223) *+{004} ="14",
"1", {\ar"0"},
"0", {\ar"2"},
"2", {\ar"1"},
"3", {\ar"1"},
"1", {\ar"4"},
"4", {\ar"2"},
"2", {\ar"5"},
"4", {\ar"3"},
"6", {\ar"3"},
"3", {\ar"7"},
"5", {\ar"4"},
"7", {\ar"4"},
"4", {\ar"8"},
"8", {\ar"5"},
"5", {\ar"9"},
"7", {\ar"6"},
"10", {\ar"6"},
"6", {\ar"11"},
"8", {\ar"7"},
"11", {\ar"7"},
"7", {\ar"12"},
"9", {\ar"8"},
"12", {\ar"8"},
"8", {\ar"13"},
"13", {\ar"9"},
"9", {\ar"14"},
"11", {\ar"10"},
"12", {\ar"11"},
"13", {\ar"12"},
"14", {\ar"13"},
\end{xy}
\]

For each $1\leq i \leq 3$ the subset $\{x\arr{i} x+f_i \; | \; x,x+f_i \in Q_0\} \subset Q_1$ is a cut. Thus all QPs $(Q^{(s)}, \epsilon^{(s)})$ have enough cuts.

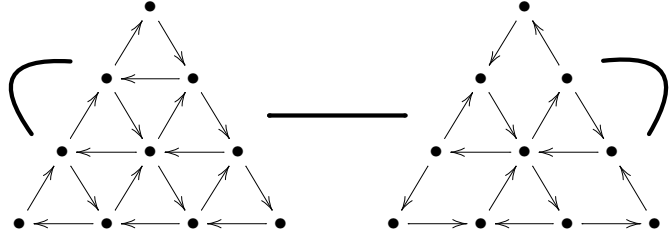
\begin{figure}[h]
\[
\begin{xy} 0;<1pt,0pt>:<0pt,-1pt>::
(0,0) *+{
\begin{xy} 0;<0.38pt,0pt>:<0pt,-0.38pt>:: 
(129,0) *+{\point 0} ="0",
(86,71) *+{\point 1} ="1",
(171,71) *+{\point 2} ="2",
(42,144) *+{\point 3} ="3",
(129,144) *+{\point 4} ="4",
(215,144) *+{\point 5} ="5",
(0,215) *+{\point 6} ="6",
(86,215) *+{\point 7} ="7",
(171,215) *+{\point 8} ="8",
(257,215) *+{\point 9} ="9",
"1", {\ar"0"},
"0", {\ar"2"},
"2", {\ar"1"},
"3", {\ar"1"},
"1", {\ar"4"},
"4", {\ar"2"},
"2", {\ar"5"},
"4", {\ar"3"},
"6", {\ar"3"},
"3", {\ar"7"},
"5", {\ar"4"},
"7", {\ar"4"},
"4", {\ar"8"},
"8", {\ar"5"},
"5", {\ar"9"},
"7", {\ar"6"},
"8", {\ar"7"},
"9", {\ar"8"},
\end{xy}
} ="1",
(140,0) *+{
\begin{xy} 0;<0.38pt,0pt>:<0pt,-0.38pt>:: 
(129,0) *+{\point 0} ="0",
(86,71) *+{\point 1} ="1",
(171,71) *+{\point 2} ="2",
(42,144) *+{\point 3} ="3",
(129,144) *+{\point 4} ="4",
(215,144) *+{\point 5} ="5",
(0,215) *+{\point 6} ="6",
(86,215) *+{\point 7} ="7",
(171,215) *+{\point 8} ="8",
(257,215) *+{\point 9} ="9",
"0", {\ar"1"},
"2", {\ar"0"},
"3", {\ar"1"},
"1", {\ar"4"},
"4", {\ar"2"},
"2", {\ar"5"},
"4", {\ar"3"},
"3", {\ar"6"},
"5", {\ar"4"},
"7", {\ar"4"},
"4", {\ar"8"},
"9", {\ar"5"},
"6", {\ar"7"},
"8", {\ar"7"},
"8", {\ar"9"},
\end{xy}
} ="2",
(40,0) *+{\mbox{ }}="12",
(100,0) *+{\mbox{ }}="21",
(-25,-20) *+{\mbox{ }}="111",
(-42,10) *+{\mbox{ }}="112",
"111",{\lat@/_20pt/ "112"},
(165,-20) *+{\mbox{ }}="221",
(185,10) *+{\mbox{ }}="222",
"221",{\lat@/^20pt/ "222"},
"12", {\lat"21"},
\end{xy}
\]
\caption{Planar mutation lattice for $(Q^{(4)}, W^{(4)})$.}
\label{mutation lattice A4}
\end{figure}

\begin{figure}[h]
\[
\begin{xy} 0;<1pt,0pt>:<0pt,-1pt>::
(0,0) *+{
\begin{xy} 0;<0.38pt,0pt>:<0pt,-0.38pt>:: 
(134,0) *+{\bullet} ="0",
(100,56) *+{\bullet} ="1",
(167,56) *+{\bullet} ="2",
(86,121) *+{\bullet} ="3",
(134,94) *+{\bullet} ="4",
(181,121) *+{\bullet} ="5",
(33,167) *+{\bullet} ="6",
(86,173) *+{\bullet} ="7",
(181,173) *+{\bullet} ="8",
(234,167) *+{\bullet} ="9",
(0,223) *+{\bullet} ="10",
(67,223) *+{\bullet} ="11",
(134,204) *+{\bullet} ="12",
(200,223) *+{\bullet} ="13",
(267,223) *+{\bullet} ="14",
"1", {\ar"0"},
"0", {\ar"2"},
"2", {\ar"1"},
"1", {\ar"3"},
"6", {\ar"1"},
"5", {\ar"2"},
"2", {\ar"9"},
"3", {\ar"4"},
"3", {\ar"6"},
"7", {\ar"3"},
"4", {\ar"5"},
"5", {\ar"8"},
"9", {\ar"5"},
"10", {\ar"6"},
"6", {\ar"11"},
"12", {\ar"7"},
"8", {\ar"12"},
"13", {\ar"9"},
"9", {\ar"14"},
"11", {\ar"10"},
"11", {\ar"12"},
"13", {\ar"11"},
"12", {\ar"13"},
"14", {\ar"13"},
\end{xy}
} ="1",
(140,0) *+{
\begin{xy} 
0;<0.38pt,0pt>:<0pt,-0.38pt>:: 
(134,0) *+{\bullet} ="0",
(100,56) *+{\bullet} ="1",
(167,56) *+{\bullet} ="2",
(86,121) *+{\bullet} ="3",
(134,94) *+{\bullet} ="4",
(181,121) *+{\bullet} ="5",
(33,167) *+{\bullet} ="6",
(86,173) *+{\bullet} ="7",
(181,173) *+{\bullet} ="8",
(234,167) *+{\bullet} ="9",
(0,223) *+{\bullet} ="10",
(67,223) *+{\bullet} ="11",
(134,204) *+{\bullet} ="12",
(200,223) *+{\bullet} ="13",
(267,223) *+{\bullet} ="14",
"0", {\ar"1"},
"2", {\ar"0"},
"1", {\ar"3"},
"6", {\ar"1"},
"5", {\ar"2"},
"2", {\ar"9"},
"3", {\ar"4"},
"3", {\ar"6"},
"7", {\ar"3"},
"4", {\ar"5"},
"5", {\ar"8"},
"9", {\ar"5"},
"6", {\ar"10"},
"12", {\ar"7"},
"8", {\ar"12"},
"14", {\ar"9"},
"10", {\ar"11"},
"11", {\ar"12"},
"13", {\ar"11"},
"12", {\ar"13"},
"13", {\ar"14"},
\end{xy} 
} ="2",
(0,115) *+{
\begin{xy} 
<0pt,0pt>;<0.4pt,0pt>:<0pt,-0.4pt>:: 
(134,0) *+{\bullet} ="0",
(100,56) *+{\bullet} ="1",
(167,56) *+{\bullet} ="2",
(67,112) *+{\bullet} ="3",
(134,112) *+{\bullet} ="4",
(200,112) *+{\bullet} ="5",
(33,167) *+{\bullet} ="6",
(100,167) *+{\bullet} ="7",
(167,167) *+{\bullet} ="8",
(234,167) *+{\bullet} ="9",
(0,223) *+{\bullet} ="10",
(67,223) *+{\bullet} ="11",
(134,223) *+{\bullet} ="12",
(200,223) *+{\bullet} ="13",
(267,223) *+{\bullet} ="14",
"1", {\ar"0"},
"0", {\ar"2"},
"2", {\ar"1"},
"3", {\ar"1"},
"1", {\ar"4"},
"4", {\ar"2"},
"2", {\ar"5"},
"4", {\ar"3"},
"6", {\ar"3"},
"3", {\ar"7"},
"5", {\ar"4"},
"7", {\ar"4"},
"4", {\ar"8"},
"8", {\ar"5"},
"5", {\ar"9"},
"7", {\ar"6"},
"10", {\ar"6"},
"6", {\ar"11"},
"8", {\ar"7"},
"11", {\ar"7"},
"7", {\ar"12"},
"9", {\ar"8"},
"12", {\ar"8"},
"8", {\ar"13"},
"13", {\ar"9"},
"9", {\ar"14"},
"11", {\ar"10"},
"12", {\ar"11"},
"13", {\ar"12"},
"14", {\ar"13"},
\end{xy}
} ="3",
(140,115) *+{
\begin{xy} 
0;<0.4pt,0pt>:<0pt,-0.4pt>:: 
(134,0) *+{\bullet} ="0",
(100,56) *+{\bullet} ="1",
(167,56) *+{\bullet} ="2",
(67,112) *+{\bullet} ="3",
(134,112) *+{\bullet} ="4",
(200,112) *+{\bullet} ="5",
(33,167) *+{\bullet} ="6",
(100,167) *+{\bullet} ="7",
(167,167) *+{\bullet} ="8",
(234,167) *+{\bullet} ="9",
(0,223) *+{\bullet} ="10",
(67,223) *+{\bullet} ="11",
(134,223) *+{\bullet} ="12",
(200,223) *+{\bullet} ="13",
(267,223) *+{\bullet} ="14",
"0", {\ar"1"},
"2", {\ar"0"},
"3", {\ar"1"},
"1", {\ar"4"},
"4", {\ar"2"},
"2", {\ar"5"},
"4", {\ar"3"},
"6", {\ar"3"},
"3", {\ar"7"},
"5", {\ar"4"},
"7", {\ar"4"},
"4", {\ar"8"},
"8", {\ar"5"},
"5", {\ar"9"},
"7", {\ar"6"},
"6", {\ar"10"},
"8", {\ar"7"},
"11", {\ar"7"},
"7", {\ar"12"},
"9", {\ar"8"},
"12", {\ar"8"},
"8", {\ar"13"},
"14", {\ar"9"},
"10", {\ar"11"},
"12", {\ar"11"},
"13", {\ar"12"},
"13", {\ar"14"},
\end{xy}
} ="4",
(260,180) *+{
\begin{xy} 
<0pt,0pt>;<0.37pt,0pt>:<0pt,-0.37pt>:: 
(108,0) *+{\bullet} ="0",
(72,60) *+{\bullet} ="1",
(180,0) *+{\bullet} ="2",
(34,126) *+{\bullet} ="3",
(141,99) *+{\bullet} ="4",
(217,60) *+{\bullet} ="5",
(0,188) *+{\bullet} ="6",
(100,167) *+{\bullet} ="7",
(178,167) *+{\bullet} ="8",
(252,126) *+{\bullet} ="9",
(33,243) *+{\bullet} ="10",
(98,243) *+{\bullet} ="11",
(178,243) *+{\bullet} ="12",
(249,243) *+{\bullet} ="13",
(283,188) *+{\bullet} ="14",
"1", {\ar"0"},
"2", {\ar"0"},
"0", {\ar"4"},
"1", {\ar"3"},
"4", {\ar"1"},
"4", {\ar"2"},
"2", {\ar"5"},
"6", {\ar"3"},
"3", {\ar"7"},
"5", {\ar"4"},
"7", {\ar"4"},
"4", {\ar"8"},
"9", {\ar"5"},
"7", {\ar"6"},
"6", {\ar"10"},
"8", {\ar"7"},
"10", {\ar"7"},
"7", {\ar"11"},
"8", {\ar"9"},
"12", {\ar"8"},
"8", {\ar"13"},
"14", {\ar"8"},
"9", {\ar"14"},
"11", {\ar"10"},
"11", {\ar"12"},
"13", {\ar"12"},
"13", {\ar"14"},
\end{xy}
} ="5",
(0,240) *+{
\begin{xy} 
<0pt,0pt>;<0.4pt,0pt>:<0pt,-0.4pt>::
(134,0) *+{\bullet} ="2",
(100,56) *+{\bullet} ="1",
(167,56) *+{\bullet} ="5",
(25,89) *+{\bullet} ="0",
(134,112) *+{\bullet} ="4",
(242,89) *+{\bullet} ="14",
(33,167) *+{\bullet} ="3",
(100,167) *+{\bullet} ="7",
(167,167) *+{\bullet} ="8",
(234,167) *+{\bullet} ="9",
(0,223) *+{\bullet} ="6",
(67,223) *+{\bullet} ="11",
(134,267) *+{\bullet} ="10",
(200,223) *+{\bullet} ="12",
(267,223) *+{\bullet} ="13",
"0", {\ar"1"},
"3", {\ar"0"},
"1", {\ar"2"},
"1", {\ar"3"},
"4", {\ar"1"},
"2", {\ar"5"},
"6", {\ar"3"},
"3", {\ar"7"},
"5", {\ar"4"},
"7", {\ar"4"},
"4", {\ar"8"},
"9", {\ar"5"},
"5", {\ar"14"},
"11", {\ar"6"},
"8", {\ar"7"},
"7", {\ar"11"},
"8", {\ar"9"},
"12", {\ar"8"},
"9", {\ar"13"},
"14", {\ar"9"},
"10", {\ar"11"},
"12", {\ar"10"},
"11", {\ar"12"},
"13", {\ar"12"},
\end{xy}
} ="6",
(140,232) *+{
\begin{xy} 
0;<0.4pt,0pt>:<0pt,-0.4pt>::  
(134,0) *+{\bullet} ="2",
(100,56) *+{\bullet} ="1",
(167,56) *+{\bullet} ="5",
(67,112) *+{\bullet} ="0",
(134,112) *+{\bullet} ="4",
(200,112) *+{\bullet} ="14",
(33,167) *+{\bullet} ="3",
(100,167) *+{\bullet} ="7",
(167,167) *+{\bullet} ="8",
(234,167) *+{\bullet} ="9",
(0,223) *+{\bullet} ="6",
(67,223) *+{\bullet} ="11",
(134,223) *+{\bullet} ="10",
(200,223) *+{\bullet} ="12",
(267,223) *+{\bullet} ="13",
"1", {\ar"0"},
"0", {\ar"3"},
"1", {\ar"2"},
"4", {\ar"1"},
"2", {\ar"5"},
"6", {\ar"3"},
"3", {\ar"7"},
"5", {\ar"4"},
"7", {\ar"4"},
"4", {\ar"8"},
"14", {\ar"5"},
"11", {\ar"6"},
"8", {\ar"7"},
"7", {\ar"11"},
"8", {\ar"9"},
"12", {\ar"8"},
"9", {\ar"13"},
"9", {\ar"14"},
"11", {\ar"10"},
"10", {\ar"12"},
"13", {\ar"12"},
\end{xy}
} ="7",
(0,370) *+{
\begin{xy} 
0;<0.4pt,0pt>:<0pt,-0.4pt>:: 
(134,0) *+{\bullet} ="2",
(100,56) *+{\bullet} ="1",
(167,56) *+{\bullet} ="5",
(25,89) *+{\bullet} ="0",
(134,112) *+{\bullet} ="4",
(242,89) *+{\bullet} ="14",
(33,167) *+{\bullet} ="3",
(100,167) *+{\bullet} ="7",
(167,167) *+{\bullet} ="8",
(234,167) *+{\bullet} ="9",
(0,223) *+{\bullet} ="6",
(67,223) *+{\bullet} ="11",
(134,267) *+{\bullet} ="10",
(200,223) *+{\bullet} ="12",
(267,223) *+{\bullet} ="13",
"0", {\ar"1"},
"3", {\ar"0"},
"2", {\ar"1"},
"1", {\ar"3"},
"4", {\ar"1"},
"1", {\ar"5"},
"5", {\ar"2"},
"3", {\ar"6"},
"3", {\ar"7"},
"11", {\ar"3"},
"5", {\ar"4"},
"7", {\ar"4"},
"4", {\ar"8"},
"9", {\ar"5"},
"5", {\ar"14"},
"6", {\ar"11"},
"8", {\ar"7"},
"7", {\ar"11"},
"8", {\ar"9"},
"12", {\ar"8"},
"9", {\ar"12"},
"13", {\ar"9"},
"14", {\ar"9"},
"10", {\ar"11"},
"12", {\ar"10"},
"11", {\ar"12"},
"12", {\ar"13"},
\end{xy}
} ="8",
(140,362) *+{
\begin{xy} 
0;<0.4pt,0pt>:<0pt,-0.4pt>:: 
(134,0) *+{\bullet} ="2",
(100,56) *+{\bullet} ="1",
(167,56) *+{\bullet} ="5",
(67,112) *+{\bullet} ="0",
(134,112) *+{\bullet} ="4",
(200,112) *+{\bullet} ="14",
(33,167) *+{\bullet} ="3",
(100,167) *+{\bullet} ="7",
(167,167) *+{\bullet} ="8",
(234,167) *+{\bullet} ="9",
(0,223) *+{\bullet} ="6",
(67,223) *+{\bullet} ="11",
(134,223) *+{\bullet} ="10",
(200,223) *+{\bullet} ="12",
(267,223) *+{\bullet} ="13",
"1", {\ar"0"},
"0", {\ar"3"},
"2", {\ar"1"},
"4", {\ar"1"},
"1", {\ar"5"},
"5", {\ar"2"},
"3", {\ar"6"},
"3", {\ar"7"},
"11", {\ar"3"},
"5", {\ar"4"},
"7", {\ar"4"},
"4", {\ar"8"},
"14", {\ar"5"},
"6", {\ar"11"},
"8", {\ar"7"},
"7", {\ar"11"},
"8", {\ar"9"},
"12", {\ar"8"},
"9", {\ar"12"},
"13", {\ar"9"},
"9", {\ar"14"},
"11", {\ar"10"},
"10", {\ar"12"},
"12", {\ar"13"},
\end{xy}
} ="9",
(40,0) *+{\mbox{ }}="12",
(100,0) *+{\mbox{ }}="21",
(40,115) *+{\mbox{ }}="34",
(100,115) *+{\mbox{ }}="43",
(55,232) *+{\mbox{ }}="67",
(100,232) *+{\mbox{ }}="76",
(55,362) *+{\mbox{ }}="89",
(100,362) *+{\mbox{ }}="98",
(180,115) *+{\mbox{ }}="45",
(228,150) *+{\mbox{ }}="54",
(177,233) *+{\mbox{ }}="75",
(207,217) *+{\mbox{ }}="57",
(-25,-20) *+{\mbox{ }}="111",
(-42,10) *+{\mbox{ }}="112",
"111",{\lat@/_20pt/ "112"},
(165,-20) *+{\mbox{ }}="221",
(185,10) *+{\mbox{ }}="222",
"221",{\lat@/^20pt/ "222"},
(-50,360) *+{\mbox{ }}="881",
(-50,390) *+{\mbox{ }}="882",
"881",{\lat@/_20pt/ "882"},
(170,355) *+{\mbox{ }}="991",
(190,388) *+{\mbox{ }}="992",
"991",{\lat@/^20pt/ "992"},
"12", {\lat"21"},
"1", {\lat"3"},
"2", {\lat"4"},
"34", {\lat"43"},
"3", {\lat"6"},
"4", {\lat"7"},
"45", {\lat"54"},
"57", {\lat"75"},
"67", {\lat"76"},
"6", {\lat"8"},
"7", {\lat"9"},
"89", {\lat"98"},
\end{xy}
\]
\caption{Planar mutation lattice for $(Q^{(5)}, W^{(5)})$.}
\label{mutation lattice A5}
\end{figure}
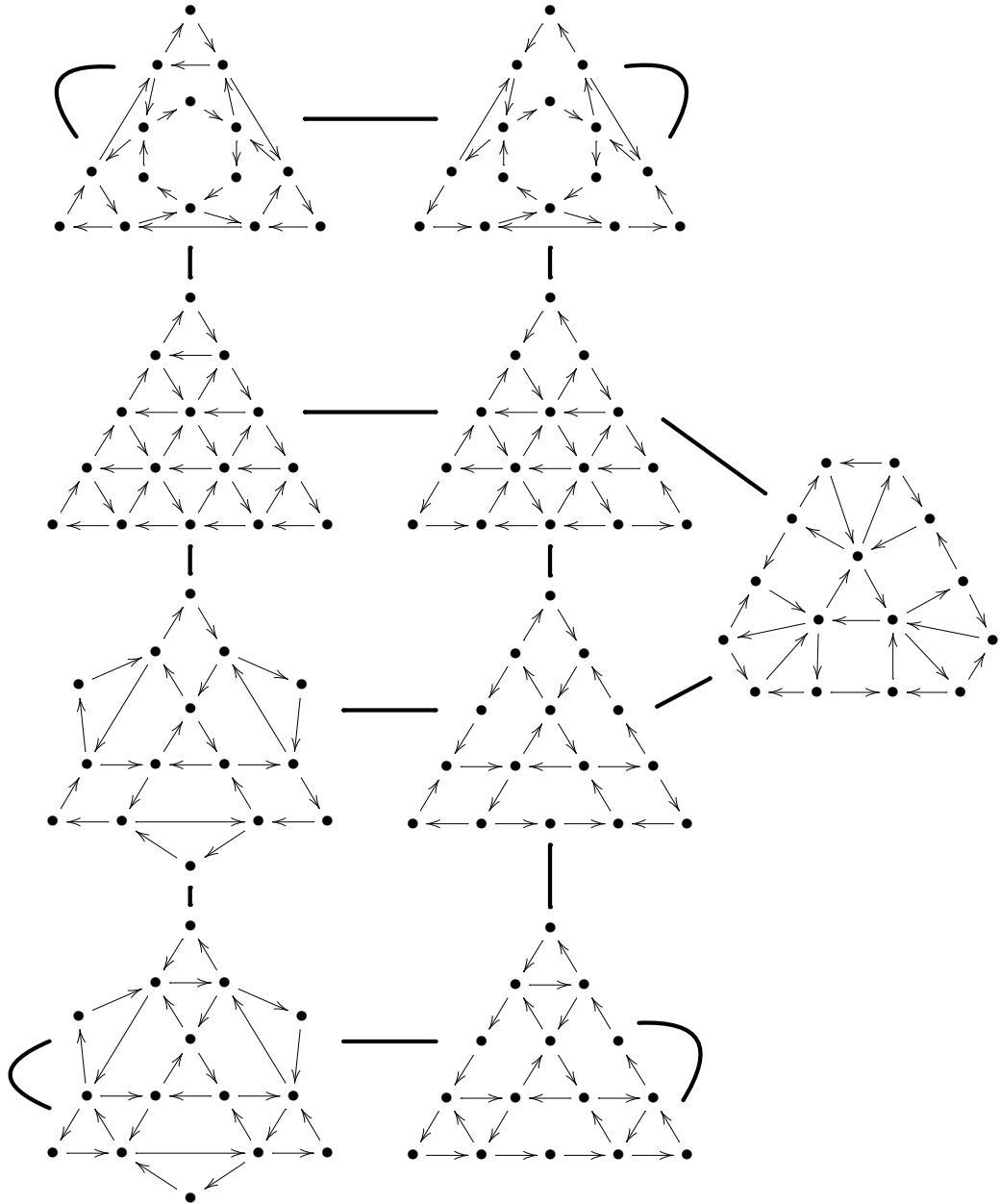

In contrast to the selfinjective cluster tilted algebras this family of planar QPs gives rise to many other by use of planar mutation. For instance the planar mutation lattices of $(Q^{(4)}, \epsilon^{(4)})$ and $(Q^{(5)}, \epsilon^{(5)})$ is displayed in figures \ref{mutation lattice A4} and \ref{mutation lattice A5} respectively.

\subsection{Examples from squares}\label{planar 4}
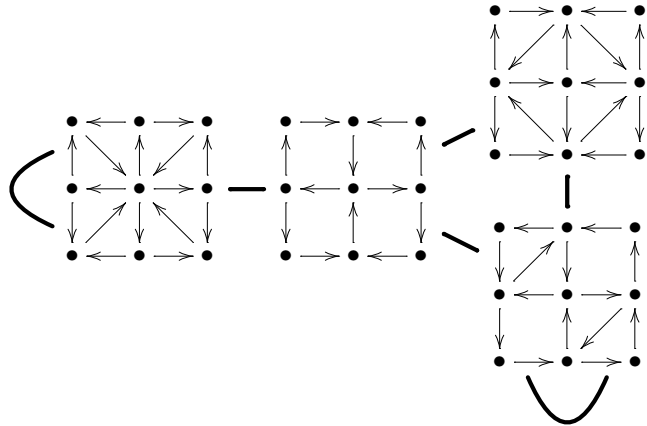
\begin{figure}[h]
\[
\begin{xy} 0;<0.8pt,0pt>:<0pt,-0.8pt>::
(100,0) *+{
\begin{xy} 0;<0.2pt,0pt>:<0pt,-0.2pt>:: 
(0,0) *+{\point 0} ="0",
(126,0) *+{\point 1} ="1",
(253,0) *+{\point 2} ="2",
(0,126) *+{\point 3} ="3",
(126,126) *+{\point 4} ="4",
(253,126) *+{\point 5} ="5",
(0,253) *+{\point 6} ="6",
(126,253) *+{\point 7} ="7",
(253,253) *+{\point 8} ="8",
"0", {\ar"1"},
"3", {\ar"0"},
"2", {\ar"1"},
"1", {\ar"4"},
"5", {\ar"2"},
"4", {\ar"3"},
"3", {\ar"6"},
"4", {\ar"5"},
"7", {\ar"4"},
"5", {\ar"8"},
"6", {\ar"7"},
"8", {\ar"7"},
\end{xy}
} ="e",
(200,50) *+{
\begin{xy} 0;<0.2pt,0pt>:<0pt,-0.2pt>:: 
(0,0) *+{\point 0} ="0",
(126,0) *+{\point 1} ="1",
(253,0) *+{\point 2} ="2",
(0,126) *+{\point 3} ="3",
(126,126) *+{\point 4} ="4",
(253,126) *+{\point 5} ="5",
(0,253) *+{\point 6} ="6",
(126,253) *+{\point 7} ="7",
(253,253) *+{\point 8} ="8",
"1", {\ar"0"},
"0", {\ar"3"},
"2", {\ar"1"},
"3", {\ar"1"},
"1", {\ar"4"},
"5", {\ar"2"},
"4", {\ar"3"},
"3", {\ar"6"},
"4", {\ar"5"},
"7", {\ar"4"},
"5", {\ar"7"},
"8", {\ar"5"},
"6", {\ar"7"},
"7", {\ar"8"},
\end{xy}
} ="1",
(0,0) *+{
\begin{xy} 0;<0.2pt,0pt>:<0pt,-0.2pt>:: 
(0,0) *+{\point 0} ="0",
(126,0) *+{\point 1} ="1",
(253,0) *+{\point 2} ="2",
(0,126) *+{\point 3} ="3",
(126,126) *+{\point 4} ="4",
(253,126) *+{\point 5} ="5",
(0,253) *+{\point 6} ="6",
(126,253) *+{\point 7} ="7",
(253,253) *+{\point 8} ="8",
"1", {\ar"0"},
"3", {\ar"0"},
"0", {\ar"4"},
"1", {\ar"2"},
"4", {\ar"1"},
"2", {\ar"4"},
"5", {\ar"2"},
"4", {\ar"3"},
"3", {\ar"6"},
"4", {\ar"5"},
"6", {\ar"4"},
"4", {\ar"7"},
"8", {\ar"4"},
"5", {\ar"8"},
"7", {\ar"6"},
"7", {\ar"8"},
\end{xy}
} ="2",
(200,-50) *+{
\begin{xy} 0;<0.2pt,0pt>:<0pt,-0.2pt>:: 
(0,0) *+{\point 0} ="0",
(135,0) *+{\point 1} ="1",
(271,0) *+{\point 2} ="2",
(0,135) *+{\point 3} ="3",
(135,135) *+{\point 4} ="4",
(271,135) *+{\point 5} ="5",
(0,271) *+{\point 6} ="6",
(135,271) *+{\point 7} ="7",
(271,271) *+{\point 8} ="8",
"0", {\ar"1"},
"3", {\ar"0"},
"2", {\ar"1"},
"1", {\ar"3"},
"4", {\ar"1"},
"1", {\ar"5"},
"5", {\ar"2"},
"3", {\ar"4"},
"3", {\ar"6"},
"7", {\ar"3"},
"5", {\ar"4"},
"4", {\ar"7"},
"7", {\ar"5"},
"5", {\ar"8"},
"6", {\ar"7"},
"8", {\ar"7"},
\end{xy}
} ="5",
"e", {\lat"1"},
"e", {\lat"2"},
"e", {\lat"5"},
"1", {\lat"5"},
(-35,-20) *+{\mbox{ }}="221",
(-35,20) *+{\mbox{ }}="222",
"221",{\lat@/_20pt/ "222"},
(220,85) *+{\mbox{ }}="111",
(180,85) *+{\mbox{ }}="112",
"111",{\lat@/^20pt/ "112"},
\end{xy}
\]
\caption{Planar mutation lattice for $(\mathbb{A}_3\square \mathbb{A}_3, W^{\square}_{\mathbb{A}_3, \mathbb{A}_3})$.}
\label{mutation lattice A3sq}
\end{figure}

Let $Q$ be a Dynkin quiver of type $\mathbb{A}_s$ with alternating orientation. By Theorem \ref{square}, $(Q\square Q, W^{\square}_{Q, Q})$ is selfinjective. Moreover it is planar. The planar mutation lattice of $(Q\square Q, W^{\square}_{Q, Q})$, when $s = 3$ is found in Figure \ref{mutation lattice A3sq}.

We proceed to define a class of planar QPs that we call square shaped. Every planar QP appearing in Figure \ref{mutation lattice A3sq} lies in this class. Moreover, so does $(Q^1\tildetensor Q^2, W^{\tildetensor}_{Q^1,Q^2})$ and $(Q^1\square Q^2, W^{\square}_{Q^1, Q^2})$ for all $Q^1$, $Q^2$ Dynkin of type $\mathbb{A}_s$, for some $s \in \N$. 

Let $Q$ be a quiver without loops constructed in the following way.
\begin{itemize}
\item The vertex set $Q_0 = \{1, \ldots, s\}^2\subset \R^2$. 
\item For each $1 \leq i,j < s$, let $P$ be the full subquiver of $Q$ with vertex set 
\[
P_0 =\{a=(i,j),b=(i+1,j),c=(i,j+1),d=(i+1,j+1)\}
\]
Then $P$ or $P^{\rm op}$ appears in the list below.
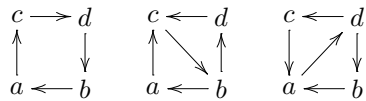
\[
\xymatrix@C0.5cm@R0.5cm{
c \ar[r]& d \ar[d] & c \ar[dr]& d \ar[l] & c \ar[d]& d \ar[l] \ar[d]\\
a \ar[u]& b \ar[l]& a \ar[u]& b \ar[u]\ar[l] & a \ar[ur]& b\ar[l]
}
\]
\end{itemize}
Let $(Q,\epsilon)$ be the planar QP arising from the pictures of $P$ and $P^{\rm op}$ above. We call such a QP square shaped. Here is one example:
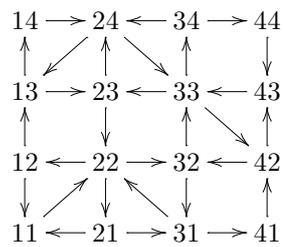
\[
\xymatrix@C0.5cm@R0.5cm{14 \ar[r]& 24 \ar[dl]\ar[dr]& 34 \ar[l]\ar[r]& 44 \ar[d]\\
13 \ar[u] \ar[r]& 23 \ar[u] \ar[d]& 33 \ar[l]\ar[u]\ar[dr]& 43 \ar[l]\\
12 \ar[u] \ar[d]& 22 \ar[l]\ar[r]\ar[d]& 32 \ar[u]\ar[d]& 42 \ar[u]\ar[l]\\ 
11 \ar[ur]& 21 \ar[l]\ar[r] & 31 \ar[ul]\ar[r] &41 \ar[u]\\}
\]

Define $\sigma : Q_0 \rightarrow Q_0$ by $\sigma(i,j) = (s-i+1,s-j+1)$. We call $(Q,\epsilon)$ symmetric if it has an automorphism that coincides with $\sigma$ on $Q_0$.

\begin{theorem}\label{selfinj square shaped}
Every symmetric square shaped QP is selfinjective with Nakayama permutation $\sigma$.
\end{theorem}

\begin{proof}
Let $(Q,\epsilon)$ be a symmetric square shaped QP. We regard $X = X_{(Q,\epsilon)}$ as a filled square in the plane. For each subquiver $P$ of $Q$ we also regard $X_P$ as a subset of $X$.

Consider the partition
\[
Q_0 = Q_0^0 \coprod Q_0^1,
\]
where $Q_0^k$ consists of all vertices $(i,j)$ such that $i+j \equiv k \;\; \mod 2$. Define the diagonal arrows of $Q$ to be $D_1 = D_0^0 \coprod D_0^1$, where $D_1^k$ is the set of arrows $a$ with $s(a),e(a) \in Q_0^k$. 
For each $k \in \{0,1\}$ set $D_0^k = s(D_1^k) \cup e(D_1^k)$ and let $D^k$ be the subquiver of $Q$ with vertices $D_0^k$ and arrows $D_1^k$. Denote the union of $D^0$ and $D^1$ by $D$.

We call the vertices in $Q_0 \setminus D_0$ mutable. Mutation at a mutable vertex is planar and produces a new square shaped QP. Let us make this observation more precise. Let $x$ be mutable and $y,z \in Q_0$ be of the form $x \pm (1,0)$ and $x \pm (0,1)$ respectively. 
Moreover, let $d^Q_{yz} \in \{0,1\}$ be the number of diagonal arrows between $y$ and $z$ in $Q$. Then $d^{\mu_xQ}_{yz} = 1-d^Q_{yz}$. Notice that if $x \in Q_0^k$, then the diagonal arrows affected by mutation at $x$ lie in $D^{k+1}_1$.

If $D$ is empty, then $(Q,\epsilon)$ is the square product of two Dynkin quivers of type $\mathbb{A}_s$. By Theorem \ref{square}, $(Q,\epsilon)$ is selfinjective with Nakayama permutation $\sigma$. Now assume that $D$ is not empty. 
We proceed to show that there is a set of mutable vertices in $Q$ invariant under $\sigma$ such that after mutation at these vertices the square shaped QP obtained will have fewer diagonal arrows than $Q$. This will complete the proof.

A subquiver $B$ of $D$ is called a border if it is either of type $\tilde{\mathbb{A}}$ with at most one vertex on the boundary $\partial X$ of $X$, or of type $\mathbb{A}$ with exactly two vertices on $\partial X$. 
Starting with a diagonal arrow we may extend it in both directions until a border is achieved. Thus $D$ is the union of all borders.

Let $B$ be a border. By possibly attaching a connected part of the $\partial X$ to $X_B$ we get a simple closed curve $c$. 
By the Jordan curve theorem, $\R^2 \setminus c$ decomposes into two connected components $A$ and $U$, where $A$ is bounded and $U$ is unbounded. 
Since $U$ intersects $\R^2 \setminus X$ we have $\R^2 \setminus X \subset U$ and so $A \subset X$. Moreover, for every point $x$ in $U\cap X$ there is a path in $U$ connecting $x$ with some point in $\R^2 \setminus X$. 
This path intersects $U \cap \partial X$ which is connected. Hence $U\cap X$ is connected. It follows that $X \setminus X_B$ has exactly two connected components. 
We denote these components by $Y^1_B$ and $Y^2_B$, where $x_0 \in Y^2_B$ for some fixed point $x_0 \in X \setminus X_Q$ not depending on $B$. 
Since $c = \partial A = \partial U$ we may recover $B$ from either $Y^1_B$ or $Y^2_B$.

Let $B$ and $B'$ be borders. We claim that $X_B \subset Y^1_{B'} \cup X_{B'}$ holds if and only if $Y^2_{B'} \subset Y^2_B$. 
If $X_B \subset Y^1_{B'} \cup X_{B'}$, then $Y^2_{B'} = X \setminus(Y^1_{B'} \cup X_{B'}) \subset X \setminus X_B$. 
Since $Y^2_{B}$ is a connected component of $X \setminus X_B$ intersecting $Y^2_{B'}$ at $x_0$ we have $Y^2_{B'} \subset Y^2_{B}$, as $Y^2_{B'}$ is connected. 
On the other hand, if $Y^2_{B'} \subset Y^2_{B}$, then $X_B \subset X \setminus Y^2_{B'} = Y^1_{B'} \cup X_{B'}$. 
Hence we get a partial order $B \leq B'$ defined by $X_B \subset Y^1_{B'} \cup X_{B'}$.

Now let $B$ be minimal border with respect to $\leq$. By symmetry we may assume that $B$ is a subquiver of $D^0$. Set $V := Q^1_0\cap Y_B^1$. We claim that all vertices in $V$ are mutable. 
Assume otherwise. Then there is an arrow $a \in D^1_1$ such that $s(a),e(a) \in V$. As explained earlier we can extend $a$ to a border $B'$ in $D^1$. 
Since $X_{B'}$ is connected and does not intersect $X_B$ we have that $X_{B'} \subset Y_B^1$, which contradicts the minimality of $B$.

Mutating at $V$ will remove the arrows in $B$ and will not create or remove any other diagonal arrows. 
Similarly mutating at $\sigma V$ will only affect the arrows in $\sigma B$. If $B = \sigma B$, then we mutate at $\sigma V = V$ to reduce the number of diagonal arrows. 
Otherwise we mutate at $(V \cup \sigma V) \setminus (V \cap \sigma V)$ to remove the arrows in $(B_1 \cup \sigma B_1) \setminus (B_1 \cap \sigma B_1)$.
\end{proof}

\begin{figure}[h]
\[
\begin{array}{cccc}
\begin{xy} 0;<0.4275pt,0pt>:<0pt,0.4275pt>:: 
(50,50) *+{\point{0}} ="0",
(100,50) *+{\point{1}} ="1",
(0,0) *+{\point{2}} ="2",
(100,0) *+{\point{3}} ="3",
(0,100) *+{\point{4}} ="4",
(50,100) *+{\point{5}} ="5",
(150,50) *+{\point{6}} ="6",
(200,0) *+{\point{7}} ="7",
(0,200) *+{\point{8}} ="8",
(50,150) *+{\point{9}} ="9",
(150,100) *+{\point{10}} ="10",
(200,100) *+{\point{11}} ="11",
(100,200) *+{\point{12}} ="12",
(200,200) *+{\point{13}} ="13",
(100,150) *+{\point{14}} ="14",
(150,150) *+{\point{15}} ="15",
"1", {\ar"0"},
"0", {\ar"2"},
"4", {\ar"0"},
"0", {\ar"5"},
"6", {\ar"1"},
"2", {\ar"3"},
"2", {\ar"4"},
"3", {\ar"7"},
"4", {\ar"8"},
"9", {\ar"4"},
"5", {\ar"9"},
"7", {\ar"6"},
"10", {\ar"6"},
"6", {\ar"11"},
"11", {\ar"7"},
"8", {\ar"9"},
"12", {\ar"8"},
"9", {\ar"14"},
"15", {\ar"10"},
"13", {\ar"11"},
"11", {\ar"15"},
"13", {\ar"12"},
"15", {\ar"13"},
"14", {\ar"15"},
\end{xy}
&
\begin{xy} 0;<0.4275pt,0pt>:<0pt,0.4275pt>:: 
(50,50) *+{\point{0}} ="0",
(100,50) *+{\point{1}} ="1",
(0,0) *+{\point{2}} ="2",
(100,0) *+{\point{3}} ="3",
(25,100) *+{\point{4}} ="4",
(75,100) *+{\point{5}} ="5",
(150,50) *+{\point{6}} ="6",
(200,0) *+{\point{7}} ="7",
(0,200) *+{\point{8}} ="8",
(50,150) *+{\point{9}} ="9",
(125,100) *+{\point{10}} ="10",
(175,100) *+{\point{11}} ="11",
(100,200) *+{\point{12}} ="12",
(200,200) *+{\point{13}} ="13",
(100,150) *+{\point{14}} ="14",
(150,150) *+{\point{15}} ="15",
"1", {\ar"0"},
"0", {\ar"4"},
"0", {\ar"5"},
"9", {\ar"0"},
"6", {\ar"1"},
"2", {\ar"3"},
"4", {\ar"2"},
"2", {\ar"8"},
"3", {\ar"7"},
"8", {\ar"4"},
"4", {\ar"9"},
"5", {\ar"9"},
"10", {\ar"6"},
"11", {\ar"6"},
"6", {\ar"15"},
"7", {\ar"11"},
"13", {\ar"7"},
"12", {\ar"8"},
"9", {\ar"14"},
"15", {\ar"10"},
"11", {\ar"13"},
"15", {\ar"11"},
"13", {\ar"12"},
"14", {\ar"15"},
\end{xy}
&
\begin{xy} 0;<0.4275pt,0pt>:<0pt,0.4275pt>:: 
(50,50) *+{\point{0}} ="0",
(100,50) *+{\point{1}} ="1",
(0,25) *+{\point{2}} ="2",
(100,0) *+{\point{3}} ="3",
(0,100) *+{\point{4}} ="4",
(50,100) *+{\point{5}} ="5",
(150,50) *+{\point{6}} ="6",
(200,25) *+{\point{7}} ="7",
(0,175) *+{\point{8}} ="8",
(50,150) *+{\point{9}} ="9",
(150,100) *+{\point{10}} ="10",
(200,100) *+{\point{11}} ="11",
(100,200) *+{\point{12}} ="12",
(200,175) *+{\point{13}} ="13",
(100,150) *+{\point{14}} ="14",
(150,150) *+{\point{15}} ="15",
"1", {\ar"0"},
"0", {\ar"2"},
"4", {\ar"0"},
"0", {\ar"5"},
"6", {\ar"1"},
"3", {\ar"2"},
"2", {\ar"4"},
"2", {\ar"7"},
"7", {\ar"3"},
"4", {\ar"8"},
"9", {\ar"4"},
"5", {\ar"9"},
"7", {\ar"6"},
"10", {\ar"6"},
"6", {\ar"11"},
"11", {\ar"7"},
"8", {\ar"9"},
"8", {\ar"12"},
"13", {\ar"8"},
"9", {\ar"14"},
"15", {\ar"10"},
"13", {\ar"11"},
"11", {\ar"15"},
"12", {\ar"13"},
"15", {\ar"13"},
"14", {\ar"15"},
\end{xy}
&
\begin{xy} 0;<0.4275pt,0pt>:<0pt,0.4275pt>:: 
(50,50) *+{\point{0}} ="0",
(100,50) *+{\point{1}} ="1",
(0,25) *+{\point{2}} ="2",
(100,0) *+{\point{3}} ="3",
(25,100) *+{\point{4}} ="4",
(75,100) *+{\point{5}} ="5",
(150,50) *+{\point{6}} ="6",
(200,25) *+{\point{7}} ="7",
(0,175) *+{\point{8}} ="8",
(50,150) *+{\point{9}} ="9",
(125,100) *+{\point{10}} ="10",
(175,100) *+{\point{11}} ="11",
(100,200) *+{\point{12}} ="12",
(200,175) *+{\point{13}} ="13",
(100,150) *+{\point{14}} ="14",
(150,150) *+{\point{15}} ="15",
"1", {\ar"0"},
"0", {\ar"4"},
"0", {\ar"5"},
"9", {\ar"0"},
"6", {\ar"1"},
"3", {\ar"2"},
"4", {\ar"2"},
"2", {\ar"7"},
"2", {\ar"8"},
"7", {\ar"3"},
"8", {\ar"4"},
"4", {\ar"9"},
"5", {\ar"9"},
"10", {\ar"6"},
"11", {\ar"6"},
"6", {\ar"15"},
"7", {\ar"11"},
"13", {\ar"7"},
"8", {\ar"12"},
"13", {\ar"8"},
"9", {\ar"14"},
"15", {\ar"10"},
"11", {\ar"13"},
"15", {\ar"11"},
"12", {\ar"13"},
"14", {\ar"15"},
\end{xy}
\end{array}
\]
\caption{Planar mutation lattice for $s=4$.}
\label{mutation lattice2}
\end{figure}
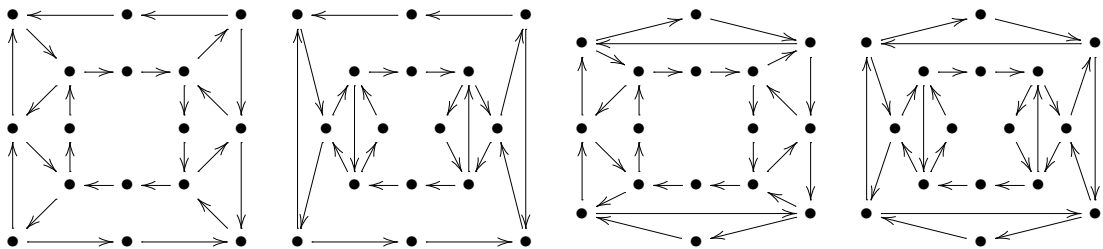

In the above proof we saw that symmetric square shaped QPs of fixed size $s$ are obtained from each other via planar mutation. 
However, for $s \geq 4$, many other planar QPs can be constructed by planar mutation as well. 
For instance the planar mutation lattice for the case $s=4$ is displayed in Figure \ref{mutation lattice2}.

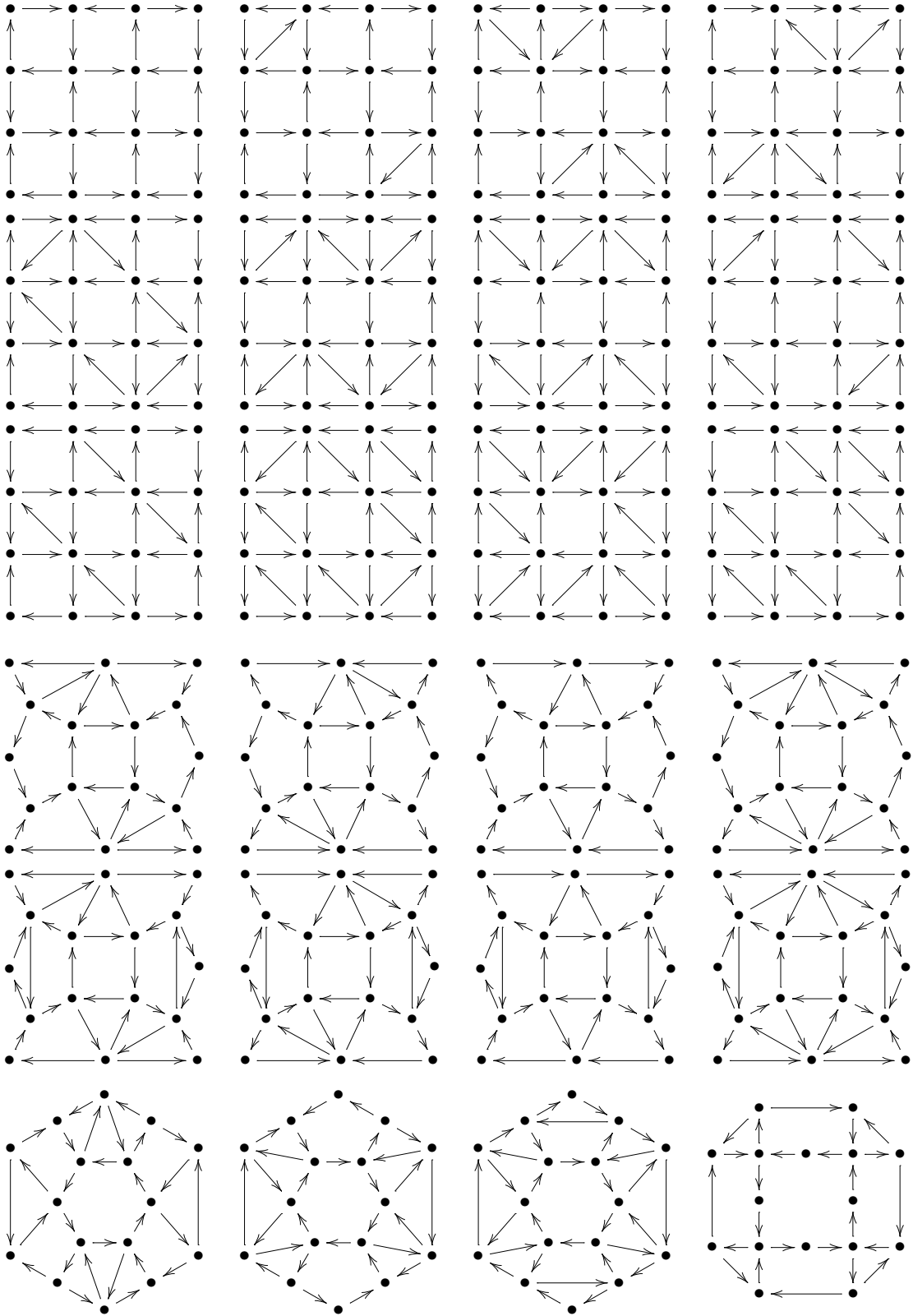
\begin{figure}
\[
\begin{array}{cccc}
\begin{xy} 0;<0.38pt,0pt>:<0pt,-0.38pt>:: 
(0,0) *+{\point{0}} ="0",
(75,0) *+{\point{1}} ="1",
(150,0) *+{\point{2}} ="2",
(225,0) *+{\point{3}} ="3",
(0,75) *+{\point{4}} ="4",
(75,75) *+{\point{5}} ="5",
(150,75) *+{\point{6}} ="6",
(225,75) *+{\point{7}} ="7",
(0,150) *+{\point{8}} ="8",
(75,150) *+{\point{9}} ="9",
(150,150) *+{\point{10}} ="10",
(225,150) *+{\point{11}} ="11",
(0,225) *+{\point{12}} ="12",
(75,225) *+{\point{13}} ="13",
(150,225) *+{\point{14}} ="14",
(225,225) *+{\point{15}} ="15",
"0", {\ar"1"},
"4", {\ar"0"},
"2", {\ar"1"},
"1", {\ar"5"},
"2", {\ar"3"},
"6", {\ar"2"},
"3", {\ar"7"},
"5", {\ar"4"},
"4", {\ar"8"},
"5", {\ar"6"},
"9", {\ar"5"},
"7", {\ar"6"},
"6", {\ar"10"},
"11", {\ar"7"},
"8", {\ar"9"},
"12", {\ar"8"},
"10", {\ar"9"},
"9", {\ar"13"},
"10", {\ar"11"},
"14", {\ar"10"},
"11", {\ar"15"},
"13", {\ar"12"},
"13", {\ar"14"},
"15", {\ar"14"},
\end{xy}
&
\begin{xy} 0;<0.38pt,0pt>:<0pt,-0.38pt>:: 
(0,0) *+{\point{0}} ="0",
(75,0) *+{\point{1}} ="1",
(150,0) *+{\point{2}} ="2",
(225,0) *+{\point{3}} ="3",
(0,75) *+{\point{4}} ="4",
(75,75) *+{\point{5}} ="5",
(150,75) *+{\point{6}} ="6",
(225,75) *+{\point{7}} ="7",
(0,150) *+{\point{8}} ="8",
(75,150) *+{\point{9}} ="9",
(150,150) *+{\point{10}} ="10",
(225,150) *+{\point{11}} ="11",
(0,225) *+{\point{12}} ="12",
(75,225) *+{\point{13}} ="13",
(150,225) *+{\point{14}} ="14",
(225,225) *+{\point{15}} ="15",
"1", {\ar"0"},
"0", {\ar"4"},
"2", {\ar"1"},
"4", {\ar"1"},
"1", {\ar"5"},
"2", {\ar"3"},
"6", {\ar"2"},
"3", {\ar"7"},
"5", {\ar"4"},
"4", {\ar"8"},
"5", {\ar"6"},
"9", {\ar"5"},
"7", {\ar"6"},
"6", {\ar"10"},
"11", {\ar"7"},
"8", {\ar"9"},
"12", {\ar"8"},
"10", {\ar"9"},
"9", {\ar"13"},
"10", {\ar"11"},
"14", {\ar"10"},
"11", {\ar"14"},
"15", {\ar"11"},
"13", {\ar"12"},
"13", {\ar"14"},
"14", {\ar"15"},
\end{xy}
&
\begin{xy} 0;<0.38pt,0pt>:<0pt,-0.38pt>:: 
(0,0) *+{\point{0}} ="0",
(75,0) *+{\point{1}} ="1",
(150,0) *+{\point{2}} ="2",
(225,0) *+{\point{3}} ="3",
(0,75) *+{\point{4}} ="4",
(75,75) *+{\point{5}} ="5",
(150,75) *+{\point{6}} ="6",
(225,75) *+{\point{7}} ="7",
(0,150) *+{\point{8}} ="8",
(75,150) *+{\point{9}} ="9",
(150,150) *+{\point{10}} ="10",
(225,150) *+{\point{11}} ="11",
(0,225) *+{\point{12}} ="12",
(75,225) *+{\point{13}} ="13",
(150,225) *+{\point{14}} ="14",
(225,225) *+{\point{15}} ="15",
"1", {\ar"0"},
"4", {\ar"0"},
"0", {\ar"5"},
"1", {\ar"2"},
"5", {\ar"1"},
"2", {\ar"3"},
"2", {\ar"5"},
"6", {\ar"2"},
"3", {\ar"7"},
"5", {\ar"4"},
"4", {\ar"8"},
"5", {\ar"6"},
"9", {\ar"5"},
"7", {\ar"6"},
"6", {\ar"10"},
"11", {\ar"7"},
"8", {\ar"9"},
"12", {\ar"8"},
"10", {\ar"9"},
"9", {\ar"13"},
"10", {\ar"11"},
"13", {\ar"10"},
"10", {\ar"14"},
"15", {\ar"10"},
"11", {\ar"15"},
"13", {\ar"12"},
"14", {\ar"13"},
"14", {\ar"15"},
\end{xy}
&
\begin{xy} 0;<0.38pt,0pt>:<0pt,-0.38pt>:: 
(0,0) *+{\point{0}} ="0",
(75,0) *+{\point{1}} ="1",
(150,0) *+{\point{2}} ="2",
(225,0) *+{\point{3}} ="3",
(0,75) *+{\point{4}} ="4",
(75,75) *+{\point{5}} ="5",
(150,75) *+{\point{6}} ="6",
(225,75) *+{\point{7}} ="7",
(0,150) *+{\point{8}} ="8",
(75,150) *+{\point{9}} ="9",
(150,150) *+{\point{10}} ="10",
(225,150) *+{\point{11}} ="11",
(0,225) *+{\point{12}} ="12",
(75,225) *+{\point{13}} ="13",
(150,225) *+{\point{14}} ="14",
(225,225) *+{\point{15}} ="15",
"0", {\ar"1"},
"4", {\ar"0"},
"1", {\ar"2"},
"1", {\ar"5"},
"6", {\ar"1"},
"3", {\ar"2"},
"2", {\ar"6"},
"6", {\ar"3"},
"3", {\ar"7"},
"5", {\ar"4"},
"4", {\ar"8"},
"5", {\ar"6"},
"9", {\ar"5"},
"7", {\ar"6"},
"6", {\ar"10"},
"11", {\ar"7"},
"8", {\ar"9"},
"12", {\ar"8"},
"10", {\ar"9"},
"9", {\ar"12"},
"13", {\ar"9"},
"9", {\ar"14"},
"10", {\ar"11"},
"14", {\ar"10"},
"11", {\ar"15"},
"12", {\ar"13"},
"14", {\ar"13"},
"15", {\ar"14"},
\end{xy}
\\
\begin{xy} 0;<0.38pt,0pt>:<0pt,-0.38pt>:: 
(0,0) *+{\point{0}} ="0",
(75,0) *+{\point{1}} ="1",
(150,0) *+{\point{2}} ="2",
(225,0) *+{\point{3}} ="3",
(0,75) *+{\point{4}} ="4",
(75,75) *+{\point{5}} ="5",
(150,75) *+{\point{6}} ="6",
(225,75) *+{\point{7}} ="7",
(0,150) *+{\point{8}} ="8",
(75,150) *+{\point{9}} ="9",
(150,150) *+{\point{10}} ="10",
(225,150) *+{\point{11}} ="11",
(0,225) *+{\point{12}} ="12",
(75,225) *+{\point{13}} ="13",
(150,225) *+{\point{14}} ="14",
(225,225) *+{\point{15}} ="15",
"0", {\ar"1"},
"4", {\ar"0"},
"2", {\ar"1"},
"1", {\ar"4"},
"5", {\ar"1"},
"1", {\ar"6"},
"2", {\ar"3"},
"6", {\ar"2"},
"3", {\ar"7"},
"4", {\ar"5"},
"4", {\ar"8"},
"9", {\ar"4"},
"6", {\ar"5"},
"5", {\ar"9"},
"7", {\ar"6"},
"10", {\ar"6"},
"6", {\ar"11"},
"11", {\ar"7"},
"8", {\ar"9"},
"12", {\ar"8"},
"9", {\ar"10"},
"9", {\ar"13"},
"14", {\ar"9"},
"11", {\ar"10"},
"10", {\ar"14"},
"14", {\ar"11"},
"11", {\ar"15"},
"13", {\ar"12"},
"13", {\ar"14"},
"15", {\ar"14"},
\end{xy}
&
\begin{xy} 0;<0.38pt,0pt>:<0pt,-0.38pt>:: 
(0,0) *+{\point{0}} ="0",
(75,0) *+{\point{1}} ="1",
(150,0) *+{\point{2}} ="2",
(225,0) *+{\point{3}} ="3",
(0,75) *+{\point{4}} ="4",
(75,75) *+{\point{5}} ="5",
(150,75) *+{\point{6}} ="6",
(225,75) *+{\point{7}} ="7",
(0,150) *+{\point{8}} ="8",
(75,150) *+{\point{9}} ="9",
(150,150) *+{\point{10}} ="10",
(225,150) *+{\point{11}} ="11",
(0,225) *+{\point{12}} ="12",
(75,225) *+{\point{13}} ="13",
(150,225) *+{\point{14}} ="14",
(225,225) *+{\point{15}} ="15",
"1", {\ar"0"},
"0", {\ar"4"},
"1", {\ar"2"},
"4", {\ar"1"},
"1", {\ar"5"},
"6", {\ar"1"},
"3", {\ar"2"},
"2", {\ar"6"},
"6", {\ar"3"},
"3", {\ar"7"},
"5", {\ar"4"},
"4", {\ar"8"},
"5", {\ar"6"},
"9", {\ar"5"},
"7", {\ar"6"},
"6", {\ar"10"},
"11", {\ar"7"},
"8", {\ar"9"},
"12", {\ar"8"},
"10", {\ar"9"},
"9", {\ar"12"},
"13", {\ar"9"},
"9", {\ar"14"},
"10", {\ar"11"},
"14", {\ar"10"},
"11", {\ar"14"},
"15", {\ar"11"},
"12", {\ar"13"},
"14", {\ar"13"},
"14", {\ar"15"},
\end{xy}
&
\begin{xy} 0;<0.38pt,0pt>:<0pt,-0.38pt>:: 
(0,0) *+{\point{0}} ="0",
(75,0) *+{\point{1}} ="1",
(150,0) *+{\point{2}} ="2",
(225,0) *+{\point{3}} ="3",
(0,75) *+{\point{4}} ="4",
(75,75) *+{\point{5}} ="5",
(150,75) *+{\point{6}} ="6",
(225,75) *+{\point{7}} ="7",
(0,150) *+{\point{8}} ="8",
(75,150) *+{\point{9}} ="9",
(150,150) *+{\point{10}} ="10",
(225,150) *+{\point{11}} ="11",
(0,225) *+{\point{12}} ="12",
(75,225) *+{\point{13}} ="13",
(150,225) *+{\point{14}} ="14",
(225,225) *+{\point{15}} ="15",
"1", {\ar"0"},
"4", {\ar"0"},
"0", {\ar"5"},
"1", {\ar"2"},
"5", {\ar"1"},
"3", {\ar"2"},
"2", {\ar"5"},
"6", {\ar"2"},
"2", {\ar"7"},
"7", {\ar"3"},
"5", {\ar"4"},
"4", {\ar"8"},
"5", {\ar"6"},
"9", {\ar"5"},
"7", {\ar"6"},
"6", {\ar"10"},
"11", {\ar"7"},
"8", {\ar"9"},
"8", {\ar"12"},
"13", {\ar"8"},
"10", {\ar"9"},
"9", {\ar"13"},
"10", {\ar"11"},
"13", {\ar"10"},
"10", {\ar"14"},
"15", {\ar"10"},
"11", {\ar"15"},
"12", {\ar"13"},
"14", {\ar"13"},
"14", {\ar"15"},
\end{xy}
&
\begin{xy} 0;<0.38pt,0pt>:<0pt,-0.38pt>:: 
(0,0) *+{\point{0}} ="0",
(75,0) *+{\point{1}} ="1",
(150,0) *+{\point{2}} ="2",
(225,0) *+{\point{3}} ="3",
(0,75) *+{\point{4}} ="4",
(75,75) *+{\point{5}} ="5",
(150,75) *+{\point{6}} ="6",
(225,75) *+{\point{7}} ="7",
(0,150) *+{\point{8}} ="8",
(75,150) *+{\point{9}} ="9",
(150,150) *+{\point{10}} ="10",
(225,150) *+{\point{11}} ="11",
(0,225) *+{\point{12}} ="12",
(75,225) *+{\point{13}} ="13",
(150,225) *+{\point{14}} ="14",
(225,225) *+{\point{15}} ="15",
"1", {\ar"0"},
"0", {\ar"4"},
"2", {\ar"1"},
"4", {\ar"1"},
"1", {\ar"5"},
"3", {\ar"2"},
"6", {\ar"2"},
"2", {\ar"7"},
"7", {\ar"3"},
"5", {\ar"4"},
"4", {\ar"8"},
"5", {\ar"6"},
"9", {\ar"5"},
"7", {\ar"6"},
"6", {\ar"10"},
"11", {\ar"7"},
"8", {\ar"9"},
"8", {\ar"12"},
"13", {\ar"8"},
"10", {\ar"9"},
"9", {\ar"13"},
"10", {\ar"11"},
"14", {\ar"10"},
"11", {\ar"14"},
"15", {\ar"11"},
"12", {\ar"13"},
"13", {\ar"14"},
"14", {\ar"15"},
\end{xy}
\\
\begin{xy} 0;<0.38pt,0pt>:<0pt,-0.38pt>:: 
(0,0) *+{\point{0}} ="0",
(75,0) *+{\point{1}} ="1",
(150,0) *+{\point{2}} ="2",
(225,0) *+{\point{3}} ="3",
(0,75) *+{\point{4}} ="4",
(75,75) *+{\point{5}} ="5",
(150,75) *+{\point{6}} ="6",
(225,75) *+{\point{7}} ="7",
(0,150) *+{\point{8}} ="8",
(75,150) *+{\point{9}} ="9",
(150,150) *+{\point{10}} ="10",
(225,150) *+{\point{11}} ="11",
(0,225) *+{\point{12}} ="12",
(75,225) *+{\point{13}} ="13",
(150,225) *+{\point{14}} ="14",
(225,225) *+{\point{15}} ="15",
"1", {\ar"0"},
"0", {\ar"4"},
"2", {\ar"1"},
"5", {\ar"1"},
"1", {\ar"6"},
"2", {\ar"3"},
"6", {\ar"2"},
"3", {\ar"7"},
"4", {\ar"5"},
"4", {\ar"8"},
"9", {\ar"4"},
"6", {\ar"5"},
"5", {\ar"9"},
"7", {\ar"6"},
"10", {\ar"6"},
"6", {\ar"11"},
"11", {\ar"7"},
"8", {\ar"9"},
"12", {\ar"8"},
"9", {\ar"10"},
"9", {\ar"13"},
"14", {\ar"9"},
"11", {\ar"10"},
"10", {\ar"14"},
"15", {\ar"11"},
"13", {\ar"12"},
"13", {\ar"14"},
"14", {\ar"15"},
\end{xy}
&
\begin{xy} 0;<0.38pt,0pt>:<0pt,-0.38pt>:: 
(0,0) *+{\point{0}} ="0",
(75,0) *+{\point{1}} ="1",
(150,0) *+{\point{2}} ="2",
(225,0) *+{\point{3}} ="3",
(0,75) *+{\point{4}} ="4",
(75,75) *+{\point{5}} ="5",
(150,75) *+{\point{6}} ="6",
(225,75) *+{\point{7}} ="7",
(0,150) *+{\point{8}} ="8",
(75,150) *+{\point{9}} ="9",
(150,150) *+{\point{10}} ="10",
(225,150) *+{\point{11}} ="11",
(0,225) *+{\point{12}} ="12",
(75,225) *+{\point{13}} ="13",
(150,225) *+{\point{14}} ="14",
(225,225) *+{\point{15}} ="15",
"0", {\ar"1"},
"4", {\ar"0"},
"2", {\ar"1"},
"1", {\ar"4"},
"5", {\ar"1"},
"1", {\ar"6"},
"3", {\ar"2"},
"6", {\ar"2"},
"2", {\ar"7"},
"7", {\ar"3"},
"4", {\ar"5"},
"4", {\ar"8"},
"9", {\ar"4"},
"6", {\ar"5"},
"5", {\ar"9"},
"7", {\ar"6"},
"10", {\ar"6"},
"6", {\ar"11"},
"11", {\ar"7"},
"8", {\ar"9"},
"8", {\ar"12"},
"13", {\ar"8"},
"9", {\ar"10"},
"9", {\ar"13"},
"14", {\ar"9"},
"11", {\ar"10"},
"10", {\ar"14"},
"14", {\ar"11"},
"11", {\ar"15"},
"12", {\ar"13"},
"13", {\ar"14"},
"15", {\ar"14"},
\end{xy}
&
\begin{xy} 0;<0.38pt,0pt>:<0pt,-0.38pt>:: 
(0,0) *+{\point{0}} ="0",
(75,0) *+{\point{1}} ="1",
(150,0) *+{\point{2}} ="2",
(225,0) *+{\point{3}} ="3",
(0,75) *+{\point{4}} ="4",
(75,75) *+{\point{5}} ="5",
(150,75) *+{\point{6}} ="6",
(225,75) *+{\point{7}} ="7",
(0,150) *+{\point{8}} ="8",
(75,150) *+{\point{9}} ="9",
(150,150) *+{\point{10}} ="10",
(225,150) *+{\point{11}} ="11",
(0,225) *+{\point{12}} ="12",
(75,225) *+{\point{13}} ="13",
(150,225) *+{\point{14}} ="14",
(225,225) *+{\point{15}} ="15",
"1", {\ar"0"},
"4", {\ar"0"},
"0", {\ar"5"},
"1", {\ar"2"},
"5", {\ar"1"},
"2", {\ar"3"},
"2", {\ar"5"},
"6", {\ar"2"},
"3", {\ar"6"},
"7", {\ar"3"},
"5", {\ar"4"},
"8", {\ar"4"},
"4", {\ar"9"},
"5", {\ar"6"},
"9", {\ar"5"},
"6", {\ar"7"},
"6", {\ar"10"},
"11", {\ar"6"},
"7", {\ar"11"},
"9", {\ar"8"},
"8", {\ar"12"},
"10", {\ar"9"},
"12", {\ar"9"},
"9", {\ar"13"},
"10", {\ar"11"},
"13", {\ar"10"},
"10", {\ar"14"},
"15", {\ar"10"},
"11", {\ar"15"},
"13", {\ar"12"},
"14", {\ar"13"},
"14", {\ar"15"},
\end{xy}
&
\begin{xy} 0;<0.38pt,0pt>:<0pt,-0.38pt>:: 
(0,0) *+{\point{0}} ="0",
(75,0) *+{\point{1}} ="1",
(150,0) *+{\point{2}} ="2",
(225,0) *+{\point{3}} ="3",
(0,75) *+{\point{4}} ="4",
(75,75) *+{\point{5}} ="5",
(150,75) *+{\point{6}} ="6",
(225,75) *+{\point{7}} ="7",
(0,150) *+{\point{8}} ="8",
(75,150) *+{\point{9}} ="9",
(150,150) *+{\point{10}} ="10",
(225,150) *+{\point{11}} ="11",
(0,225) *+{\point{12}} ="12",
(75,225) *+{\point{13}} ="13",
(150,225) *+{\point{14}} ="14",
(225,225) *+{\point{15}} ="15",
"1", {\ar"0"},
"0", {\ar"4"},
"2", {\ar"1"},
"5", {\ar"1"},
"1", {\ar"6"},
"3", {\ar"2"},
"6", {\ar"2"},
"2", {\ar"7"},
"7", {\ar"3"},
"4", {\ar"5"},
"4", {\ar"8"},
"9", {\ar"4"},
"6", {\ar"5"},
"5", {\ar"9"},
"7", {\ar"6"},
"10", {\ar"6"},
"6", {\ar"11"},
"11", {\ar"7"},
"8", {\ar"9"},
"8", {\ar"12"},
"13", {\ar"8"},
"9", {\ar"10"},
"9", {\ar"13"},
"14", {\ar"9"},
"11", {\ar"10"},
"10", {\ar"14"},
"15", {\ar"11"},
"12", {\ar"13"},
"13", {\ar"14"},
"14", {\ar"15"},
\end{xy}
\end{array}
\]

\[
\begin{array}{cccc}
\begin{xy} 0;<0.38pt,0pt>:<0pt,-0.38pt>:: 
(0,0) *+{\point{0}} ="0",
(25,50) *+{\point{1}} ="1",
(115,0) *+{\point{2}} ="2",
(225,0) *+{\point{3}} ="3",
(0,114) *+{\point{4}} ="4",
(75,75) *+{\point{5}} ="5",
(150,75) *+{\point{6}} ="6",
(200,50) *+{\point{7}} ="7",
(25,175) *+{\point{8}} ="8",
(75,150) *+{\point{9}} ="9",
(150,150) *+{\point{10}} ="10",
(227,112) *+{\point{11}} ="11",
(0,225) *+{\point{12}} ="12",
(115,225) *+{\point{13}} ="13",
(200,175) *+{\point{14}} ="14",
(225,225) *+{\point{15}} ="15",
"0", {\ar"1"},
"2", {\ar"0"},
"1", {\ar"2"},
"1", {\ar"4"},
"5", {\ar"1"},
"2", {\ar"3"},
"2", {\ar"5"},
"6", {\ar"2"},
"3", {\ar"7"},
"4", {\ar"8"},
"5", {\ar"6"},
"9", {\ar"5"},
"7", {\ar"6"},
"6", {\ar"10"},
"11", {\ar"7"},
"8", {\ar"9"},
"12", {\ar"8"},
"10", {\ar"9"},
"9", {\ar"13"},
"13", {\ar"10"},
"10", {\ar"14"},
"14", {\ar"11"},
"13", {\ar"12"},
"14", {\ar"13"},
"13", {\ar"15"},
"15", {\ar"14"},
\end{xy}
&
\begin{xy} 0;<0.38pt,0pt>:<0pt,-0.38pt>:: 
(0,0) *+{\point{0}} ="0",
(25,50) *+{\point{1}} ="1",
(115,0) *+{\point{2}} ="2",
(225,0) *+{\point{3}} ="3",
(0,114) *+{\point{4}} ="4",
(75,75) *+{\point{5}} ="5",
(150,75) *+{\point{6}} ="6",
(200,50) *+{\point{7}} ="7",
(25,175) *+{\point{8}} ="8",
(75,150) *+{\point{9}} ="9",
(150,150) *+{\point{10}} ="10",
(227,112) *+{\point{11}} ="11",
(0,225) *+{\point{12}} ="12",
(115,225) *+{\point{13}} ="13",
(200,175) *+{\point{14}} ="14",
(225,225) *+{\point{15}} ="15",
"1", {\ar"0"},
"0", {\ar"2"},
"1", {\ar"4"},
"5", {\ar"1"},
"3", {\ar"2"},
"2", {\ar"5"},
"6", {\ar"2"},
"2", {\ar"7"},
"7", {\ar"3"},
"4", {\ar"8"},
"5", {\ar"6"},
"9", {\ar"5"},
"7", {\ar"6"},
"6", {\ar"10"},
"11", {\ar"7"},
"8", {\ar"9"},
"8", {\ar"12"},
"13", {\ar"8"},
"10", {\ar"9"},
"9", {\ar"13"},
"13", {\ar"10"},
"10", {\ar"14"},
"14", {\ar"11"},
"12", {\ar"13"},
"15", {\ar"13"},
"14", {\ar"15"},
\end{xy}
&
\begin{xy} 0;<0.38pt,0pt>:<0pt,-0.38pt>:: 
(0,0) *+{\point{1}} ="1",
(25,50) *+{\point{0}} ="0",
(115,0) *+{\point{2}} ="2",
(225,0) *+{\point{3}} ="3",
(0,114) *+{\point{4}} ="4",
(75,75) *+{\point{5}} ="5",
(150,75) *+{\point{6}} ="6",
(200,50) *+{\point{7}} ="7",
(25,175) *+{\point{8}} ="8",
(75,150) *+{\point{9}} ="9",
(150,150) *+{\point{10}} ="10",
(227,112) *+{\point{11}} ="11",
(0,225) *+{\point{12}} ="12",
(115,225) *+{\point{13}} ="13",
(200,175) *+{\point{15}} ="15",
(225,225) *+{\point{14}} ="14",
"0", {\ar"1"},
"0", {\ar"4"},
"5", {\ar"0"},
"1", {\ar"2"},
"2", {\ar"3"},
"2", {\ar"5"},
"6", {\ar"2"},
"3", {\ar"7"},
"4", {\ar"8"},
"5", {\ar"6"},
"9", {\ar"5"},
"7", {\ar"6"},
"6", {\ar"10"},
"11", {\ar"7"},
"8", {\ar"9"},
"12", {\ar"8"},
"10", {\ar"9"},
"9", {\ar"13"},
"13", {\ar"10"},
"10", {\ar"15"},
"15", {\ar"11"},
"13", {\ar"12"},
"14", {\ar"13"},
"15", {\ar"14"},
\end{xy}
&
\begin{xy} 0;<0.38pt,0pt>:<0pt,-0.38pt>:: 
(0,0) *+{\point{0}} ="0",
(25,50) *+{\point{1}} ="1",
(115,0) *+{\point{2}} ="2",
(225,0) *+{\point{3}} ="3",
(0,114) *+{\point{4}} ="4",
(75,75) *+{\point{5}} ="5",
(150,75) *+{\point{6}} ="6",
(200,50) *+{\point{7}} ="7",
(25,175) *+{\point{8}} ="8",
(75,150) *+{\point{9}} ="9",
(150,150) *+{\point{10}} ="10",
(227,112) *+{\point{11}} ="11",
(0,225) *+{\point{12}} ="12",
(115,225) *+{\point{13}} ="13",
(200,175) *+{\point{14}} ="14",
(225,225) *+{\point{15}} ="15",
"0", {\ar"1"},
"2", {\ar"0"},
"1", {\ar"2"},
"1", {\ar"4"},
"5", {\ar"1"},
"3", {\ar"2"},
"2", {\ar"5"},
"6", {\ar"2"},
"2", {\ar"7"},
"7", {\ar"3"},
"4", {\ar"8"},
"5", {\ar"6"},
"9", {\ar"5"},
"7", {\ar"6"},
"6", {\ar"10"},
"11", {\ar"7"},
"8", {\ar"9"},
"8", {\ar"12"},
"13", {\ar"8"},
"10", {\ar"9"},
"9", {\ar"13"},
"13", {\ar"10"},
"10", {\ar"14"},
"14", {\ar"11"},
"12", {\ar"13"},
"14", {\ar"13"},
"13", {\ar"15"},
"15", {\ar"14"},
\end{xy}
\\
\begin{xy} 0;<0.38pt,0pt>:<0pt,-0.38pt>:: 
(0,0) *+{\point{0}} ="0",
(25,50) *+{\point{1}} ="1",
(115,0) *+{\point{2}} ="2",
(225,0) *+{\point{3}} ="3",
(0,114) *+{\point{4}} ="4",
(75,75) *+{\point{5}} ="5",
(150,75) *+{\point{6}} ="6",
(200,50) *+{\point{7}} ="7",
(25,175) *+{\point{8}} ="8",
(75,150) *+{\point{9}} ="9",
(150,150) *+{\point{10}} ="10",
(227,112) *+{\point{11}} ="11",
(0,225) *+{\point{12}} ="12",
(115,225) *+{\point{13}} ="13",
(200,175) *+{\point{14}} ="14",
(225,225) *+{\point{15}} ="15",
"0", {\ar"1"},
"2", {\ar"0"},
"1", {\ar"2"},
"4", {\ar"1"},
"5", {\ar"1"},
"1", {\ar"8"},
"2", {\ar"3"},
"2", {\ar"5"},
"6", {\ar"2"},
"3", {\ar"7"},
"8", {\ar"4"},
"5", {\ar"6"},
"9", {\ar"5"},
"7", {\ar"6"},
"6", {\ar"10"},
"7", {\ar"11"},
"14", {\ar"7"},
"8", {\ar"9"},
"12", {\ar"8"},
"10", {\ar"9"},
"9", {\ar"13"},
"13", {\ar"10"},
"10", {\ar"14"},
"11", {\ar"14"},
"13", {\ar"12"},
"14", {\ar"13"},
"13", {\ar"15"},
"15", {\ar"14"},
\end{xy}
&
\begin{xy} 0;<0.38pt,0pt>:<0pt,-0.38pt>:: 
(0,0) *+{\point{0}} ="0",
(25,50) *+{\point{1}} ="1",
(115,0) *+{\point{2}} ="2",
(225,0) *+{\point{3}} ="3",
(0,114) *+{\point{4}} ="4",
(75,75) *+{\point{5}} ="5",
(150,75) *+{\point{6}} ="6",
(200,50) *+{\point{7}} ="7",
(25,175) *+{\point{8}} ="8",
(75,150) *+{\point{9}} ="9",
(150,150) *+{\point{10}} ="10",
(227,112) *+{\point{11}} ="11",
(0,225) *+{\point{12}} ="12",
(115,225) *+{\point{13}} ="13",
(200,175) *+{\point{14}} ="14",
(225,225) *+{\point{15}} ="15",
"1", {\ar"0"},
"0", {\ar"2"},
"4", {\ar"1"},
"5", {\ar"1"},
"1", {\ar"8"},
"3", {\ar"2"},
"2", {\ar"5"},
"6", {\ar"2"},
"2", {\ar"7"},
"7", {\ar"3"},
"8", {\ar"4"},
"5", {\ar"6"},
"9", {\ar"5"},
"7", {\ar"6"},
"6", {\ar"10"},
"7", {\ar"11"},
"14", {\ar"7"},
"8", {\ar"9"},
"8", {\ar"12"},
"13", {\ar"8"},
"10", {\ar"9"},
"9", {\ar"13"},
"13", {\ar"10"},
"10", {\ar"14"},
"11", {\ar"14"},
"12", {\ar"13"},
"15", {\ar"13"},
"14", {\ar"15"},
\end{xy}
&
\begin{xy} 0;<0.38pt,0pt>:<0pt,-0.38pt>:: 
(26,50) *+{\point{0}} ="0",
(1,0) *+{\point{1}} ="1",
(113,0) *+{\point{2}} ="2",
(226,0) *+{\point{3}} ="3",
(0,114) *+{\point{4}} ="4",
(76,75) *+{\point{5}} ="5",
(151,75) *+{\point{6}} ="6",
(201,50) *+{\point{7}} ="7",
(26,175) *+{\point{8}} ="8",
(76,150) *+{\point{9}} ="9",
(151,150) *+{\point{10}} ="10",
(228,114) *+{\point{11}} ="11",
(1,225) *+{\point{12}} ="12",
(115,225) *+{\point{13}} ="13",
(226,225) *+{\point{14}} ="14",
(201,175) *+{\point{15}} ="15",
"0", {\ar"1"},
"4", {\ar"0"},
"5", {\ar"0"},
"0", {\ar"8"},
"1", {\ar"2"},
"2", {\ar"3"},
"2", {\ar"5"},
"6", {\ar"2"},
"3", {\ar"7"},
"8", {\ar"4"},
"5", {\ar"6"},
"9", {\ar"5"},
"7", {\ar"6"},
"6", {\ar"10"},
"7", {\ar"11"},
"15", {\ar"7"},
"8", {\ar"9"},
"12", {\ar"8"},
"10", {\ar"9"},
"9", {\ar"13"},
"13", {\ar"10"},
"10", {\ar"15"},
"11", {\ar"15"},
"13", {\ar"12"},
"14", {\ar"13"},
"15", {\ar"14"},
\end{xy}
&
\begin{xy} 0;<0.38pt,0pt>:<0pt,-0.38pt>:: 
(0,0) *+{\point{0}} ="0",
(25,50) *+{\point{1}} ="1",
(112.5,0) *+{\point{2}} ="2",
(225,0) *+{\point{3}} ="3",
(0,112.5) *+{\point{4}} ="4",
(75,75) *+{\point{5}} ="5",
(150,75) *+{\point{6}} ="6",
(200,50) *+{\point{7}} ="7",
(25,175) *+{\point{8}} ="8",
(75,150) *+{\point{9}} ="9",
(150,150) *+{\point{10}} ="10",
(225,112.5) *+{\point{11}} ="11",
(0,225) *+{\point{12}} ="12",
(112.5,225) *+{\point{13}} ="13",
(200,175) *+{\point{14}} ="14",
(225,225) *+{\point{15}} ="15",
"0", {\ar"1"},
"2", {\ar"0"},
"1", {\ar"2"},
"4", {\ar"1"},
"5", {\ar"1"},
"1", {\ar"8"},
"3", {\ar"2"},
"2", {\ar"5"},
"6", {\ar"2"},
"2", {\ar"7"},
"7", {\ar"3"},
"8", {\ar"4"},
"5", {\ar"6"},
"9", {\ar"5"},
"7", {\ar"6"},
"6", {\ar"10"},
"7", {\ar"11"},
"14", {\ar"7"},
"8", {\ar"9"},
"8", {\ar"12"},
"13", {\ar"8"},
"10", {\ar"9"},
"9", {\ar"13"},
"13", {\ar"10"},
"10", {\ar"14"},
"11", {\ar"14"},
"12", {\ar"13"},
"14", {\ar"13"},
"13", {\ar"15"},
"15", {\ar"14"},
\end{xy}
\\
\end{array}
\]
\[
\begin{array}{cccc}
\begin{xy} 0;<0pt,0.4114pt>:<0.7125pt,0pt>:: 
(180,120) *+{\point{0}} ="0",
(210,90) *+{\point{1}} ="1",
(165,75) *+{\point{2}} ="2",
(240,60) *+{\point{3}} ="3",
(60,120) *+{\point{4}} ="4",
(120,90) *+{\point{5}} ="5",
(165,45) *+{\point{6}} ="6",
(210,30) *+{\point{7}} ="7",
(30,90) *+{\point{8}} ="8",
(75,75) *+{\point{9}} ="9",
(120,30) *+{\point{10}} ="10",
(180,0) *+{\point{11}} ="11",
(0,60) *+{\point{12}} ="12",
(75,45) *+{\point{13}} ="13",
(30,30) *+{\point{14}} ="14",
(60,0) *+{\point{15}} ="15",
"1", {\ar"0"},
"4", {\ar"0"},
"0", {\ar"5"},
"2", {\ar"1"},
"1", {\ar"3"},
"3", {\ar"2"},
"5", {\ar"2"},
"2", {\ar"6"},
"6", {\ar"3"},
"3", {\ar"7"},
"5", {\ar"4"},
"4", {\ar"8"},
"9", {\ar"5"},
"7", {\ar"6"},
"6", {\ar"10"},
"11", {\ar"7"},
"8", {\ar"9"},
"12", {\ar"8"},
"9", {\ar"12"},
"13", {\ar"9"},
"10", {\ar"11"},
"10", {\ar"13"},
"15", {\ar"10"},
"11", {\ar"15"},
"12", {\ar"13"},
"14", {\ar"12"},
"13", {\ar"14"},
"14", {\ar"15"},
\end{xy}
&
\begin{xy} 0;<0pt,0.4114pt>:<0.7125pt,0pt>:: 
(180,120) *+{\point{0}} ="0",
(165,75) *+{\point{1}} ="1",
(210,90) *+{\point{2}} ="2",
(240,60) *+{\point{3}} ="3",
(60,120) *+{\point{4}} ="4",
(120,90) *+{\point{5}} ="5",
(165,45) *+{\point{6}} ="6",
(210,30) *+{\point{7}} ="7",
(30,90) *+{\point{8}} ="8",
(75,75) *+{\point{9}} ="9",
(120,30) *+{\point{10}} ="10",
(180,0) *+{\point{11}} ="11",
(0,60) *+{\point{12}} ="12",
(30,30) *+{\point{13}} ="13",
(75,45) *+{\point{14}} ="14",
(60,0) *+{\point{15}} ="15",
"0", {\ar"1"},
"2", {\ar"0"},
"0", {\ar"4"},
"5", {\ar"0"},
"1", {\ar"2"},
"1", {\ar"5"},
"6", {\ar"1"},
"2", {\ar"3"},
"3", {\ar"7"},
"4", {\ar"5"},
"4", {\ar"8"},
"9", {\ar"4"},
"5", {\ar"9"},
"7", {\ar"6"},
"10", {\ar"6"},
"6", {\ar"11"},
"11", {\ar"7"},
"8", {\ar"9"},
"12", {\ar"8"},
"9", {\ar"14"},
"11", {\ar"10"},
"14", {\ar"10"},
"10", {\ar"15"},
"15", {\ar"11"},
"13", {\ar"12"},
"14", {\ar"13"},
"13", {\ar"15"},
"15", {\ar"14"},
\end{xy}
&
\begin{xy} 0;<0pt,0.4114pt>:<0.7125pt,0pt>:: 
(180,120) *+{\point{0}} ="0",
(165,75) *+{\point{1}} ="1",
(210,90) *+{\point{2}} ="2",
(240,60) *+{\point{3}} ="3",
(60,120) *+{\point{4}} ="4",
(120,90) *+{\point{5}} ="5",
(165,45) *+{\point{6}} ="6",
(210,30) *+{\point{7}} ="7",
(30,90) *+{\point{8}} ="8",
(75,75) *+{\point{9}} ="9",
(120,30) *+{\point{10}} ="10",
(180,0) *+{\point{11}} ="11",
(0,60) *+{\point{12}} ="12",
(30,30) *+{\point{13}} ="13",
(75,45) *+{\point{14}} ="14",
(60,0) *+{\point{15}} ="15",
"0", {\ar"1"},
"2", {\ar"0"},
"0", {\ar"4"},
"5", {\ar"0"},
"1", {\ar"2"},
"1", {\ar"5"},
"6", {\ar"1"},
"3", {\ar"2"},
"2", {\ar"7"},
"7", {\ar"3"},
"4", {\ar"5"},
"4", {\ar"8"},
"9", {\ar"4"},
"5", {\ar"9"},
"7", {\ar"6"},
"10", {\ar"6"},
"6", {\ar"11"},
"11", {\ar"7"},
"8", {\ar"9"},
"8", {\ar"12"},
"13", {\ar"8"},
"9", {\ar"14"},
"11", {\ar"10"},
"14", {\ar"10"},
"10", {\ar"15"},
"15", {\ar"11"},
"12", {\ar"13"},
"14", {\ar"13"},
"13", {\ar"15"},
"15", {\ar"14"},
\end{xy}
&
\begin{xy}
<0pt,20pt>;<0.4275pt,20pt>:<0pt,20.4275pt>:: 
(50,0) *+{\point{0}} ="0",
(150,0) *+{\point{1}} ="1",
(100,50) *+{\point{2}} ="2",
(150,50) *+{\point{3}} ="3",
(0,50) *+{\point{4}} ="4",
(50,50) *+{\point{5}} ="5",
(150,100) *+{\point{6}} ="6",
(200,50) *+{\point{7}} ="7",
(0,150) *+{\point{8}} ="8",
(50,100) *+{\point{9}} ="9",
(150,150) *+{\point{10}} ="10",
(200,150) *+{\point{11}} ="11",
(50,150) *+{\point{12}} ="12",
(100,150) *+{\point{13}} ="13",
(50,200) *+{\point{14}} ="14",
(150,200) *+{\point{15}} ="15",
(0,-10) *+{ } ="16",
"1", {\ar"0"},
"4", {\ar"0"},
"0", {\ar"5"},
"3", {\ar"1"},
"1", {\ar"7"},
"2", {\ar"3"},
"5", {\ar"2"},
"3", {\ar"6"},
"7", {\ar"3"},
"5", {\ar"4"},
"4", {\ar"8"},
"9", {\ar"5"},
"6", {\ar"10"},
"11", {\ar"7"},
"8", {\ar"12"},
"14", {\ar"8"},
"12", {\ar"9"},
"10", {\ar"11"},
"10", {\ar"13"},
"15", {\ar"10"},
"11", {\ar"15"},
"13", {\ar"12"},
"12", {\ar"14"},
"14", {\ar"15"},
\end{xy}
\end{array}
\]
\end{figure}

\subsection{Examples from $n$-gons}\label{planar n}

The two families of selfinjective QPs given in section \ref{ex from ctilt} are planar. 
For each $n \geq 4$ define the following planar QP
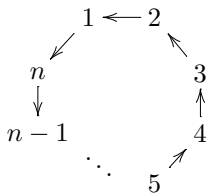
\[
\begin{xy} 0;<0.38pt,0pt>:<0pt,-0.38pt>:: 
(0,55) *+{n} ="0",
(50,0) *+{1} ="1",
(115,0) *+{2} ="2",
(160,55) *+{3} ="3",
(160,115) *+{4} ="4",
(115,165) *+{5} ="5",
(0,115) *+{n-1} ="6",
(60,140) *+{\ddots} ="7",
"1", {\ar"0"},
"0", {\ar"6"},
"2", {\ar"1"},
"3", {\ar"2"},
"4", {\ar"3"},
"5", {\ar"4"},
\end{xy}
\]
If $n$ is even, then define also
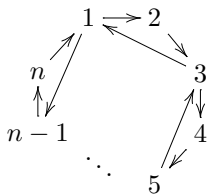
\[
\begin{xy} 0;<0.38pt,0pt>:<0pt,-0.38pt>:: 
(0,55) *+{n} ="0",
(50,0) *+{1} ="1",
(115,0) *+{2} ="2",
(160,55) *+{3} ="3",
(160,115) *+{4} ="4",
(115,165) *+{5} ="5",
(0,115) *+{n-1} ="6",
(60,140) *+{\ddots} ="7",
"0", {\ar"1"},
"6", {\ar"0"},
"1", {\ar"2"},
"3", {\ar"1"},
"1", {\ar"6"},
"2", {\ar"3"},
"3", {\ar"4"},
"5", {\ar"3"},
"4", {\ar"5"},
\end{xy}
\]
For each even $n \geq 4$ the above pair forms a single planar mutation class.

\section{Concluding remarks}
We started out by showing that the study of $2$-representation finite algebras can be reduced to selfinjective QPs and their cuts. 
In the Section \ref{planar} we saw that even in the quite restrictive setting of planar QPs many selfinjective ones appear. 
At present, all examples of selfinjective planar QPs that we are aware of can be obtained from triangles, square shaped ones and $n$-gons. 
Also, as far as we know they all have enough cuts. Thus we pose the following questions

\begin{question}
\begin{itemize}
\item[(1)] Is every selfinjective planar QP mutation-equivalent to one of the triangles in Section \ref{planar 3}, squares in Section \ref{planar 4} or $n$-gons in Section \ref{planar n}?
\item[(2)] Does every selfinjective planar QP have enough cuts?
\end{itemize}
\end{question}

\end{document}